\documentclass{article}

\usepackage{amsfonts,amssymb,amsthm}
\usepackage[totalwidth=17cm,totalheight=24cm]{geometry}

\catcode`\@=11\relax
\newwrite\@unused
\def\typeout#1{{\let\protect\string\immediate\write\@unused{#1}}}
\def\@nnil{\@nil}
\def\@empty{}
\def\@psdonoop#1\@@#2#3{}
\def\@psdo#1:=#2\do#3{\edef\@psdotmp{#2}\ifx\@psdotmp\@empty \else
    \expandafter\@psdoloop#2,\@nil,\@nil\@@#1{#3}\fi}
\def\@psdoloop#1,#2,#3\@@#4#5{\def#4{#1}\ifx #4\@nnil \else
       #5\def#4{#2}\ifx #4\@nnil \else#5\@ipsdoloop #3\@@#4{#5}\fi\fi}
\def\@ipsdoloop#1,#2\@@#3#4{\def#3{#1}\ifx #3\@nnil 
       \let\@nextwhile=\@psdonoop \else
      #4\relax\let\@nextwhile=\@ipsdoloop\fi\@nextwhile#2\@@#3{#4}}
\def\@tpsdo#1:=#2\do#3{\xdef\@psdotmp{#2}\ifx\@psdotmp\@empty \else
    \@tpsdoloop#2\@nil\@nil\@@#1{#3}\fi}
\def\@tpsdoloop#1#2\@@#3#4{\def#3{#1}\ifx #3\@nnil 
       \let\@nextwhile=\@psdonoop \else
      #4\relax\let\@nextwhile=\@tpsdoloop\fi\@nextwhile#2\@@#3{#4}}
\def\psdraft{
        \def\@psdraft{0}
}
\def\psfull{
        \def\@psdraft{100}
}
\psfull
\newif\if@prologfile
\newif\if@postlogfile
\newif\if@noisy
\def\pssilent{
        \@noisyfalse
}
\def\psnoisy{
        \@noisytrue
}
\psnoisy
\newif\if@bbllx
\newif\if@bblly
\newif\if@bburx
\newif\if@bbury
\newif\if@height
\newif\if@width
\newif\if@rheight
\newif\if@rwidth
\newif\if@clip
\newif\if@verbose
\def\@p@@sclip#1{\@cliptrue}
\def\@p@@sfile#1{
                   \def\@p@sfile{#1}
}
\def\@p@@sfigure#1{\def\@p@sfile{#1}}
\def\@p@@sbbllx#1{
                \@bbllxtrue
                \dimen100=#1
                \edef\@p@sbbllx{\number\dimen100}
}
\def\@p@@sbblly#1{
                \@bbllytrue
                \dimen100=#1
                \edef\@p@sbblly{\number\dimen100}
}
\def\@p@@sbburx#1{
                \@bburxtrue
                \dimen100=#1
                \edef\@p@sbburx{\number\dimen100}
}
\def\@p@@sbbury#1{
                \@bburytrue
                \dimen100=#1
                \edef\@p@sbbury{\number\dimen100}
}
\def\@p@@sheight#1{
                \@heighttrue
                \dimen100=#1
                \edef\@p@sheight{\number\dimen100}
}
\def\@p@@swidth#1{
                \@widthtrue
                \dimen100=#1
                \edef\@p@swidth{\number\dimen100}
}
\def\@p@@srheight#1{
                \@rheighttrue
                \dimen100=#1
                \edef\@p@srheight{\number\dimen100}
}
\def\@p@@srwidth#1{
                \@rwidthtrue
                \dimen100=#1
                \edef\@p@srwidth{\number\dimen100}
}
\def\@p@@ssilent#1{ 
                \@verbosefalse
}
\def\@p@@sprolog#1{\@prologfiletrue\def\@prologfileval{#1}}
\def\@p@@spostlog#1{\@postlogfiletrue\def\@postlogfileval{#1}}
\def\@cs@name#1{\csname #1\endcsname}
\def\@setparms#1=#2,{\@cs@name{@p@@s#1}{#2}}
\def\ps@init@parms{
                \@bbllxfalse \@bbllyfalse
                \@bburxfalse \@bburyfalse
                \@heightfalse \@widthfalse
                \@rheightfalse \@rwidthfalse
                \def\@p@sbbllx{}\def\@p@sbblly{}
                \def\@p@sbburx{}\def\@p@sbbury{}
                \def\@p@sheight{}\def\@p@swidth{}
                \def\@p@srheight{}\def\@p@srwidth{}
                \def\@p@sfile{}
                \def\@p@scost{10}
                \def\@sc{}
                \@prologfilefalse
                \@postlogfilefalse
                \@clipfalse
                \if@noisy
                        \@verbosetrue
                \else
                        \@verbosefalse
                \fi
}
\def\parse@ps@parms#1{
                \@psdo\@psfiga:=#1\do
                   {\expandafter\@setparms\@psfiga,}}
\newif\ifno@bb
\newif\ifnot@eof
\newread\ps@stream
\def\bb@missing{
        \if@verbose{
                \typeout{psfig: searching \@p@sfile \space  for bounding box}
        }\fi
        \openin\ps@stream=\@p@sfile
        \no@bbtrue
        \not@eoftrue
        \catcode`\%=12
        \loop
                \read\ps@stream to \line@in
                \global\toks200=\expandafter{\line@in}
                \ifeof\ps@stream \not@eoffalse \fi
                \@bbtest{\toks200}
                \if@bbmatch\not@eoffalse\expandafter\bb@cull\the\toks200\fi
        \ifnot@eof \repeat
        \catcode`\%=14
}       
\catcode`\%=12
\newif\if@bbmatch
\def\@bbtest#1{\expandafter\@a@\the#1
\long\def\@a@#1
\long\def\bb@cull#1 #2 #3 #4 #5 {
        \dimen100=#2 bp\edef\@p@sbbllx{\number\dimen100}
        \dimen100=#3 bp\edef\@p@sbblly{\number\dimen100}
        \dimen100=#4 bp\edef\@p@sbburx{\number\dimen100}
        \dimen100=#5 bp\edef\@p@sbbury{\number\dimen100}
        \no@bbfalse
}
\catcode`\%=14
\def\compute@bb{
                \no@bbfalse
                \if@bbllx \else \no@bbtrue \fi
                \if@bblly \else \no@bbtrue \fi
                \if@bburx \else \no@bbtrue \fi
                \if@bbury \else \no@bbtrue \fi
                \ifno@bb \bb@missing \fi
                \ifno@bb \typeout{FATAL ERROR: no bb supplied or found}
                        \no-bb-error
                \fi
                \count203=\@p@sbburx
                \count204=\@p@sbbury
                \advance\count203 by -\@p@sbbllx
                \advance\count204 by -\@p@sbblly
                \edef\@bbw{\number\count203}
                \edef\@bbh{\number\count204}
}
\def\in@hundreds#1#2#3{\count240=#2 \count241=#3
                     \count100=\count240        
                     \divide\count100 by \count241
                     \count101=\count100
                     \multiply\count101 by \count241
                     \advance\count240 by -\count101
                     \multiply\count240 by 10
                     \count101=\count240        
                     \divide\count101 by \count241
                     \count102=\count101
                     \multiply\count102 by \count241
                     \advance\count240 by -\count102
                     \multiply\count240 by 10
                     \count102=\count240        
                     \divide\count102 by \count241
                     \count200=#1\count205=0
                     \count201=\count200
                        \multiply\count201 by \count100
                        \advance\count205 by \count201
                     \count201=\count200
                        \divide\count201 by 10
                        \multiply\count201 by \count101
                        \advance\count205 by \count201
                     \count201=\count200
                        \divide\count201 by 100
                        \multiply\count201 by \count102
                        \advance\count205 by \count201
                     \edef\@result{\number\count205}
}
\def\compute@wfromh{
                \in@hundreds{\@p@sheight}{\@bbw}{\@bbh}
                \edef\@p@swidth{\@result}
}
\def\compute@hfromw{
                \in@hundreds{\@p@swidth}{\@bbh}{\@bbw}
                \edef\@p@sheight{\@result}
}
\def\compute@handw{
                \if@height 
                        \if@width
                        \else
                                \compute@wfromh
                        \fi
                \else 
                        \if@width
                                \compute@hfromw
                        \else
                                \edef\@p@sheight{\@bbh}
                                \edef\@p@swidth{\@bbw}
                        \fi
                \fi
}
\def\compute@resv{
                \if@rheight \else \edef\@p@srheight{\@p@sheight} \fi
                \if@rwidth \else \edef\@p@srwidth{\@p@swidth} \fi
}
\def\compute@sizes{
        \compute@bb
        \compute@handw
        \compute@resv
}
\def\psfig#1{\vbox {
        \ps@init@parms
        \parse@ps@parms{#1}
        \compute@sizes
        \ifnum\@p@scost<\@psdraft{
                \if@verbose{
                        \typeout{psfig: including \@p@sfile \space }
                }\fi
                \special{ps::[begin]    \@p@swidth \space \@p@sheight \space
                                \@p@sbbllx \space \@p@sbblly \space
                                \@p@sbburx \space \@p@sbbury \space
                                startTexFig \space }
                \if@clip{
                        \if@verbose{
                                \typeout{(clip)}
                        }\fi
                        \special{ps:: doclip \space }
                }\fi
                \if@prologfile
                    \special{ps: plotfile \@prologfileval \space } \fi
                \special{ps: plotfile \@p@sfile \space }
                \if@postlogfile
                    \special{ps: plotfile \@postlogfileval \space } \fi
                \special{ps::[end] endTexFig \space }
                \vbox to \@p@srheight true sp{
                        \hbox to \@p@srwidth true sp{
                                \hss
                        }
                \vss
                }
        }\else{
                \vbox to \@p@srheight true sp{
                \vss
                        \hbox to \@p@srwidth true sp{
                                \hss
                                \if@verbose{
                                        \@p@sfile
                                }\fi
                                \hss
                        }
                \vss
                }
        }\fi
}}
\catcode`\@=12\relax

\newcommand{\dr}{\partial}
\newcommand{\C}{{\mathbb C}}
\newcommand{\N}{{\mathbb N}}
\newcommand{\R}{{\mathbb R}}
\newcommand{\Z}{{\mathbb Z}}
\newcommand{\II}{I\hspace{-0.1cm}I}
\newcommand{\III}{I\hspace{-0.1cm}I\hspace{-0.1cm}I}
\newcommand{\tr}{\mbox{tr}}
\newcommand{\cat}{\mbox{CAT}}
\newcommand{\ric}{\mbox{ric}}
\newcommand{\dev}{\mbox{dev}}
\newcommand{\card}{\mbox{card}}
\newcommand{\sh}{\mbox{sh}}
\newcommand{\ch}{\mbox{CH}}
\newcommand{\arctg}{\mbox{arctg}}
\newcommand{\SO}{\mbox{SO}}
\newcommand{\CP}{\mbox{{\C}P}}
\newcommand{\can}{\mbox{can}}
\newcommand{\argth}{\mbox{argth}}
\newcommand{\hess}{\mbox{Hess}}
\newcommand{\inte}{\mbox{int}}

\newcommand{\ihm}{ideal hyperbolic manifold }

\newcommand{\deltab}{\overline{\delta}}
\newcommand{\gammab}{\overline{\gamma}}
\newcommand{\sigmab}{\overline{\sigma}}
\newcommand{\Sib}{\overline{\Sigma}}
\newcommand{\Thetab}{\overline{\Theta}}
\newcommand{\Omegab}{\overline{\Omega}}

\newcommand{\cb}{\overline{c}}
\newcommand{\eb}{\overline{e}}
\newcommand{\fb}{\overline{f}}
\newcommand{\gb}{\overline{g}}
\newcommand{\hb}{\overline{h}}
\newcommand{\pb}{\overline{p}}
\newcommand{\tb}{\overline{t}}
\newcommand{\vb}{\overline{v}}
\newcommand{\Bb}{\overline{B}}
\newcommand{\Cb}{\overline{C}}
\newcommand{\Gb}{\overline{G}}
\newcommand{\Kb}{\overline{K}}
\newcommand{\Pb}{\overline{P}}
\newcommand{\Mb}{\overline{M}}
\newcommand{\Sb}{\overline{S}}
\newcommand{\ricb}{\overline{\ric}}
\newcommand{\Db}{\overline{D}}
\newcommand{\Rb}{\overline{R}}
\newcommand{\thetab}{\overline{\theta}}
\newcommand{\omegab}{\overline{\omega}}

\newcommand{\Ct}{\tilde{C}}
\newcommand{\Ft}{\tilde{F}}
\newcommand{\Mt}{\tilde{M}}
\newcommand{\Nt}{\tilde{N}}
\newcommand{\Pt}{\tilde{P}}
\newcommand{\Qt}{\tilde{Q}}
\newcommand{\St}{\tilde{S}}
\newcommand{\dMt}{\tilde{\dr M}}
\newcommand{\phit}{\tilde{\phi}}
\newcommand{\gammat}{\tilde{\gamma}}
\newcommand{\sigmat}{\tilde{\sigma}}
\newcommand{\Sigmat}{\tilde{\Sigma}}

\newcommand{\Hsh}{{\mathcal H}_{\mbox{sh}}}
\newcommand{\Hex}{{\mathcal H}_{\mbox{ex}}}
\newcommand{\Hsm}{{\mathcal H}_{\mbox{sm}}}
\newcommand{\PMC}{\Phi_{\cM, \cC}}

\newcommand{\Thetasm}{\Theta_{\mbox{sm}}}

\newcommand{\bM}{{\bf M}}

\newcommand{\Rr}{\stackrel{\circ}{R}}

\newtheorem{prop}{Proposition}[section]
\newtheorem{lemma}[prop]{Lemma}
\newtheorem{sublemma}[prop]{Sub-lemma}

\newtheorem{thm}[prop]{Theorem}
\newtheorem{cor}[prop]{Corollary}
\newtheorem{remark}[prop]{Remark}

\newtheorem{df}[prop]{Definition}
\newtheorem{pty}[prop]{Property}
\newtheorem{question}[prop]{Question}

\newcommand{\pg}{\paragraph}

\newenvironment{thn}[1]{\vskip 0.2cm \noindent{\bf Theorem #1.} \it}{\rm
\vspace{0.2cm}} 
\newenvironment{crn}[1]{\vskip 0.2cm \noindent{\bf Corollary #1.} \it}{\rm
\vspace{0.2cm}} 
\newenvironment{lmn}[1]{\vskip 0.2cm \noindent{\bf Lemma #1.} \it}{\rm
\vspace{0.2cm}} 
\newenvironment{qn}[1]{\vskip 0.2cm \noindent{\bf Question #1.} \it}{\rm
\vspace{0.2cm}} 
\newenvironment{sketch}{\vskip 0.2cm \noindent{\bf Brief sketch of the
    proof~. ~~~}}{$\qed$ \vspace{0.2cm}} 

\newcommand{\btm}{\begin{thm}}
\newcommand{\etm}{\end{thm}}
\newcommand{\bpt}{\begin{pty}}
\newcommand{\ept}{\end{pty}}
\newcommand{\blm}{\begin{lemma}}
\newcommand{\elm}{\end{lemma}}
\newcommand{\bsl}{\begin{sublemma}}
\newcommand{\esl}{\end{sublemma}}
\newcommand{\bcr}{\begin{cor}}
\newcommand{\ecr}{\end{cor}}
\newcommand{\bdf}{\begin{df}}
\newcommand{\edf}{\end{df}}
\newcommand{\bprop}{\begin{prop}}
\newcommand{\eprop}{\end{prop}}
\newcommand{\bas}{\begin{asser}}
\newcommand{\eas}{\end{asser}}
\newcommand{\beq}{\begin{equation}}
\newcommand{\eeq}{\end{equation}}
\newcommand{\bpv}{\begin{proof}}
\newcommand{\epv}{\end{proof}}
\newcommand{\bpvs}{\begin{sketch}}
\newcommand{\epvs}{\end{sketch}}
\newcommand{\bit}{\begin{itemize}}
\newcommand{\eit}{\end{itemize}}
\newcommand{\bpn}{\begin{pfn}}
\newcommand{\epn}{\end{pfn}}
\newcommand{\btn}{\begin{thn}}
\newcommand{\etn}{\end{thn}}
\newcommand{\bcn}{\begin{crn}}
\newcommand{\ecn}{\end{crn}}
\newcommand{\bqn}{\begin{qn}}
\newcommand{\eqn}{\end{qn}}
\newcommand{\bln}{\begin{lmn}}
\newcommand{\eln}{\end{lmn}}
\newcommand{\brk}{\begin{remark}}
\newcommand{\erk}{\end{remark}}
\newcommand{\bq}{\begin{question}}
\newcommand{\eq}{\end{question}}

\newenvironment{pfn}[1]{\vskip 0.2cm \noindent{\it Proof #1.}}{$\square$
\vspace{0.2cm}}

\newcommand{\cA}{{\mathcal A}}
\newcommand{\cC}{\mathcal{C}}
\newcommand{\cD}{\mathcal{D}}
\newcommand{\cG}{\mathcal{G}}
\newcommand{\cH}{{\mathcal H}}
\newcommand{\cL}{{\mathcal L}}
\newcommand{\cM}{{\mathcal M}}
\newcommand{\cN}{{\mathcal N}}
\newcommand{\cP}{{\mathcal P}}
\newcommand{\cS}{{\mathcal S}}
\newcommand{\cT}{{\mathcal T}}

\newcommand{\cAb}{\overline{\mathcal A}}

\newcommand{\Met}{\mathcal{M}et}
\newcommand{\Imm}{\mathcal{I}mm}
\newcommand{\CMet}{\mathcal{CM}et}
\newcommand{\CImm}{\mathcal{CI}mm}
\newcommand{\gab}{\overline{\gamma}}
\newcommand{\hyp}{\mathbf{H}^3}
\newcommand{\dhyp}{\partial\hyp}
\newcommand{\isom}{\mathrm{Isom}}
\newcommand{\db}{\overline{\partial}}

\newcommand{\hbu}{\stackrel{\bullet}{h}}
\newcommand{\hbbu}{\stackrel{\bullet}{\hb}}
\newcommand{\Pd}{\stackrel{\bullet}{P}}
\newcommand{\Sd}{\stackrel{\bullet}{S}}
\newcommand{\Vb}{\stackrel{\bullet}{V}}
\newcommand{\alphad}{\stackrel{\bullet}{\alpha}}
\newcommand{\thetad}{\stackrel{\bullet}{\theta}}
\newcommand{\thetabbu}{\stackrel{\bullet}{\thetab}}

\begin{document}

\title{Hyperideal polyhedra in hyperbolic manifolds}

\author{Jean-Marc Schlenker\thanks{
Laboratoire Emile Picard, UMR CNRS 5580,
Universit{\'e} Paul Sabatier,
118 route de Narbonne,
31062 Toulouse Cedex 4,
France.
\texttt{schlenker@picard.ups-tlse.fr; http://picard.ups-tlse.fr/\~{
}schlenker}. }}

\date{December 2002; revised, June 2003}

\maketitle

\begin{abstract}

Let $(M, \dr M)$ be a 3-manifold with incompressible boundary that
admits a convex 
co-compact hyperbolic metric. We consider the
hyperbolic metrics on $M$ such that $\dr M$ looks locally like a
hyperideal polyhedron, and we characterize the possible dihedral
angles. 

We find as special cases the results of Bao and Bonahon
\cite{bao-bonahon} on hyperideal polyhedra, and those of Rousset
\cite{rousset1} on fuchsian hyperideal polyhedra. Our results can also
be stated in terms of circle configurations on $\dr M$, they provide an
extension of the Koebe theorem on circle packings.

The proof uses some elementary properties of the hyperbolic volume, in
particular the Schl{\"a}fli formula and the fact that the volume of
(truncated) hyperideal simplices is a concave function of the dihedral
angles. 

\bigskip

\begin{center} {\bf R{\'e}sum{\'e}} \end{center}

Soit $(M, \dr M)$ une vari{\'e}t{\'e} de dimension 3 {\`a} bord incompressible, qui
admet une 
m{\'e}trique hyperbolique convexe co-compacte. On consid{\`e}re les
m{\'e}triques hyperboliques sur $M$ pour lesquelles 
le bord ressemble localement {\`a} un poly{\`e}dre hyperbolique hyperid{\'e}al, et
on caract{\'e}rise les angles di{\`e}dres possibles. 

On retrouve comme cas particulier les r{\'e}sultats r{\'e}cents de Bao et
Bonahon \cite{bao-bonahon} pour les poly{\`e}dres hyperid{\'e}aux, et de Rousset
\cite{rousset1} 
pour les poly{\`e}dres hyperid{\'e}aux fuchsiens. Nos r{\'e}sultats peuvent aussi
s'exprimer en terme de configurations de cercles sur le bord de $M$, ils
donnent une extension du th{\'e}or{\`e}me de Koebe sur les empilements de
cercles. 

La preuve repose sur les propri{\'e}t{\'e}s {\'e}l{\'e}mentaires du volume hyperbolique,
en particulier sur la formule de Schl{\"a}fli et sur le fait que le volume
des simplexes hyperid{\'e}aux (tronqu{\'e}s) est une fonction concave des angles
di{\`e}dres.  

\end{abstract}

\tableofcontents

\pg{Hyperbolic manifolds with boundary}

In all this paper, we will consider a 3-manifold with boundary $(M, \dr
M)$. We suppose that it admits a complete, convex co-compact hyperbolic
metric. This is a topological assumption, which could be stated in
purely topological terms (see e.g. \cite{thurston-notes}). 

A natural question is to understand all the hyperbolic metrics on $M$ in
terms of quantities that can be read on the boundary. For complete,
convex co-compact metrics, this is achieved by the (hyperbolic version
of) the Ahlfors-Bers theorem \cite{ahlfors}, which states that those
metrics are uniquely determined by the conformal structure induced on
$\dr M$.

When we consider hyperbolic metrics such that $\dr M$ is smooth and
strictly convex 
(i.e. the boundary is at finite distance) there are some related
results. First, the induced metrics on the boundary are exactly the
metrics with curvature $K>-1$, and each is obtained in exactly one way
\cite{L4,hmcb}. In addition, the third fundamental forms of the boundar
(see section 1) are exactly the metrics with curvature $K<1$ and closed,
contractible geodesics of length $L>2\pi$ \cite{these,iie,hmcb}; this
means that $\tilde{\dr M}$ is globally $\cat(1)$.

\pg{Hyperideal boundaries}

It looks like those statements should not be restricted to the case
where the boundary is smooth; moreover, they should allow some situations
where the boundary of $M$ has some points at infinity. For instance, one
can consider "ideal hyperbolic manifolds" in the following sense.
First note that, given a hyperbolic metric $g$ with convex boundary on
$M$, there is a unique complete, convex co-compact hyperbolic manifold
$E(M)$ in which $(M, g)$ can be isometrically embedded in such a way
that the induced morphism $\pi_1M\rightarrow \pi_1 E(M)$ is an
isomorphism.
Then $(\tilde{E(M)}, g)=H^3$, and $\pi_1M$ has a natural action on $H^3$
by isometries. 

\bdf \label{df:ideal}
Let $g$ be a hyperbolic metric with convex boundary on $M$. We say that
$(M, g)$ is an {\bf ideal hyperbolic manifold} if:
\begin{itemize}
\item for each convex ball $\Omega\subset H^3$ and each isometric
  embedding $\phi:\Omega\rightarrow E(M)$ , the
intersection of $M$ 
with $\phi(\Omega)$ the image by $\phi$ of the intersection with
$\Omega$ of an ideal 
polyhedron $P\subset H^3$.
\item $\dr M$ contains no closed curve which is a geodesic of $M$.
\end{itemize}

\edf

This definition is a natural extension of the notion of ideal polyhedra
in $H^3$. As for ideal polyhedra, the third fundamental forms of those
manifolds is a measure located on the edges; understanding it is
equivalent to understanding the dihedral angles. Another rather simple example
is given by what can be called the "fuchsian" case, when $E(M)$ is the
quotient of $H^3$ by the $\pi_1$ of a closed surface which acts on $H^3$
fixing a totally geodesic plane $P_0$, and the universal cover of $M$, seen as
a subset in $H^3$, is invariant under the reflection in $P_0$. 

For hyperbolic
polyhedra, the description of the third fundamental form reduces to the
condition that the dual graph of the polyhedron is $\cat(1)$, see
\cite{rivin-annals}; more precisely, that its closed paths have length
$L\geq 2\pi$, with equality exactly when they bound a face. 
In the more general case of an ideal manifold, the
result is similar \cite{ideal}.

We will consider analogs of the ideal manifolds, but replacing the
notion of ideal polyhedron by the more general notion of hyperideal
polyhedron. A hyperideal polyhedron can be defined in at least two
equivalent ways.
\begin{itemize}
\item As the intersection of a finite set of half-spaces in $H^3$, with
  the condition that, for each end $E$, either all faces adjacent to $E$
  intersect in one ideal point, or there exists a plane which is
  orthogonal to all the faces adjacent to $E$.
\item Using the projective model of $H^3$ as the open unit ball $B^3$ in
  $\R^3$, the 
  hyperideal polyhedra are the intersections with $H^3$ of the (convex)
  polyhedra in $\R^3$ with all vertices outside $B^3$, but with all
  edges intersecting $B^3$.
\end{itemize}
Note that those two definitions allow for some ideal vertices,
i.e. vertices on the boundary at infinity of $H^3$. The other vertices
are called "strictly hyperideal", and a hyperideal polyhedron with no
ideal vertex is called "strictly hyperideal". 

The same definition as for ideal manifolds can be used to define
"hyperideal manifolds", 
i.e. hyperbolic manifold with a boundary that looks locally like a
hyperideal hyperbolic polyhedron. 

\bdf \label{df:hyperideal}
Let $g$ be a hyperbolic metric with convex boundary on $M$. We say that
$(M, g)$ is a {\bf hyperideal hyperbolic manifold} if:
\begin{itemize}
\item for each convex ball $\Omega\subset H^3$ and each isometric
  embedding $\phi:\Omega\rightarrow E(M)$ , the
intersection of $M$ 
with $\phi(\Omega)$ is the image by $\phi$ of the intersection with
$\Omega$ of a hyperideal 
polyhedron $P\subset H^3$.
\item $\dr M$ contains no closed curve which is a geodesic of $M$.
\end{itemize}
\edf

This definition, and in particular the second point, is designed to
exclude some "bad" situations where $\dr M$ has non-empty intersection
with the convex core of $M$. More details on this can be found in
section 7. Given a 3-manifold $M$ with boundary, a "hyperideal metric"
on $M$ is a hyperbolic metric such that $(M, g)$ is a hyperideal
hyperbolic manifold. We will sometimes call this a "hyperideal
hyperbolic structure" on $M$.

The main goal of this paper is to understand the possible dihedral
angles of hyperideal manifolds, and to obtain a result similar to the
result obtained for hyperideal polyhedra in \cite{bao-bonahon}.

\pg{Dihedral angles}

Before stating the main results, we have to define some sequences of
edges which play a special role. We consider now a cellulation $\sigma$
of $\dr M$, i.e. a decomposition of $\dr M$ in the union of a finite
number of embedded images of the interior of polygons in $\R^2$, with
disjoint interior, such that the intersection of two adjacent polygons
is an edge of each. We will call $\sigma_1$ the 1-skeleton of $\sigma$,
which is a graph, and $\sigma_1^*$ the dual graph.

\bdf \label{df:circuit}
A {\bf circuit} in $\sigma$ is a sequence $e_0, e_1, \cdots, e_n=e_0$
which corresponds to the successive edges of a closed path in
$\sigma_1^*$ which is homotopically trivial in $M$. A circuit is
{\bf elementary} if the dual closed path in $\sigma_1^*$ bounds a face.
\edf

\bdf \label{df:path}
A {\bf simple path} is a sequence of edges $e_1, \cdots, e_n$ in $\sigma_1$
corresponding to the successive edges of a path in $\sigma_1^*$, which:
\begin{itemize}
\item begins and ends at boundary points of a face $f$ of $\sigma_1^*$.
\item is not included in the boundary of $f$.
\item is homotopic in $M$ to a segment in $f$. 
\end{itemize}
\edf

We can now state our main result. 

\btm \label{tm:angles}
Suppose that $M$ has incompressible boundary. 
Let $\sigma$ be a cellulation of $\dr M$, and let $w:\sigma_1\rightarrow
(0,\pi)$ be a map on the set of edges of $\sigma$. There exists a
hyperideal hyperbolic structure on $M$, with boundary combinatorics
given by $\sigma$ and exterior dihedral angles given by $w$, if and only
if:
\begin{itemize}
\item the sum of the values of $w$ on each circuit in $\sigma_1$ is
greater than $2\pi$, and strictly greater if the circuit is
non-elementary.
\item The sum of the values of $w$ on each simple path in $\sigma_1$ is
strictly larger than $\pi$.
\end{itemize}
This hyperideal structure is then unique. 
\etm

This is an extension of the main result of \cite{ideal}, which concerns
ideal manifolds only.

\pg{Outline of the proof}

The proof is related to the method used in \cite{ideal}; the
starting point is the Schl{\"a}fli formula, which describes the first-order
variations of the volume of a polyhedron in terms of the variation of its
dihedral angles (see section 1). A simple consequence, obtained in
section 3, is that the volume of
hyperideal simplices (the definition is below) is a strictly concave
function of the dihedral angles. This fact was well known to be true for
ideal simplices, and this is the basis for several important constructions
concerning ideal polyhedra (see e.g. \cite{thurston-notes}, chapter 7,
\cite{rivin-annals,Ri2}). 

\begin{lmn}{\ref{lm:concave-1}}
For each $i\in \{0,\cdots, 4\}$, 
the volume $V$ is a strictly concave function on the space of hyperideal
simplices
having exactly $i$ ideal vertices $v_1, \cdots, v_i$, parametrized by
the dihedral angles. 
\end{lmn}

A consequence of the concavity of the volume is obtained using an
interesting technique, based on deformations among singular hyperbolic
structures to get hyperbolic metrics; those ideas can be traced back to
the work of Thurston \cite{thurston-notes} on the Andreev theorem
\cite{Andreev-ideal}, and 
then of Colin de Verdi{\`e}re \cite{CdeV}, Br{\"a}gger \cite{bragger}, and Rivin
\cite{Ri2}. Applying those ideas to hyperideal polyhedra, one 
obtains the following result, which we will prove in section 4 since
this proof is partly different from the one given by Bao and Bonahon. It
is a special case of theorem \ref{tm:angles}, but also a tool in its proof.

\begin{lmn}{\ref{lm:angles-poly} (Bao, Bonahon \cite{bao-bonahon})}
Let $\sigma$ be a cellulation of $S^2$, and let $w:\sigma_1\rightarrow
(0,\pi)$ be a map on the set of edges of $\sigma$. There exists a
hyperideal polyhedron with combinatorics
given by $\sigma$ and exterior dihedral angles given by $w$ if and only
if:
\begin{itemize}
\item the sum of the values of $w$ on each circuit in $\sigma_1$ is
greater than $2\pi$, and strictly greater if the circuit is
non-elementary.
\item The sum of the values of $w$ on each simple path in $\sigma_1$ is
  strictly larger than $\pi$.
\end{itemize}
This hyperideal polyhedron is then unique. 
\end{lmn}

Since the sum of a finite number of concave functions is concave, lemma
\ref{lm:concave-1} 
can be used to prove that the volume of any hyperideal
polyhedron is also a concave function of the dihedral angles, and this
actually also applies to hyperideal manifolds. 

\begin{lmn}{\ref{lm:poly-vol}}
Let $\sigma$ be a cellulation of $S^2$. The volume is a strictly concave
function of the dihedral angles, on the space of hyperideal polyhedra
with combinatorics given by $\sigma$.
\end{lmn}

We will define in section 4 hyperideal cellulations  as decompositions
of a hyperideal manifolds in isometric 
images of hyperideal polyhedra, with some non-degeneracy
conditions. Section 4 contains the proof of the: 

\begin{lmn}{\ref{lm:triang}}
Any hyperideal manifold admits a hyperideal cellulation.
\end{lmn}

Using this, the concavity of the volume of polyhedra and the Schl{\"a}fli
formula, we will obtain in section 5 an infinitesimal rigidity
statement. 

\begin{lmn}{\ref{lm:rigidity}}
Let $M$ be a hyperideal manifold. Any first-order deformation of its
dihedral angles is obtained by a unique first-order deformation of $M$. 
\end{lmn}

This can be considered as the key lemma of this paper. The infinitesimal
rigidity of hyperbolic manifolds with convex boundary is a problem that
can be tackled in different ways. For instance, for polyhedra, it is proved in
\cite{bao-bonahon} by the Cauchy-Legendre method; in \cite{hmcb}, the
result is obtained 
using some analytic estimates and transformations defined by Pogorelov to
translate rigidity problems from $H^3$ to $\R^3$. The method used here is
completely different. 

Another important element, as in many deformation proofs, is to obtain a
compactness result. This is done in section 6, where the following
result is proved. 

\begin{lmn}{\ref{lm:compact}}
Let $(g_n)_{n\in \N}$ be a sequence of hyperideal structures on $M$,
with the same  boundary combinatorics. For each $n$, let $\alpha_n$ be
the function which associates to each boundary edge of $(M, g_n)$ its
exterior dihedral angle, and suppose that $\alpha_n\rightarrow \alpha$,
where $\alpha$ still satisfies the hypothesis of theorem
\ref{tm:angles}. Then, after taking a subsequence, $g_n$ converges to
a hyperideal structure $g$ on $M$.
\end{lmn}

Finally, a techically important point is to prove that the space of
dihedral angle assignations appearing in the hypothesis of theorem
\ref{tm:angles} is connected. The key point is that, when $M$ has
incompressible boundary, the conditions on
the dihedral angles on the various boundary components of $\dr M$ behave
independently, so that it is sufficient to prove the connectedness of
the space of angles in the "fuchsian" case, i.e. when one considers a
manifold which is topologically $S\times \R$, where $S$ is a compact
surface of genus at least $2$, with an isometric involution fixing a
compact surface. But this case is well understood thanks to a result of
Rousset, who proved in this case the analog of theorem \ref{tm:angles},
although by very different methods. 

\begin{thn}{\ref{tm:fuchsian} (M. Rousset \cite{rousset1})}
Let $S$ be a compact surface of genus at least 2, let $\sigma$ be a
cellulation of $S$, and let $w:\sigma_1\rightarrow 
(0,\pi)$ be a map on the set of edges of $\sigma$. There exists a
hyperideal fuchsian realization of $S$, with boundary combinatorics
given by $\sigma$ and exterior dihedral angles given by $w$, if and only
if:
\begin{itemize}
\item The sum of the values of $w$ on each circuit in $\sigma_1$ is
greater than $2\pi$, and strictly greater if the circuit is
non-elementary.
\item The sum of the values of $w$ on each simple path in $\sigma_1$ is
strictly larger than $\pi$.
\end{itemize}
This hyperideal realization is then unique. 
\end{thn}

\pg{The Koebe circle packing theorem}

The hyperideal polyhedra are related to the Koebe circle packing
theorem. Recall that a {\bf circle packing} in $S^2$ is a finite set of
circles bounding disjoint open disks, and the {\bf incidence graph} of a
circle packing is the (combinatorially defined) graph with one vertex
for each circle, and an edge between two circles if and only if the
circles are tangent.

\btm[Koebe \cite{koebe}] 
Let $\Gamma$ be the 1-skeleton of a triangulation of $S^2$. There is a
unique circle packing in $S^2$ with incidence graph $\Gamma$.
\etm

The uniqueness here is up to M\"obius transformations. 
This theorem has an extension to graphs which are the 1-skeleton of a
cellulation of $S^2$, but the uniqueness demands some additional
hypothesis. Given a circle packing, an {\bf interstice} is a connected
component of the complement of the (closed) disks bounded by the circles.

\btm[Koebe \cite{koebe}] \label{tm:koebe2}
Let $\Gamma$ be the 1-skeleton of a cellulation of $S^2$. There is a
unique circle packing $C$ in $S^2$ such that:
\begin{itemize}
\item the incidence graph of $C$ is $\Gamma$.
\item for each interstice $I$, there is a circle orthogonal to each
  circle of $C$ adjacent to $I$. 
\end{itemize}
\etm

Thurston \cite{thurston-notes} realized that the Koebe circle packing is
related to ideal polyhedra and the Andreev theorem. There is however a
simpler relationship between circle packings and hyperideal polyhedra. 
Let $\sigma$ be a cellulation of $S^2$. Consider the function
$w:\sigma_1\rightarrow (0,\pi)$ defined by $w(e)=\pi-\epsilon$ for each
edge $e$ of $\sigma$. By lemma \ref{lm:angles-poly}, if $\epsilon$ is
small enough, there is a unique hyperideal polyhedron $P_w$ with
combinatorics given by $\sigma$ and exterior dihedral angles by
$w$. 

The intersections of the faces of $P_w$ (or more precisely of the planes
containing the faces) with $\dr_\infty H^3$ are circles; when two circles
correspond to faces with a common edge $e$, they intersect
with angles $w(e)$. Moreover, the truncated faces (see sections 1 and 3)
correspond to another family of circles, which intersect the circles of
the first family orthogonally.

Taking the limit as $w\rightarrow \pi$ on each edge, we find two
families $F_1, F_2$ of circles on $S^2$, such that:
\begin{itemize}
\item the circles of $F_1$ correspond to the faces of
  $\sigma$. Those circles intersect if and only if the corresponding
  faces of $\sigma$ share an edge, and they are then tangent.
\item the circles of $F_2$ correspond to the vertices of $\sigma$. They
  intersect if and only if the corresponding vertices are adjacent, and
  they are then tangent.
\item a circle $c_1$ of $F_1$ intersects a circle $c_2$ of $F_2$ if and
  only if the face corresponding to $c_1$ contains the vertex
  corresponding to $c_2$. The intersection is then orthogonal.
\end{itemize}
This is another description of theorem \ref{tm:koebe2}.
Note that, by the Schl{\"a}fli formula (see equations
(\ref{eq:schlafli}) in section 3), the limit taken here corresponds to
letting the volume of the polyhedra go to its maximal value.
The same line of reasoning leads from theorem \ref{tm:angles} to the
following extension of the Koebe circle packing theorem. We first state
the simpler form where we only consider triangulations. 

\btm 
Suppose that $M$ has incompressible boundary. Let $\Gamma$ be the
1-skeleton of a triangulation of $\dr M$. There exists a unique couple
$(c,C)$, where $c$ is a $\C P^1$-structure on $\dr M$ induced by a
complete, convex co-compact hyperbolic metric on $M$, and $C$ is a
circle packing of $(\dr M, c)$ with incidence graph $\Gamma$.
\etm

\btm \label{tm:koebe}
Suppose that $M$ has incompressible boundary.
Let $\sigma$ be a cellulation of $\dr M$. There exists a unique triple
$(c, F_1, F_2)$, where: 
\begin{itemize}
\item $c$ is a $\C P^1$-structure on $\dr M$ induced by a
complete, convex co-compact hyperbolic metric on $M$.
\item $F_1$ and $F_2$ are circle packings of $(\dr M, c)$.
\item the circles of $F_1$ correspond to the faces of
  $\sigma$. Those circles intersect if and only if the corresponding
  faces of $\sigma$ share an edge, and they are then tangent.
\item the circles of $F_2$ correspond to the vertices of $\sigma$. They
  intersect if and only if the corresponding vertices are adjacent, and
  they are then tangent.
\item a circle $c_1$ of $F_1$ intersects a circle $c_2$ of $F_2$ if and
  only if the face corresponding to $c_1$ contains the vertex
  corresponding to $c_2$. The intersection is then orthogonal.
\end{itemize}
\etm

This result contains as a special case, obtained by
considering the "fuchsian" situation, some known results on circle
packings on surfaces of genus at least $2$ with a hyperbolic metric (see
\cite{thurston-notes,CdeV}). A slightly more general
example is obtained by considering a closed surface 
$S$ of genus $g\geq 2$, and the 3-manifold $M:=S\times \R$. Then $M$
has two boundary components, which we call $S_1$ and $S_2$, each
diffeomorphic to $S$. Let $\Gamma_1$, and 
$\Gamma_2$ be the 1-skeletons of finite triangulations, in $S_1$ and $S_2$
respectively. There is then a unique quasi-fuchsian hyperbolic metric on
$M$, inducing $\C P^1$-structures $c_1$, $c_2$ on $S_1$ and $S_2$,
a unique circle 
packing $C_1$ in $S_1$ for $c_1$ and a unique circle packing $C_2$ in
$S_2$ for $c_2$, such that the incidence graph of $C_1$ is $\Gamma_1$,
the incidence graph of $C_2$ is $\Gamma_2$.

Theorem \ref{tm:koebe} is actually a limit case of a more general
statement on configurations of circles, theorem \ref{tm:circles}, 
which is a direct translation of
theorem \ref{tm:angles}. The configurations appearing there have two
families of circles, and are more general than those usually associated
to (generalizations of) the Andreev theorem on ideal polyhedra. To each
such circle configuration one associates a volume, and theorem
\ref{tm:koebe} is obtained as the limit case when the volume is
maximal. The proof of theorem \ref{tm:koebe}, along with some additional
details, can be found in section 9.

\pg{Induced metrics}

Once we know that the hyperideal manifolds, with given boundary
combinatorics, are parametrized by their dihedral angles, the fact that
the volume has non-degenerate hessian can be translated, using the
Schl{\"a}fli formula, into an {\bf infinitesimal rigidity} statement: any
non-zero first-order deformation induces a non-trivial deformation of
the metric induced on the boundary. 

\begin{lmn}{\ref{lm:rig-inf}}
Let $(M, g)$ be a hyperideal manifold. It has no first-order
deformation (among hyperideal manifolds) which does not change the
induced metric on $\dr M$.
\end{lmn}

Using this lemma, we can recover rather simply a result describing the
induced metrics on hyperideal polyhedra (see \cite{shu}).

\begin{thn}{\ref{tm:poly-metrics}}
Let $h$ be a complete hyperbolic metric on $S^2$ minus a finite number
of points. There is a unique hyperideal polyhedron on $H^3$ whose
induced metric is $h$.
\end{thn}

The same arguments also yield an analogous result for the fuchsian case.

\begin{thn}{\ref{tm:fuchsian-metrics}}
Let $S$ be a compact surface with non-empty boundary of genus at least
$2$, and let $h$ be a complete hyperbolic metric on $S$ minus a finite
number of points. There is a
unique fuchsian hyperideal manifold $(M, g)$ such that the induced
metric on both connected components of $\dr M$ is $h$.  
\end{thn}

\section{Hyperideal manifolds}

We recall here some basic facts about hyperbolic geometry, in particular
hyperideal polyhedra.

\pg{Hyperbolic 3-space and the de Sitter space}

Hyperbolic 3-space can be constructed as a quadric in the Minkowski
4-space $\R^4_1$, with the induced metric:
$$ H^3 := \{ x\in \R^4_1 ~|~ \langle x,x \rangle =-1 ~ \wedge x_0>0
\}~. $$
But $\R^4_1$ also contains another quadric, the de
Sitter space of dimension 3:
$$ S^3_1  := \{ x\in \R^4_1 ~|~ \langle x,x \rangle =1\}~. $$
By construction, it is Lorentzian and has an action of $\SO(3,1)$ which
is transitive on orthonormal frames, so it has constant curvature; one
can easily check that its curvature is $1$. $S^3_1$ contains many
space-like totally geodesic 2-planes, each isometric to $S^2$ with its
canonical round metric. Each separates $S^3_1$ into two ``hemispheres'',
each isometric to a model which we will denote by $S^3_{1,+}$. For instance,
in the quadric above, 
the set of points $x\in S^3_1$ with positive first coordinate $x_0>0$ is a
hemisphere.

\pg{The third fundamental form and the dual metric}

Consider a smooth surface $S\subset H^3$. The Riemannian metric on $H^3$
defines by restriction a Riemannian metric on $S$, which is called the
induced metric, or first fundamental form, of $S$. We will denote it by
$I$. 

There is another metric, the third fundamental form, which is defined on
a smooth, strictly convex surface $S$ in $H^3$. To define it, let $N$ be
a unit normal vector field to $S$, and let $D$ be the Levi-Civit{\`a}
connection of $H^3$; the second fundamental form of $S$ is defined by:
$$ \forall s\in S, ~\forall x,y\in T_sS, ~\II(x,y)=-I(D_xN,
y)=-I(x,D_yN)~, $$ 
and the third fundamental form by:
$$ \forall s\in S, ~\forall x,y\in T_sS, ~\III(x,y)=I(D_xN, D_yN)~. $$
The same definition applies in $S^3_1$. 

There is a polyhedral analog of the third fundamental form. For a compact
polyhedron $P\subset H^3$, it can be defined by gluing, for each vertex $v$
of $P$, the interior of a spherical polygon which is the dual of the
link of $P$ 
at $v$. The duality which is used here is the projective duality in
$S^2$, so that the polygon has an edge of length $\alpha$ for each edge
adjacent to $v$ with exterior dihedral angle $\alpha$. This third
fundamental form is also often called the "dual metric" of the
polyhedron. The definition of
the dual metric for ideal or hyperideal polyhedra is outlined below
using the hyperbolic-de Sitter duality.

\pg{The hyperbolic-de Sitter duality}

The understanding of the third fundamental form of surfaces in $H^3$ relies
heavily on an important duality between $H^3$ and the de Sitter space
$S^3_1$. It associates to each point $x\in H^3$ a space-like, totally
geodesic plane in $S^3_1$, and to each point $y\in S^3_1$ an oriented
totally geodesic plane in $H^3$. 

It can be defined using the quadric models of $H^3$ and $S^3_1$. Let
$x\in H^3$; define $d_x$ as the line in $\R^4_1$ going through $0$ and
$x$. $d_x$ is a time-like line; call $d_x^*$ the orthogonal space in
$\R^4_1$, which is a space-like 3-plane. So $d_x^*$ intersects $S^3_1$
in a space-like totally geodesic 2-plane, which we call $x^*$, and which
is the dual of $x$. Conversely, given a space-like totally geodesic
plane $p\in S^3_1$, it is the intersection with $S^3_1$ of a space-like
3-plane $P\ni 0$. Let $d$ be its orthogonal, which is a time-like line;
the dual $p^*$ of $p$ is the intersection $d\cap H^3$.

The same construction works in the opposite direction. Given a point
$y\in S^3_1$, we call $d_y$ the oriented line going through $0$ and $y$,
and $d_y^*$ its orthogonal, which is an oriented time-like 3-plane. Then
$y^*:=d_y^*\cap H^3$ is an oriented totally geodesic plane. 

We can then define the duality on surfaces. Given a smooth, oriented
surface $S\subset H^3$, its dual $S^*$ is the set of points in $S^3_1$
which are the dual of the oriented planes which are tangent to
$S$. If $S$ is smooth and strictly convex, then $S^*$ is smooth,
space-like, and strictly convex. Conversely, given a smooth, space-like
surface $\Sigma\subset 
S^3_1$, its dual is the set $\Sigma^*$ of points in $H^3$ which are the
duals of the planes tangent to $\Sigma$. 

One of the main properties of this duality, on smooth surfaces, is that
the induced metric on $\Sigma^*$ is the third fundamental form of
$\Sigma$, while the third fundamental form of $\Sigma^*$ in $S^3_1$ is
the induced metric on $\Sigma$.

\pg{The dual metric of hyperbolic polyhedra}

The hyperbolic-de Sitter duality works also for compact polyhedra in
$H^3$ (for which it was introduced in \cite{Ri,RH}). 
Given a compact polyhedron $P\subset H^3$, its dual is the
convex, space-like polyhedron $P^*\subset S^3_1$ with:
\begin{itemize}
\item as vertices, the duals of the planes containing the faces of $P$;
\item as edges, segments of the geodesics dual to the geodesics
  containing the edges of $P$;
\item as faces, polygons in the planes dual to the vertices of $P$.
\end{itemize}
$P^*$ can be defined as the set of points dual to the support planes of
$P$.
The main point is that the induced metric on $P^*$ is the dual metric of
$P$. It is a spherical cone-metric, with singular points corresponding to
the vertices of $P^*$, where the singular curvature is negative. 

When $P$ has ideal vertices, the set of points duals to the support
planes of $P$ has a "hole" for each ideal vertex $v$ of $P$. This hole
corresponds, in the projective model of $H^3$, to the face of the dual
of $P$ (defined in projective terms) tangent to $\dr D^3$ at $v$. The
length $l$ of the boundary of this face is equal to the sum of the exterior
angles of the link of $P$ at $v$, which is Euclidean; thus $l=2\pi$. For
instance, if $P$ is an ideal polyhedron, the set of dual points of its
support planes is a graph, which we will call its {\bf dual graph}.
It defines a cellulation of $S^2$, with each face of length $2\pi$. 
Figure 1 represents, in the projective model of $H^2$, 
the dual of an ideal polygon, which is a finite
set of points; it should help understand what the dual of an ideal
polyhedron in $H^3$ is. 

\begin{figure}[h]
\centerline{\psfig{figure=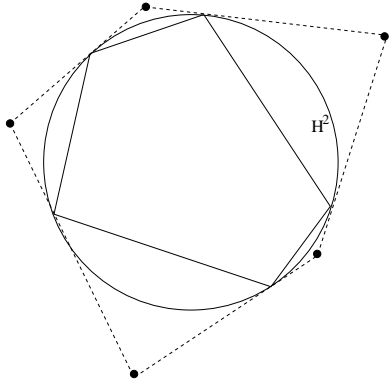,height=5cm}}
\caption{The dual of an ideal polygon}
\end{figure}

When $P$ has some ideal points, we will define its dual metric $\III$ as
the metric obtained by gluing in the "holes"  corresponding to the ideal
points a hemisphere (i.e. isometric to a hemisphere of $S^2$ with its
canonical round metric). 
The result is a metric space
which obviously has negative singular curvature at its
singular points, because the singular points correspond to the vertices
of the graph, and the total angle around those points is $\pi$ times
the number of faces. $\III$ is the ``natural'' third fundamental form
for instance in a limit sense, as follows: 

\bpt
Let $(\Omega_n)_{n\in \N}$ be an increasing sequence of open subsets of
the interior of $P$ with smooth, convex boundary, such that $\cup_n
\Omega_n$ is the interior of $P$. Then 
the third fundamental forms of $\dr \Omega_n$ converge to the dual metric of
$P$. 
\ept

We leave the proof to the reader.

For a hyperideal polyhedron $P$, one can still consider the set of support
planes, which is a graph in the de Sitter space. The total length
of the edges of a face, however, is strictly greater than $2\pi$ for
faces corresponding to strictly hyperideal vertices. To define the dual
metric, one has to glue in each of the faces a "singular hemisphere",
obtained as the quotient by a rotation of angle $\theta>2\pi$ of the
universal cover of the complement of the center in a hemisphere. 

\pg{Truncated hyperideal polyhedra}
Given a hyperideal polyhedron
$P\subset H^3$, there is, for each strictly hyperideal $v$ vertex of $P$, a
plane $p$ which is 
orthogonal to all the faces of $P$ adjacent to $v$. $p$ is the dual of $v$ in
the hyperbolic-de Sitter duality. Following \cite{bao-bonahon}, we call  
$P_t$ the
polyhedron obtained by cutting all the ends of $P$ by this plane. $P_t$
is a polyhedron of finite volume in $H^3$, and is compact if and only if
$P$ is strictly hyperideal. The dual metric of $P$ can also be defined
as the dual metric of $P_t$.

The hyperbolic polyhedra which can be obtained as truncated hyperideal
polyhedra are quite special. The set of their faces can be split into
two subset, the "real" faces and the "cuts". Each "cut" is adjacent to
"real faces" only, while each non-ideal vertex is adjacent to exactly one
"cut". Moreover, the dihedral angle between the "cuts" and the "real
faces" is always $\pi/2$. 

The dual metric of a hyperideal polyhedron 
is therefore quite special too. It is
constructed by gluing pieces which are "singular hemispheres", i.e. the
quotient of the universal cover of the complement of the center in a
spherical hemisphere by a rotation of angle strictly larger than
$2\pi$. Actually the dual metric of $P_t$ is the same as the dual metric
of $P$, defined by gluing singular hemispheres in the faces of the dual
graph. More details about this can be found in \cite{rousset1}.

\pg{The dual metrics of hyperideal manifolds}
The previous considerations on duality were local, so they extend from
polyhedra to the boundaries of ideal or hyperideal manifolds. Here we
consider a hyperideal manifold $M$. We can define the dual graph of its
universal cover as the set of points in $S^3_1$ dual to the support
planes of $\Mt\subset H^3$, and the dual graph of $M$ as the quotient of
the dual graph of $\Mt$ by the action of $\pi_1M$ on $H^3$ and $S^3_1$. 

The dual metric --- which we will still call the third fundamental form of
the boundary --- can still be defined by gluing in each face of the dual
graph a singular hemisphere. It
lifts to a $\cat(1)$ metric on the boundary of the universal
cover of $M$. Indeed, this splits into a local curvature condition ---
which is satisfied by the local convexity, because the singular
curvature is negative at each vertex --- and a global condition on the
length of the closed geodesics, which is also true here because of lemma
\ref{lm:necessary} below. 

This dual metric can also be defined as the dual metric of the
hyperbolic manifold with convex, polyhedral boundary obtained by
truncating the strictly hyperideal ends of a hyperideal manifold.

\pg{Necessary conditions on the angles}
The third fundamental form defined in this way has the important
properties below. 

\blm \label{lm:necessary}
Let $M$ be a hyperideal manifold. Its dual metric 
$\III$ lifts to a $\cat(1)$ metric on the boundary of the universal
cover of $M$.
\elm

Consider the hyperbolic manifold with boundary $N$ obtained by
truncating the strictly hyperideal ends of $M$. 
The universal cover $\Nt$ of $N$ has a canonical embedding into $H^3$,
because $N$ carries a hyperbolic manifold with convex boundary. Consider
its boundary $\dr\Nt$, it is a convex surface in $H^3$. Therefore, its
dual metric is $\cat(1)$. Indeed:
\begin{itemize}
\item it can be considered as the induced metric on the dual surface in
  $S^3_1$, and then the non-smooth version of the Gauss formula shows
  that it is locally $\cat(1)$.
\item its closed geodesics have length $L>2\pi$, as shown by various
  geometric arguments, see e.g. \cite{CD,RH,these,shu,cpt}.
\end{itemize}

As a consequence of the description of the dual metrics of the
hyperideal manifolds, it is clear that the simple paths and the circuits
which appear in theorem \ref{tm:angles} correspond to geodesics segments,
resp. to closed geodesics, of the dual
metric. Therefore, since the dual metrics are $\cat(1)$, it is already
clear that the conditions on the lengths of the circuits and of the
simple paths in theorem \ref{tm:angles} are necessary.

\pg{The dual metric and the dihedral angles}
A rather important point is that, while knowing the dual metric (and the
combinatorics of the polyhedral surface) determines the dihedral
angles, the converse is not true --- the dihedral angles
determine the length of the edges of the dual surface, but it does not
determine the shape of the dual faces with more than 4 edges. This is
already the source of interesting questions for convex polyhedra in
$H^3$; for instance, it is still an open problem whether a
convex hyperbolic polyhedron can be infinitesimally deformed without
changing its 
dihedral angles, see e.g. \cite{dap}, or \cite{stoker} for an analogous
problem in the Euclidean case. 


\section{Hyperideal triangles}

\pg{The projective model of $H^2$}

The hyperbolic plane $H^2$ has a well-known model, sometimes called the
Klein model; it is a map $\phi$ from $H^2$ to the open disk $D^2$ of
radius $1$ in $\R^2$ with the striking properties that geodesics are
sent to segments. It can be obtained in several different ways, two of
which are as follows.

Consider the hyperboloid model of $H^2$ in Minkowski 3-space,
$\R^3_1$, wich is $\R^3$ with the Lorentz metric $-x_0^2+x_1^2+x_2^2$:
$$ H^2 = \{ x\in \R^3_1 ~ |  ~ \langle x, x\rangle = -1 ~\mbox{and}~
x_0>0 \}~. $$
The hyperbolic geodesics are then the intersections of $H^2$ with the
2-planes containing $0$. Project $H^2$ on the plane $P_0=\{ x_0=1\}$ in the
direction of the origin, by sending $(x_0, x_1, x_2)\in H^2$ to $(1,
x_1/x_0, x_2/x_0)$. The image of $H^2$ is $D^2$, and the image of the
intersection of $H^2$ with a plane $P\ni 0$ is the intersection of $P$
with the unit disk in $P_0$. 

The second approach uses the cross-ratio. If $x,y,a$ and $b$ are four
real numbers, their cross-ratio $[x,y;a,b]$ is defined as:
$$ [x,y;a,b] := \frac{(x-a)(b-y)}{(y-a)(b-x)}~. $$
Given four points on a line in $\R^2$, their cross-ratio is defined in
the same way, by replacing $x-a$ by the oriented distance between $a$
and $x$, etc. An important property is that the cross-ratio is invariant
under projective transformations.

If $x,y\in D^2$ are distinct, let $l$ be the line containing them,
and let $a,b$ be the intersection points of $l$ with $S^1=\dr
D^2$. Suppose that $a,x,y,b$ are in this order on $l$, and define the
Hilbert distance $d_H(x,y)$ as half the log of the cross-ratio of
those four points on $l$:
$$ d_H(x,y) = -\frac{1}{2} \log [x,y;a,b]~. $$
It is a result of Hilbert that $(D^2, d_H)$ is isometric to $H^2$,
and that its geodesics are the segments of $D^2$.

To each geodesic $g$ in $H^2$, one can associate a unique point in
$\R^2\setminus \overline{D^2}$, which we will call the {\bf dual} of
$g$, by a simple and classical projective construction. If $g$ is, in
the projective model, the intersection of $D^2$ with a line $l$, let
$a,b$ be the intersections of $l$ with $S^1=\dr D^2$; the dual of $g$ is
the intersection point of the tangents to $S^1$ at $a$ and $b$. Note
that this point can be at infinity if $l\ni 0$. Given a line $l\subset
\R^2$ which intersects $D^2$, we will also call "dual of $l$" the point
obtained in $\R^2\setminus \overline{D^2}$ by this construction. 
Given a point $v\in\R^2\setminus \overline{D^2}$, there is a unique line
$l$ such that the dual of $l$ is $v$; we will call $l$ the dual of
$v$. The dual of a point $v$ will be denoted by $v^*$, and the dual of
a line $l$ by $l^*$. 

There is another, equivalent definition of the dual of a point $v\in
\R^2\setminus \overline{D^2}$ in terms of the cross-ratio. For each
point $x\in D^2$, let $l$ be the line containing $v$ and $x$, and let
$a,b$ be the intersection points between $l$ and $S^1$, chosen so that
$x,a,v,b$ appear in this order on $l$. $v^*$  is then the set of all
points $x\in D^2$ such that $[x,v;a,b]=-1$. To check that it does not
contradict the previous definition, check the case when $v$ is a point
at infinity and use the projective invariance of both definitions. 

Of course, this duality is the same as the one defined in the previous section
by considering $H^3$ and $S^3_1$ as quadrics. 

\bprop \label{pr:ortho}
Let $v\in \R^2\setminus D^2$, and let $d$ be a line going through $v$
and intersecting $D^2$. Then $v^*$ intersects $d$ orthogonally for the
hyperbolic metric on $D^2$.
\eprop

\bpv
By projective invariance of the cross-ratio (and thus of the hyperbolic
metric) we only need to prove the proposition when $v$ is a point at
infinity, for instance if it corresponds to the vertical direction in
$\R^2$. In this case, for any $x\in D^2$, if $a$ and $b$ are the
intersection of $\dr D^2$ with the vertical line containing $x$, then:
$$ x\in v^* ~ \Leftrightarrow ~ \frac{ax}{xb}=1~. $$ 
Thus $v^*$ is the horizontal diameter $d$ of $D^2$. An elementary symmetry
argument shows that $d$ is orthogonal, for the hyperbolic metric $d_H$,
to the vertical lines.
\epv

\bprop \label{pr:no-inter}
Let $s$ be a segment in $\R^2$, with endpoints $x,y\in \R^2\setminus
\overline{D^2}$, but with $s$ intersecting $D^2$. Then the lines $x^*$
and $y^*$,dual to $x$ and $y$ respectively,do not intersect in $D^2$.
\eprop

\bpv
This is a consequence of the "geometric" definition of the dual of a
point. Using the projective invariance, we can suppose for instance that
both $x$ and $y$ lie on the $x$-axis, on opposite sides of $D^2$. Then
$x^*$ and $y^*$ lie on opposite sides of the $y$-axis, and can therefore
not intersect.
\epv

\pg{Some definitions}

The analog in dimension 2 of hyperideal polyhedra are the hyperideal
polygons; in particular, the hyperideal triangles.

\bdf 
A {\bf hyperideal triangle} is a triangle in $\R^2$ whose vertices are outside
the open disk $D^2$ and whose edges all intersect $D^2$. It is {\bf strictly
hyperideal} if its vertices are outside the closure $\overline{D^2}$.
\edf

The following definition is made possible by proposition
\ref{pr:no-inter}. 

\bdf \label{df:truncated}
Let $T$ be a strictly hyperideal triangle. The {\bf truncated hyperideal
  triangle} associated to $T$ is the hexagon obtained by
taking the intersection of $T$ with the half-planes bounded by the lines
dual to its vertices (and not containing the vertices).
\edf

\begin{figure}[h]
\centerline{\psfig{figure=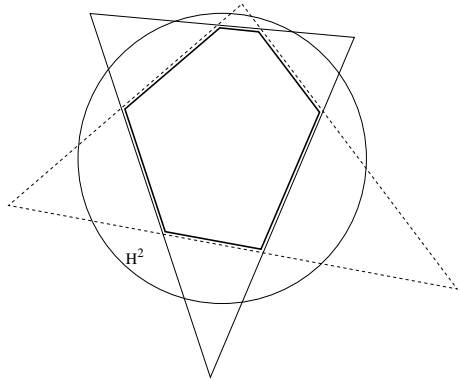,height=5cm}}
\caption{A hyperideal triangle and its truncated version}
\end{figure}

Proposition \ref{pr:ortho} indicates that the truncated hyperideal
triangles are actually right-angle hexagons. They have two kinds of
kinds of edges: the remaining parts of the edges of the non-truncated
triangle $T$, which we will call {\bf real edges}, and the intersections
of $T$ with the lines dual to its vertices, which we will call {\bf
  cuts}. The {\bf edges lengths} of a truncated hyperideal triangle are
the lengths of its real edges; the edge lengths of a hyperideal
triangle are the edge lengths of the corresponding truncated triangle.

\pg{Edge lengths}

We will need to understand what are the possible edge lengths of
hyperideal triangles. 

\bprop \label{pr:rig-1}
The edge lengths of any strictly hyperideal triangle are positive real
numbers.  
For each $l_1, l_2, l_3\in \R^+\setminus \{ 0\}$, there is at most one
hyperideal triangle with edge lengths $l_1, l_2, l_3$. Moreover, given a
strictly 
hyperideal triangle, it has no infinitesimal deformation which does not
change its edge lengths.
\eprop

\bpv
The first statement is a consequence of proposition \ref{pr:no-inter}. 

For the second statement, let $v\in \R^2\setminus
\overline{D^2}$ and choose a number $l>0$. Now let $x\in \R^2$ be such
that the line $(vx)$ intersects $S^1$ at two points $a,b$ such that
$v,a,b,x$ lie on $(vx)$ in that order. Notice then that
$[v,x;a,b]=e^{-2l}$ if and only if $x$ lies on an ellipse which
is tangent to $S^1$ at the points of $S^1\cap v^*$. Checking this can be
done by choosing $v$ as a point at infinity and then using the
projective invariance again. 

Now let $l_1, l_2$ and $l_3$ be positive real numbers, and choose a
point $x_1$ in $\R^2\setminus \overline{D^2}$. The points $x_2\in
\R^2\setminus \overline{D^2}$ such that $[x_1,x_2]$ intersects $D^2$ and
that its length is $l_3$ form part of an ellipse (tangent to $S^1$ at
$S^1\cap x_1^*$); choose one of them. 

The points $x\in \R^2\setminus \overline{D^2}$ such that $[x_1,x]$
intersects $D^2$ and that its length is $l_2$ then form a subset $E_1$
of an ellipse tangent to $S^1$ at $S^1\cap x_1^*$, while the points
$y\in \R^2\setminus \overline{D^2}$ such that $[x_2,x]$ 
intersects $D^2$ and that its length is $l_1$ form a subset $E_2$
of an ellipse tangent to $S^1$ at $S^1\cap x_2^*$. $E_1$ and $E_2$
intersect at at most two points (see figure 3). 
Taking either of those intersections as $x_3$
yields a triangle (with vertices $x_1, x_2$ and $x_3$)
with edge lengths $l_1, l_2$ and $l_3$.

To prove the uniqueness of this triangle and the third statement, just
notice that the choices of the points $x_1$ and $x_2$ are equivalent
under the action of the group of orientation-preserving projective
transformations of $\R^2$ fixing $S^1$, while the two possible choices
of $x_3$ lead to triangles which are equivalent under an orientation
inversing projective transformation fixing $S^1$. 
\epv

\begin{figure}[h]
\centerline{\psfig{figure=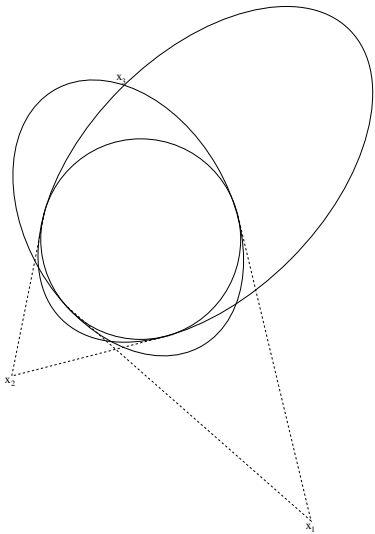,height=7cm}}
\caption{How to position $x_3$ knowing $x_1$ and $x_2$}
\end{figure}

\pg{A projective model of the de Sitter space}
The construction of a projective model for the hyperbolic plane, which
was described above, extends with few modifications to the de Sitter
plane. Recall that the de Sitter plane $S^2_1$ is a constant curvature $1$
Lorentz surface; it can be obtained, much like the hyperbolic plane, as
a quadric in $\R^3_1$, with the induced metric:
$$ S^2_1 \simeq \{ x\in \R^3_1 ~| ~ \langle x,x\rangle =1\}~. $$
Its geodesics are the intersections with $S^2_1$ of the planes
containing $0$ in $\R^3_1$. They are of three kind: the space-like,
light-like and time-like geodesics, corresponding respectively to the
space-like, light-like and time-like planes in $\R^3_1$. 

The projective model which we will consider is actually only for "half"
the de Sitter plane. It can be obtained by projecting the part of
$S^2_1$ which lies in the half-space $\{ x_0>0\}$ to the plane $P_0:= \{
x_0=1\}$ along the direction of $0$, i.e. by sending a point $(x_0, x_1,
x_2)$ with $x_0>0$ to the point $(1, x_1/x_0, x_2/x_0)\in P_0$. Note
that another possible construction uses the Hilbert distance (see
e.g. \cite{shu}).

\pg{Triangles in the de Sitter space}

Some questions concerning the dihedral angles of polyhedra, which we
will encounter below, can be understood using elementary properties of
the triangles in the de Sitter plane. We state some such results here
for future reference.

\bdf \label{df:space}
A {\bf space-like} triangle in the de Sitter plane is the image by the
projective model of a triangle in $\R^2$ which bounds an open convex set
which contains $\overline{D^2}$.
\edf

Note that the edges of the "space-like" triangles are segments of
space-like geodesics; but there are triangles in $S^2_1$ with space-like
edges which are not "space-like triangles" as we have just defined
them. 

\bprop \label{pr:triangles}
The sum of the edge lengths of any space-like triangle is strictly
greater than $2\pi$. Conversely, given $l_1, l_2, l_3\in
(0, \pi)$ with sum strictly greater than $2\pi$, there is a unique
space-like triangle in $S^2_1$ with those edge lengths.
\eprop

\bpv
The space-like triangles are exactly the triangles, in the de Sitter
plane, which are duals of compact hyperbolic triangles. So, using this
duality, proving the proposition is equivalent to showing that:
\begin{itemize}
\item the sum of the exterior angles of any hyperbolic triangle is
  strictly larger than $2\pi$.
\item for each triple $(\alpha_1, \alpha_2, \alpha_3)\in (0, \pi)^3$ such that
$\alpha_1+\alpha_2+\alpha_3>2\pi$, there is
a unique hyperbolic triangle with $\alpha_1, \alpha_2$ and $\alpha_3$ as
its exterior angles.
\end{itemize}
This is easily seen to be true by switching to the interior angles. The
uniqueness is of course up to the global hyperbolic isometries. 
\epv

\pg{A topology on the space of hyperideal triangles}

There are several natural ways to put a topology on the space of
hyperideal triangles. One can use for instance:
\begin{itemize}
\item the lengths of the edges (which uniquely characterize the
  triangles up to orientation).
\item the positions of the vertices in $\R^2$, modulo the action of the
  projective transformations leaving $D^2$ invariant.
\item the angles at the vertices, which also characterize the triangles
  up to isometry. 
\end{itemize}
We will use here the third solution. It has the advantage that the
triangles with one or more ideal vertices appear simply as the triangles
with one or more angles equal to $0$. In terms of the lengths of the
edges, the triangles with one or more ideal vertices have some infinite
lengths and this makes things a little more troublesome. 

\section{Hyperideal simplices}

\pg{Definitions}

We will be using in this section the projective model of $H^3$. Since it
is quite analoguous to the projective model of $H^2$ defined in
the previous section, we will not give any further detail on its
construction or on its basic properties. Here $D^3$ is the open ball of
radius $1$ in $\R^3$, which in the projective model is the image of
$H^3$. 

\bdf \label{df:simplices}
A {\bf hyperideal simplex} is a simplex in $\R^3$ with all its vertices
in $\R^3\setminus D^3$, and such that all its edges intersect $D^3$. It
is {\bf strictly hyperideal} if none of its vertices lie on $S^2=\dr
D^3$. An hyperideal simplex is {\bf degenerate} if it lies in a
hyperbolic plane. 
\edf

We will call $\cS$ the set of hyperideal simplices, considered up to
global hyperbolic isometries. $\cS$ can be decomposed as $\cS=\cS_0\cup
\cS_1\cup \cS_2\cup \cS_3\cup \cS_4$, where $\cS_i$ is the set of
hyperideal simplices with exactly $i$ ideal vertices; for instance,
$\cS_0$ is the set of strictly hyperideal vertices, while $\cS_4$ is the
set of ideal simplices.

Given a point $v\in \R^3\setminus \overline{D^3}$, its dual is a
hyperbolic totally geodesic plane which can be
defined, as in the previous section, in many ways including the
following two:
\begin{itemize}
\item if $C$ is the set of points in $x\in S^2$ such that the line
  $(xv)$ is tangent to $S^2$ at $x$, then $C$ is a circle on $S^2$ ---
  i.e. the boundary of a geodesic ball. It is thus the boundary at
  infinity of a plane, which we define to be $v^*$. 
\item $v^*$ is the set of points $x\in D^3$ such that $[x,v;a,b]=-1$,
  where $a$ and $b$ are the intersections with $S^2$ of the line $(v,x)$.
\end{itemize}

The analog of proposition \ref{pr:no-inter} can be proved in the same
way: 

\bprop \label{pr:no-inter-2}
Let $s$ be a segment in $\R^3$, with endpoints $x,y\in \R^3\setminus
\overline{D^3}$, but with $s$ intersecting $D^3$. Then the duals of $x$
and $y$ do not intersect. 
\eprop

\pg{Truncated simplices}

Given a hyperideal simplex $S$, we define (following Bao and Bonahon
\cite{bao-bonahon}) the associated {\bf truncated
hyperideal simplex} $\Sb$ as the intersection of $S$ with the half-spaces
bounded by the duals of its vertices (and which do not contain those
vertices). 

By proposition \ref{pr:no-inter-2}, the truncated
strictly hyperideal simplices are compact; they  have two kinds of faces,
hexagons, which we will call "real faces", and the triangles which are
the intersections of $S$ with the planes dual to its vertices. They also
have two kinds of edges: the  "real edges" are the intersections
of the edges of $S$ with $\Sb$, and the "edge lengths" of $S$ are the
lengths of the "real edges" of $\Sb$. 

On the other hand, truncation does not change the ideal ends of a
simplex, and if $S$ is a simplex with at least one ideal vertex,
then its associated truncated simplex has finite volume but is not
compact.

\pg{Edge lengths}

As a direct consequence of proposition \ref{pr:rig-1}, we find an
analoguous rigidity statement in dimension $3$. Before proving it, we
need a preliminary statement on the set of points at a given distance
from a hyperideal or ideal point.

Note that, given two points $x,y\in \R^3\setminus \Db$, one can define
the distance between them. This can be done in at least two ways, which
can directly be checked to be equivalent.
\begin{itemize}
\item using the Hilbert distance, as in section 2.
\item from the distance between their dual planes, when the segment
  $[x,y]$ intersects $\Db$.
\end{itemize}
Similarly, there is a natural notion of distance between a point in
$\Db$ and a point in $\R^3\setminus \Db$, which can be defined from the
Hilbert metric or from the distance from the hyperbolic point to the
dual of the hyperideal point.

There is also a notion of distance between an ideal point $x$ an a
hyperbolic point $y$, but it is not canonically defined and depends on
the choice of a horosphere $H$ centered at $x$; once $H$ is given, the
distance between $x$ and $y$ is defined as the distance between $y$ and
$H$. Replacing $H$ by another horosphere $H'$ changes the distances from all
points in $H^3$ to $x$ by the same constant (which is the distance
between $H$ and $H'$).

\bprop \label{pr:equidist}
Consider the projective model of $H^3$ and $S^3_{1,+}$. 
\begin{enumerate}
\item Let $x\in S^3_{1,+}$ be a hyperideal point, and let $d>0$. The set
  of points in $S^3_{1,+}$ at constant distance $d$ from $x$ is, in the
  projective model, an ellipsoid of revolution,
  which is tangent to $\dr_\infty H^3$ along
  the intersection of $\dr_\infty H^3$ with the tangent cone with vertex
  $x$. 
\item Let $x\in \dr_\infty H^3$ be an ideal point, and let $H$ be a
  horosphere at $x$. The set of points $y\in S^3_{1,+}$ such that the
  distance from $y$ to $x$ relative to $H$ is fixed is an
  ellipsoid of revolution, which is outside $B(0,1)$, but tangent to
  $S^2$ at $x$. 
\item Suppose that $x_1, x_2, x_3$ are either ideal or hyperideal
  points, and let $d_1, d_2, d_3\in \R$. For each ideal point in $\{
  x_1, x_2, x_3\}$, choose a horosphere $H_i$ centered at $x_i$. Suppose
  that the intersection of the ellipsoids $E_i$ of points in the projective
  model, at distance 
  $d_i$ from $x_i$ (relative to $H_i$) contains more than two
  points. Then $x_1, x_2$ and $x_3$ lie on a line. 
\end{enumerate}
\eprop

\bpv
For (1), apply a projective transformation sending $x$ to
infinity, for instance to the point at infinity corresponding to the
"vertical" direction, without moving $S^2$. Then $H^3$ remains in a
ball, while the cone 
tangent to $\dr_\infty H^3$ is sent to a vertical cylinder tangent to
$\dr_\infty H^3$ along the intersection of $\dr_\infty H^3$ with the
dual plane $x^*\subset H^3$. 

We are interested in the set of points at constant distance from $x$,
i.e. in the set of points $y$ such that, if $D$ is the line through $x$
and $y$ and $a,b$ are the intersections of $D$ with $\dr_\infty H^3$
(such that $x,a,b,y$ appear in this order on $D$) the cross-ratio
$[x,y;a,b]$ is equal to a constant $C$. After the projective
transformation sending $x$ to infinity, this equation becomes:
$ay=C.by$, which can be translated as: $(C-1)yb=C.ab$. This is clearly
the equation of an ellipsoid of revolution, tangent to $\dr_\infty H^3$
along $\dr_\infty H^3\cap x^*$. 

\medskip

For (2), let $y\in S^3_{1,+}$, and note that the distance from $x$ to
$y$ relative to $H$ is $C$ if and only if the distance from $x$ to the
plane $y^*$ relative to $H$ is $C$. Thus the set $S'$ of points at distance
$C$ from $x$ relative to $H$ is the dual of the set $S$ of points $z\in H^3$
at distance $C$ from $x$ relative to $H$. But $S$ is clearly a
horosphere centered at $x$. So, in the projective model, $S$ is an
ellipsoid of revolution tangent to $\dr_\infty H^3$ at $x$. It is then a
simple matter of projective geometry to check that $S'$ is also, in the
projective model, an ellipsoid of revolution tangent to $\dr_\infty H^3$
at $x$. 

\medskip

For (3), note that we can suppose that, maybe after changing their labels,
$x_1$ and $x_2$ are either both ideal points or both strictly hyperideal
points. Moreover, after applying a projective tranformation, we can
suppose that $x_1$ and $x_2$ are either on the same vertical line going
through $0$, or symmetric with respect to the vertical line containing
$0$. In both cases, the sets of points at distance $d_i$ from $x_i$
(maybe relative to $H_i$), for $i\in \{ 1,2\}$, is a circle $C$ in a plane
in the projective model. 

Since the $E_i$, $1\leq i\leq 3$, are quadrics, their intersection
contains either at most 2 points, or a whole curve. So, to prove the
statement, we can suppose that $E_1\cap E_2\cap E_3=C$. 
Now it is not difficult to check that the set of points at constant
distance from a circle like $C$ is a line in the projective model;
as a consequence, $x_1, x_2$ and $x_3$ are on a line. 
\epv

We can describe more precisely the set of points at
given distance from an ideal point $x$. First, a simple computation
shows that, in the projective model, the horospheres centered at $x$ are
simply the ellipsoids with radii $\lambda, \lambda$ and $\lambda^2$, for
$\lambda \in (0,1)$, which are tangent to $S^2$ at the end of the small
axis. Then the projective duality shows that the sets of points in
$S^3_1$ at constant distance from $x$ are the duals of those ellipsoids,
which have the same description except that now $\lambda>1$, so that the
tangency to $S^2$ occurs at the end of the large axis.

We call
$e_{12}, \cdots, e_{34}$ the edges of the simplex $S_0$
(considered as a combinatorial object). 

\bprop \label{pr:rig-2}
Let $l_{12}, \cdots, l_{34}$ be positive real numbers. There is at most
one 
strictly hyperideal simplex $S$ such that the length of $e_{ij}$ is
$l_{ij}$. There is no 
first-order deformation of $S$ which does not change its edge lengths. 
\eprop

\bpv
Proposition \ref{pr:rig-1} shows that there is at most one hyperideal
triangle with edge lengths $l_{12}, l_{13}$ and $l_{23}$. Let $x_1, x_2,
x_3$ be the three vertices obtained in this way. Finding the last vertex
$x_4$ is equivalent to finding the intersection of the sets of points
$E_1, E_2$ and $E_3$ at
given distances $l_{14}, l_{24}$ and $l_{34}$, respectively, from $x_1,
x_2$ and $x_3$. 

According to point (3) of proposition \ref{pr:equidist}, since $x_1,
x_2$ and $x_3$ do not lie on a line, the intersection of $E_1, E_2$ and
$E_3$ contains at most two points, which are exchanged by a hyperbolic
transformation fixing $x_1, x_2$ and $x_3$. This proves the uniqueness
of the simplex. 

The same argument shows the infinitesimal rigidity statement, using the
infinitesimal rigidity part of proposition \ref{pr:rig-1}.
\epv

\bcr \label{cr:lengths}
For each function $l:\{ e_{12}, \cdots, e_{34}\}\rightarrow
\R^+\setminus \{ 0\}$, there is exactly one strictly hyperideal simplex
with edge lengths given by $l$.
\ecr

\bpv
We consider the map $F$ sending a strictly hyperideal simplex (with
vertices labeled from $1$ to $4$) to the set of its edge lengths. The
statement follows from the following points, each of which can readily
be checked.
\begin{enumerate}
\item the set of strictly hyperideal simplices, and the space of
  functions $l:\{ e_{12}, \cdots, e_{34}\}\rightarrow
  \R^+\setminus \{ 0\}$, both have dimension $6$. 
\item $F$ is locally injective (i.e. its differential is everywhere
  injective) according to the infinitesimal rigidity statement of
  proposition \ref{pr:rig-2}.
\item $F$ is proper, i.e. if a sequence of simplices $(S_i)_{i\in \ N}$
  is such that the corresponding length functions $(l_i)_{i\in \N}$
  converges, then $(S_i)$ converges. 
\item the space of strictly  hyperideal simplices is connected, while
  the space of admissible distance functions is simply connected.
\end{enumerate}
\epv

\pg{The hyperbolic-de Sitter duality}

The (classical) construction which was given in section 2 of a
projective model for the hyperbolic plane and the 2-dimensional de
Sitter space extends, with minor modifications, to dimension 3. So does
the duality between the hyperbolic and de Sitter space. The dual of
a point in $H^3$ is now a 2-dimensional totally geodesic space-like
plane in $S^3_1$, while the dual of a point in $S^3_1$ is an oriented
totally geodesic plane in $H^3$. 

One important property of this duality concerns the polyhedra; the next
proposition was discovered by Rivin and Hodgson \cite{Ri,RH}. It is an
extension of the properties described in the previous section in
dimension 2. 

\bprop
Let $P$ be a compact polyhedron in $H^3$. Its dual is a compact,
space-like polyhedron in $S^3_1$. To each edge $e$ of $P$ corresponds an
edge $e^*$ of $P^*$, and the exterior dihedral angle at $e$ is the
length of $e^*$. 
\eprop

The proof is can be found in a number of sources, e.g. \cite{RH,shu}, so
we leave it to the reader.

This duality is valid not only for compact polyhedra, but also
for ideal or hyperideal polyhedra, and also for smooth, strictly convex
surfaces --- the dual objects are then the smooth, strictly convex
space-like surfaces in $S^3_1$, see e.g. \cite{these}. 

The group of orientation-preserving isometries of $H^3$, $SO(3,1)$, is
also the group of orientation-preserving isometries of $S^3_1$ which do
not exchange the two boundary components of $S^3_1$. Therefore, any
group action on $H^3$ has an extension to $S^3_1$. Since the duality
described above is defined "geometrically", it is not difficult to show
that it "commutes" to the action of the isometries of $H^3$
resp. $S^3_1$. 
This means that, given a polyhedral surface in $H^3$ which is invariant
under a subgroup $\Gamma$ of $SO(3,1)$, its dual is invariant under the
action of $\Gamma$ on $S^3_1$. 

\medskip

\bdf 
A {\bf dual hyperideal simplex} is a simplex in $\R^3$ with all vertices
and edges in 
$\R^3\setminus \overline{D^3}$ but with all faces intersecting $D^3$
(maybe at one point).
\edf

Note that we might use the same terminology to describe the intersection
of a dual hyperideal simplex with $\R^3\setminus \overline{D^3}$, and also the
corresponding non-complete polyhedron in de Sitter space.

\bprop
The dual of the hyperideal simplices are the dual hyperideal simplices.
\eprop

\bpv
This follows from the elementary properties of the hyperbolic-de
Sitter duality. If $S$ is a hyperideal simplex, each of its vertices is
either strictly hyperideal, or ideal; so the dual planes intersect
$H^3$, either on a  disk (for strictly hyperideal vertices) or at
a point (for ideal vertices). Moreover, each edge of $S$ intersects
$H^3$, so that the dual edges are space-like geodesics in $S^3_1$, which
remain outside $H^3$ in the projective model. The converse is
proved in the same way. 
\epv

The following proposition is a special case of results of
\cite{shu}, but we will outline its proof --- which is elementary ---
for completeness. 

\bprop \label{pr:dual-hyperideal}
The possible edge lengths of a dual hyperideal simplex $S$ are
the functions 
$l:\{ e_{12}, \cdots, e_{34}\}\rightarrow (0, \pi)$ such that, for each
face $f$ of $S$, the sum over the edges of $f$ of the values of
$l$ is at least $2\pi$. For each such function, there is a unique
simplex with those edge lengths.
\eprop

The uniqueness in this statement is of course up to the global
hyperbolic isometries. 

\bpv
Each face of a dual hyperideal simplex is isometric to a triangle:
\begin{itemize}
\item either in the
projective model of $H^2$ and $S^2_{1,+}$, it is then the dual of a
hyperbolic triangle (i.e. which contains $H^2$ in its interior). 
\item or in the degenerate space $H^2_{1,0}$, which is isometric to a
  light-like totally geodesic plane in $S^3_1$, and it then contains the
  limit point in its interior.
\end{itemize}
The edge lengths of those triangles are:
\begin{itemize}
\item in the first case, exactly the triples $(l_1,l_2,l_3)$ of numbers
  in $(0,\pi)$ such that $l_1+l_2+l_3> 2\pi$.  
\item in the second case, the triples  $(l_1,l_2,l_3)$ of numbers in
  $(0,\pi)$ such that $l_1+l_2+l_3= 2\pi$. 
\end{itemize}
Moreover, each
triple determines a unique triangle. Thus, the lengths of the edges of a
dual hyperideal simplex satisfy the hypothesis of the proposition, and
moreover its faces are uniquely determined. It is then a simple matter
to check that there is a unique way of gluing those faces to obtain a
simplex.
\epv

\pg{Dihedral angles of hyperideal simplices}

In a dual way, we need to understand the dihedral angles of the
hyperideal simplices. The next lemma is a very special case of the main
result of \cite{bao-bonahon}. We will however give a direct proof for
completeness. The proof is a direct consequence of the previous
proposition and the hyperbolic-de Sitter duality.

\blm \label{lm:dihedral}
Let $S$ be a hyperideal simplex. Its exterior dihedral angles are such
that, for each vertex $s$ of $S$, the sum of the angles of the edges
containing $s$ is greater than $2\pi$, and equal to $2\pi$ if and only
if $s$ is ideal. Moreover, given a map
$\alpha:\{ e_{12}, \cdots, e_{34}\}\rightarrow (0, \pi)$ such that, for each
vertex $s$ of $S_0$, the sum of the values of $\alpha$ on the edges of
$S_0$ incident to $s$ is at least $2\pi$, there exists a
unique hyperideal simplex such that the exterior dihedral angle at each edge
$e_{ij}$ is $\alpha(e_{ij})$. 
\elm

Note that the path length condition which appears in lemma
\ref{lm:angles-poly}, or in theorem \ref{tm:angles}, is redundant here;
a simple argument shows that it is always satisfied under the hypothesis
of lemma \ref{lm:dihedral}.

Another consequence is that the strictly hyperideal or ideal
nature of the vertices can be read from the dihedral angles:

\brk \label{rk:ideaux}
Let $S$ be a hyperideal simplex. A vertex $v$ of $S$ is ideal if and
only if the sum of the exterior dihedral angles of the edges adjacent to
$v$ is $2\pi$. 
\erk

This clearly remains true for hyperideal polyhedra, or for hyperideal
manifolds.

\pg{The Schl{\"a}fli formula}

It is a key element in this paper, so we recall it here. The reader
might find a proof e.g. in \cite{milnor-schlafli,geo2}.

\blm \label{lm:schlafli}
Let $P$ be compact hyperbolic polyhedron, with edge lengths $(L_i)$ and
dihedral angles $(\theta_i)$. In a first-order deformation of $P$, the
variation of its volume is given by:
\beq \label{eq:schlafli}
dV = -\frac{1}{2} \sum_i L_id\theta_i~.
\eeq
\elm

To understand the Schl{\"a}fli formula for hyperideal simplices with some
ideal vertices, we introduce linear map as follows:
$$
\begin{array}{llcr}
\phi_4: & \R^4 & \rightarrow & \R^6 \\
& (x,y,z,t) & \mapsto & (x+y, x+z, x+t, y+z, y+t, z+t) 
\end{array}
$$
We then call $\phi_1$ the restriction of $\phi_4$ to $\R\times \{(0,0,0)\}$
(identified with $\R$), $\phi_2$ the restriction of $\phi_4$ to $\R^2\times
\{(0,0)\}$, and $\phi_3$ the restriction of $\phi_4$ to $\R^3\times \{ 0\}$. 

Now let $S$ be a hyperideal simplex with one ideal vertex $v$
exactly. Let $\Sb$ be the associated truncated hyperideal simplex, so
that $\Sb$ is a finite volume polyhedron in $H^3$ with exactly one ideal
vertex. 
Choose a horosphere $H$ centered at $v$ (which we will require to be
small enough), and define the {\bf edge
  lengths} of $S$  to be the lengths of the real edges of $\Sb$ of
finite length, and the lengths of the segment of the other edges up to
their intersection with $H$ (this intersection exists if $H$ is small
enough). Of course the resulting edge lengths depend on the choice of
$H$; but changing $H$ only adds a constant to the three edges incident
to $v$, so that the lengths of the edges of $S$ are well defined as an
element of $\R^6/\phi_1(\R)$. 

In the same way, if $S$ has two ideal simplices, its edge lengths are
defined as an element of $\R^6/\phi_2(\R^2)$; for 3 ideal simplices they
are in $\R^6/\phi_3(\R^3)$, and for 4 ideal vertices, in
$\R^6/\phi_4(\R^4)$. 

With those natural definitions, we can give an extension --- also classically
known --- of lemma \ref{lm:schlafli}. 

\blm \label{lm:schlafli-2}
Let $P$ be hyperideal polyhedron, with edge lengths $(L_i)$ and
dihedral angles $(\theta_i)$. In a first-order deformation of $P$ which
leaves its ideal vertices on the sphere at infinity, the
variation of its volume is given by (\ref{eq:schlafli}):
$$
dV = -\frac{1}{2} \sum_i L_id\theta_i~.
$$
\elm

Note that this formula makes sense although the edge lengths are defined
in general only up to the addition of a constant for each ideal vertex,
because the sum of the dihedral angles of the edges containing an ideal
vertex are constrained to be $2\pi$, so that the sum of their
differential vanishes --- the additive constant in the lengths therefore
doesn't make any difference.

\pg{Infinitesimal rigidity for simplices}

The definitions of edge lengths for hyperideal simplices also lead to an
extension of proposition \ref{pr:rig-2}. 
 
\bprop \label{pr:rig-3}
Let $i\in \{ 0,1,2,3,4\}$, and let $S\in \cS_i$ be a hyperideal
simplex, with edge lengths $l:=(l_{ij})\in \R^6/\phi_i(\R^i)$. There is
no non-trivial first-order deformation of $S$ in $\cS_i$ which does not
change $l$. 
\eprop

\bpv
The proof depends on the number $n$ of ideal vertices. In each case we
consider an infinitesimal deformation $\Sd$ of $S$ which does not change the
edge lengths, and then show that it is trivial.

{$n=0$:} this was already done in proposition \ref{pr:rig-2}.

{$n=1$:} let $x_1$ be the ideal vertex, so that $x_2, x_3$ and $x_4$
are strictly hyperideal. By proposition \ref{pr:rig-1}, we can apply an
infinitesimal isometry to get in the situation where $x_1, x_2$ and
$x_3$ are fixed by $\Sd$. Moreover, we can fix the horosphere centered
at $x_1$, so that the distance between $x_1$ and $x_4$ is
well-defined. Then proposition \ref{pr:equidist} shows that $x_4$ has to
be in the intersection of three ellipsoids (corresponding to the
fixed distances to $x_1, x_2$ and $x_3$) so that $\Sd$ has to be trivial. 

{$n=2$:} let $x_1$ and $x_2$ be ideal. By adding an infinitesimal isometry, we
can suppose that $\Sd$ fixes $x_1, x_2$ and $x_3$. We 
also fix horospheres centered at $x_1$ and $x_2$, and apply the same
argument as for $n=1$.

{$n=3$:} we use the same line of reasoning; note that, up to
isometries, there is only one ideal triangle, so we can suppose that
$x_1, x_2$ and $x_3$ are fixed by $\Sd$. In addition there is a unique
way to choose horospheres centered at $x_1, x_2$ and $x_3$ which are
pairwise tangent. $x_4$ is then restrained to be in the intersection of
three ellipsoids and so $\Sd$ is trivial.

{$n=4$:} fix $x_1, x_2$ and $x_3$ (this is possible since they form
an ideal triangle), and the corresponding horospheres so that they are
pairwise disjoint. Then we can still suppose that the horosphere
centered at $x_4$ is tangent e.g. to the horosphere centered at
$x_1$. This leaves two conditions, corresponding to the distance between
the horosphere centered at $x_4$ and those centered at $x_2$ and
$x_3$. It is a rather simple matter of Euclidean geometry
to check that those conditions are
non-degenerate, i.e. that there is no possible displacement of $x_4$
which satisfies them at first order. 
\epv

Note that the last case of the previous proposition, for ideal
simplices, is well-known. One way to prove it (see \cite{ideal} for a
more general result on polyhedra or manifolds with ideal boundary) is in
the direction opposite to what we do below: first show that the volume is
a concave function of the dihedral angles, then use the Schl{\"a}fli formula
to deduce the infinitesimal rigidity with respect to the edge lengths
(or induced metric).

\pg{The volume of simplices}

Another key element, related to the Schl{\"a}fli formula and the rigidity
phenomenon which we have just seen, is the fact that
the volume of the hyperideal simplices (like the volume of the ideal
simplices, see e.g. \cite{Ri2}) is concave. 

\bdf
The {\bf volume} of a hyperideal simplex is defined as the volume of the
corresponding truncated hyperideal simplex. 
\edf

Note that, for each $i\in \{0,\cdots,4\}$, $\cS_i$ is parametrized by
the dihedral angles, subject to the condition that, for each ideal
vertex, the sum of the angles of the incident edges is $2\pi$. 
By lemma  \ref{lm:dihedral} and remark \ref{rk:ideaux}, $\cS_i$
is affinely equivalent to a convex polytope in $\R^{6-i}$. We call $V$ the
volume, seen as a function on each $\cS_i, 0\leq i\leq 4$.

\bprop \label{pr:regulier}
There exists a regular hyperideal simplex such that $\hess(V)$ is
negative definite. 
\eprop

\bpv
Let $S_0$ be a regular hyperideal simplex, with edge lengths equal to
$l_0$ and dihedral angles equal to $\theta_0$. 
By the Schl{\"a}fli formula, lemma \ref{lm:schlafli-2}, 
the matrix of $\hess(V)$ is equal to:
$$ \hess(V) = - \frac{1}{2} \left[\frac{\dr l_i}{\dr
    \theta_j}\right]_{i,j}~, $$ 
where $\theta_j$ is the exterior dihedral angle at the edge $j$. 
So we only have to prove that, for some regular hyperideal simplex, the
matrix $(\dr \theta_i/\dr l_j)_{1\leq i,j\leq 6}$ is positive definite. 

Using the symmetry of the regular simplices, this can be done as
follows. 
\begin{itemize}
\item consider a hyperideal simplex with five edge lengths equal to
  $l_0$, and one, say the length of $e$, equal to $l$. 
\item use the hyperideal version of some classical triangle formulas to
  compute the angles of the faces. There are three such angles, one,
  $\alpha_0$, which is the angle of the faces of the regular simplex
  with edge lengths $l_0$, and two others, $\alpha$ and $\beta$.
\setlength{\unitlength}{0.5cm}
\begin{figure}[h] \label{fig:simplex}
\begin{picture}(8,6)(-18,-3)
\put(-4,0){\dashbox{0.1}(8,0)}
\put(-4,0){\line(4,3){4}}
\put(-4,0){\line(4,-3){4}}
\put(4,0){\line(-4,3){4}}
\put(4,0){\line(-4,-3){4}}
\put(0,-3){\line(0,1){6}}
\put(-1.6,0.2){$\overline{\theta_e}$}
\put(0.1,1.2){$\theta_e$}
\put(-2.2,1.6){$\theta$}
\put(2,1.6){$\theta$}
\put(-2.4,-1.9){$\theta$}
\put(2,-1.9){$\theta$}
\put(-0.5,-2.6){$\alpha$}
\put(0.1,-2.6){$\alpha$}
\put(-0.5,2.3){$\alpha$}
\put(0.1,2.3){$\alpha$}
\put(-3.5,-0.4){\line(0,1){0.8}}
\put(-3.4,-0.1){$\beta$}
\put(3.6,0){\line(0,1){0.25}}
\put(2.7,0.1){$\alpha_0$}
\end{picture} 
\caption{Deformations of a regular hyperideal simplex}
\end{figure}
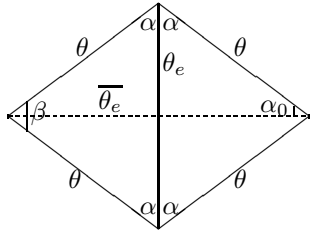
\item use $\alpha$ and $\beta$ to compute the dihedral angles of the
  simplex, using the same triangle formulas but for the links of the
  vertices (which are hyperbolic triangles). There are three angles to
  compute, the angle $\theta_e$ at $e$, the angle $\theta_{\eb}$ at the
  edge opposite to $e$, and the angle $\theta$ at the four other edges. 
\item differentiate the three dihedral angles with respect to $l$, and
  then set $l=l_0$. 
\item find the eigenvalues of the matrix containing the derivatives of
  the dihedral angles with respect to the edge lengths, and check that
  they are all strictly positive.
\end{itemize}
Since this computation is a little tedious, we omit it here and give, in
the appendix A, a little {\it maple} program to do it. 
\epv

\blm \label{lm:concave-1}
For each $i\in \{ 0,1,2,3,4\}$, the volume $V$ is a strictly concave
function on $\cS_i$. 
\elm

\bpv
Let $S\in \cS_i$ be a hyperideal simplex. Suppose that there is a
direction in 
$T_S\cS_i$ which is in the kernel of $\hess(V)$. Then by the Schl{\"a}fli
formula (\ref{eq:schlafli}), the corresponding first-order variation of
the edge lengths vanishes, and this is impossible by proposition
\ref{pr:rig-3}. Therefore, $\hess(V)$ has constant signature over
each  $\cS_i$, with maximal rank. so it only remains to check that
$\hess(V)$ is negative definite at a point. 

By proposition \ref{pr:regulier}, $\hess(V)$ is negative definite at some
regular hyperideal simplex. Therefore, $\hess(V)$ is negative definite
on $\cS_0$, so $V$ is strictly concave on $\cS_0$. 

But $\cS_1$ can be identified with one of the codimension 1 faces of
$\cS_0$; so $V:\cS_1\rightarrow \R_+$, as a limit of strictly concave
functions, is concave. Since we have seen above that its hessian is
non-degenerate, it is strictly concave. In the same way, $\cS_2$ is a
codimension 1 face of $\cS_1$, so $V$ restricted to it is strictly
concave, and the same can then be said of $\cS_3$, and then of $\cS_4$.
\epv

Of course the concavity of $V$ over $\cS_4$ is not new, it is a
special case of a result of \cite{Ri2}, namely the concavity of the
volume of ideal polyhedra. Actually it was pointed out to me by Igor
Rivin that the concavity of $V$ over $\cS_3$ is also a consequence of
\cite{Ri2}, since gluing two copies of the truncated simplex along their
cut leads to an ideal polyhedron.

Note that the first part of the proof can be used in a wider context; it
shows that for any (not necessarily convex) polyhedron in a non-flat
3-dimensional space form, the isometric deformations correspond exactly
to the kernel of the Hessian of the volume.


\section{Hyperideal polyhedra}

This section contains an extension of the result of the previous section
from hyperideal simplices to hyperideal polyhedra. The main result of
this section is the following lemma:

\blm \label{lm:angles-poly}
Let $\sigma$ be a cellulation of $S^2$, and let $w:\sigma_1\rightarrow
(0,\pi)$ be a map on the set of edges of $\sigma$. There exists a
hyperideal polyhedron with combinatorics
given by $\sigma$ and exterior dihedral angles given by $w$ if and only
if:
\begin{itemize}
\item the sum of the values of $w$ on each circuit in $\sigma_1$ is
greater than $2\pi$, and strictly greater if the circuit is
non-elementary.
\item The sum of the values of $w$ on each simple path in $\sigma_1$ is
  strictly larger than $\pi$.
\end{itemize}
This hyperideal polyhedron is then unique. 
\elm

This result was proved by Bao and Bonahon \cite{bao-bonahon} using a
"direct" deformation approach, with the key infinitesimal
rigidity lemma proved using the Legendre-Cauchy method.
Another proof was given recently by Rousset
\cite{rousset1}, who reduced this result to the description of the dual
metrics of compact polyhedra previously achieved by Rivin and Hodgson
\cite{Ri,RH} (following related work of Andreev
\cite{Andreev}). Rousset \cite{rousset1} also extended the description of
the dihedral angles of hyperideal polyhedra to the fuchsian case. 

We give here yet another proof of lemma \ref{lm:angles-poly} --- with
parts in common with the proof of \cite{bao-bonahon}, in particular in
the use of the Andreev theorem at the end. The infinitesimal rigidity is
based here on the remarkable properties of the volume of the hyperideal
simplices and polyhedra. This is similar to what was done for ideal
polyhedra previously (see \cite{CdeV,bragger,Ri2}. The main point here
is that this also works for hyperideal polyhedra, thanks to lemma
\ref{lm:concave-1}. We mostly use here the same notations as in
\cite{ideal}. 

Another important result of this section, which is related to lemma
\ref{lm:angles-poly}, is that the volume is a concave function not only
for simplices, but also for polyhedra.

\blm \label{lm:poly-vol}
Let $\sigma$ be a cellulation of $S^2$, and let $V_i$ be a subset of the
set of vertices of $\sigma$. The volume is a strictly concave
function of the dihedral angles, on the space of hyperideal polyhedra
with combinatorics given by $\sigma$, with ideal vertices at the points
of $V_i$, and strictly hyperideal vertices at the other vertices of
$\sigma$. 
\elm

Note that, to prove lemmas \ref{lm:angles-poly} and \ref{lm:poly-vol},
it is sufficient to prove them when $\sigma$ is a triangulation of
$S^2$, i.e. when all faces of $\sigma$ are triangles. The general result
then follows by adding the constraint that some exterior dihedral angles
--- on the edges which have been added to make $\sigma$ a triangulation
--- are equal to $0$. Given such a triangulation $\sigma$, it is not
difficult to find a triangulation $\Sigma$ of the ball, i.e. a
decomposition of the ball into simplices, such that $\sigma$ is the
"trace" of  $\Sigma$ on the boundary. This can be done for instance by
choosing a vertex $x_0$ of $\sigma$, and adding one simplex with vertex
$x_0$ for each face of $\sigma$ not adjacent to $x_0$. In this section,
$\sigma$ and $\Sigma$ will be fixed.

\bdf 
A {\bf sheared hyperbolic structure} on the ball is a singular
hyperbolic metric on $B^3$, defined by the choice, for each simplex of
$\Sigma$, of a diffeomorphism onto a hyperideal simplex in $H^3$, up to the
isotopies fixing the vertices. The space
of sheared hyperbolic structures is denoted by $\Hsh$.
\edf

Note that a sheared hyperbolic structure does not, in general, define a
hyperbolic structure on $B^3$, or even on the complement of the edges of
the triangulation. The obvious reason is that it is in general not
possible to glue two hyperideal simplices along a face of each --- this
is a difference with the case of ideal triangulations. Actually it is
not difficult to show that, given two hyperideal simplices $S$ and
$S'$ and distinct vertices $x_1, x_2, x_3$ of $S$ and $x'_1, x'_2, x'_3$
of $S'$, it is possible to glue $S$ to $S'$ in a way that identifies
$x_i$ to $x'_i$ (for $1\leq i\leq 3$) if and only if all the edges of $S$
with endpoints in $\{ x_1, x_2, x_3\}$ either have infinite length, or
have the same length as the corresponding edge of $S'$, and conversely. 

For each simplex of $\Sigma$, the possible hyperbolic metrics on it are
parametrized by the dihedral angles; this determines an affine structure
on $\Hsh$.

\bdf
Let $\Theta:\Sigma_1\rightarrow \R_+$ be a map on the edges of $\Sigma$,
which takes values in $(0,\pi)$ on the edges of $\sigma$. We call 
$\Hsh(\Theta)$ the subspace of $\Hsh$ of sheared hyperbolic structures
on $\Sigma$ such that the exterior dihedral angle at each boundary
edge $e$ is $\Theta(e)$, while the total angle around each interior
edge $e'$ is $\Theta(e')$.
\edf

Let $h\in \Hsh$ be a sheared hyperbolic structure on $\Sigma$. Suppose
that it is possible to glue the simplices of $\Sigma$, i.e. the
condition stated above on the length of the edges is satisfied. 
Clearly,
the singularities of $h$ are concentrated on the interior edges of
$\Sigma$. Those singularities have two parts, which can be expressed in
terms of the holonomy $\rho(e)$ of the developing map at an edge $e$:
\begin{itemize}
\item the total angle around $e$ can be different from $2\pi$; this can
  be expressed in terms of the component of $\rho(e)$ around $e$.
\item for each edge having as endpoints two ideal vertices, 
  there might be a translation component of $\rho(e)$ along $e$;
  this means that, when one "makes one turn around $e$", one ends up
  some point away from the starting point. So one can associate to each
  interior edge of $\Sigma$ a number, corresponding to this translation
  length, which we call the {\bf shear} of
  $h$ at $e$. Note that it does not depend on the orientation chosen for
  $e$. 
\end{itemize}

\bdf
An {\bf exact hyperbolic structure} on $\Sigma$ is a sheared
hyperbolic structure such that:
\begin{itemize}
\item the hyperideal simplices can be glued along their common faces. 
\item the shear at all interior edges of $\Sigma$ vanishes. 
\end{itemize}
The space of exact hyperbolic structures is
denoted by $\Hex$. We call $\Hex(\Theta):=\Hex\cap\Hsh(\Theta)$. 
\edf

The "smooth" hyperbolic structures on $\Sigma$ are, by definition, the
exact hyperbolic structures such that the total angle around each
interior vertex of $\Sigma$ is $2\pi$. Thus it is an affine submanifold
of the space of exact hyperbolic structures on $\Sigma$.

A key point is that the exact hyperbolic structures are exactly the 
critical points of the volume, seen as a function on $\Hsh(\Theta)$. 

\bprop \label{pr:holonomy}
Let $h\in \Hsh$ be a sheared hyperbolic structure, with boundary
dihedral angles given by $\Theta$. Then $h$ is in $\Hex$
if and only if $h$ is a critical point of $V$ restricted to
$\Hsh(\Theta)$. 
\eprop

We will only sketch the proof here, and we refer the reader to
\cite{ideal} for the details. The proof 
was done there for ideal simplices, but applies just as well for
hyperideal simplices. The only difference is that in the present case
there are less constaints on the dihedral angles of the simplices, so
that, for instance, if all simplices were supposed to be strictly
hyperideal, the proof would be much simpler.

\bpv[Sketch of the proof]
Suppose that $h\in \Hex$. Let $v_1, \cdots, v_n$ be the ideal vertices
of $\Sigma$. Since the shear of $h$ at all interior edges of $\Sigma$
vanishes, one can choose for each $i$ and each simplex $s$ having $v_i$
as one of its vertex a (part of) horosphere $H$ in $s$ centered at
$v_i$, in a way such that, for the different choices of $s$, the
horospheres coincide on the codimension 1 faces of $\Sigma$. 

One can then apply the Schl{\"a}fli formula of lemma \ref{lm:schlafli-2} to
check that, in any deformation of $h$ which does not change the total angles
at the vertices of $\Sigma$, the volume remains constant (at first
order).

Conversely, suppose that $h$ is a critical point of $V$ restricted to
$\Hsh(\Theta)$. We have to prove that $h\in Hex$, i.e. that corresponding
faces of two simplices can be glued, and that the shear at all interior edges
vanishes.

Let $e$ be an edge of $\Sigma$, with endpoints two strictly hyperideal
vertices; suppose that $e$ is in two simplices $S_1, S_2$ for which its length
is different. Consider the deformation of the angles of $S_1$ and $S_2$ which
increases the angle of $S_1$ at $e$ at speed $1$, decreases the angle of $S_2$
at $e$ at speed one, and does not change any other angle. This deformation is
clearly in the tangent space $T\Hsh(\Theta)$, but by the Schl\"afli formula it
changes the volume, a contradiction. So the faces can be glued. 

Now let $e$ be an interior 
edge of $\Sigma$ with endpoints two ideal vertices.
We call $e_-$ and $e_+$ the endpoints of $e$, and use a special type of
deformation of the simplices which are adjacent to $e$. To describe
those deformations, we call $s$ one of those simplices adjacent to $e$,
and $t_-$ and $t_+$ the 2-faces of $s$ which are adjacent to $e_-$ and
$e_+$, respectively, but do not contain $e$. We orient the vertices of
those triangles in a way compatible with the orientation of $s$. The
deformation is as follows.
\begin{itemize}
\item the angle at $e$ and at the opposite edge do not vary.
\item the angle at the edge of $t_+$ "before" $e_+$, and at the edge of
  $t_-$ "after" $e_-$, varies at speed $+1$.
\item the angle at the edge of $t_+$ "after" $e_+$, and at the edge of
  $t_-$ "before" $e_-$, varies at speed $-1$.
\end{itemize}
The same description applies to all the simplices containing $e$.
This deformation is compatible with the angle conditions on the
simplices even when all vertices are ideal. Moreover, a direct
computation using the Schl{\"a}fli formula and 
an adequate choice of horospheres in the simplices
with ideal vertices (as in \cite{ideal}) shows that the first-order
variation of the volume is proportional to the shear of $h$ at $e$, and
that it vanishes if and only if the shear of $h$ at $e$ vanishes. This
proves the proposition.
\epv

This in turns implies a rigidity result for hyperideal polyhedra, with
respect to their dihedral angles. It is one of the basic tools in the
proof of lemma \ref{lm:angles-poly}. This rigidity result can also be
obtained  by other methods, for instance it was proved in \cite{bao-bonahon}
using the so-called Cauchy method, which was developped by Cauchy
\cite{Cauchy} and Legendre \cite{legendre}\footnote{The fact that it was
  mostly due to 
  Legendre was recently discovered by I. Sabitov.} to prove the global
rigidity of polyhedra in $\R^3$.  
But the proof given here has the
advantage of extending from hyperideal polyhedra 
to manifolds with hyperideal boundary.

\bcr \label{cr:rig-poly}
Let $P$ be a hyperideal hyperbolic polyhedron. There is no non-trivial
infinitesimal deformation of $P$ which changes neither its combinatorics
nor its dihedral angles. 
\ecr

\bpv
Let $h$ be the exact hyperbolic structure on $\Sigma$ corresponding
to $P$. An infinitesimal deformation of $P$ which does not change its
dihedral angles would be equivalent to a first-order deformation of $h$,
in $\Hsh(\Theta)$ (where $\Theta$ is given by the dihedral angles of $P$
on the boundary edges of $\Sigma$, and is equal to $2\pi$ on the
interior edges) such that the volume remains critical at first
order. This would contradict the strict concavity of the volume.
\epv

\medskip

The proof of lemma \ref{lm:angles-poly} also requires a compactness
result, which we now state. 

\bprop \label{pr:compact-poly}
Let $(P_n)_{n\in \N}$ be a sequence of hyperideal polyhedra with the
same combinatorics. Suppose that the sequence of dihedral angles
$(\theta_n)_{n\in \N}$ converges to a limit $\theta$. Then:
\begin{itemize}
\item either some subsequence of $(P_n)$ converges.
\item or there exists a non-elementary circuit in $\sigma_1$ on which
  the sum of the values of $\theta$ is $2\pi$.
\item or there exists a simple path in $\sigma_1$ on which
  the sum of the values of $\theta$ is $\pi$.
\end{itemize}
\eprop

We will skip the proof here, since this proposition is a consequence of
a compactness lemma in \cite{shu} (which is stated in a more general
setting). Note 
however that it is also a special case of a compactness result for
manifolds with hyperideal boundary, lemma \ref{lm:compact} below. 

Corollary \ref{cr:rig-poly} and proposition \ref{pr:compact-poly} show
that, for each combinatorics, the map sending a hyperideal polyhedron to
its dihedral angles is a covering. To prove that it is one-to-one, we
will show that there exists a specific set of dihedral angles which has
a unique inverse image. 

\bprop \label{pr:connect-poly}
Let $\sigma$ be a cellulation of $S^2$, along which a subset $V_i$ of
the vertices of $\sigma$. 
There exists at least one angle
assignation satisfying the hypothesis of lemma \ref{lm:angles-poly}
which is realized as the dihedral angles of a unique hyperideal
polyhedron, the ideal vertices of which are the elements of $V_i$.
\eprop

Of course the uniqueness here is again up to the global hyperbolic
isometries. 

\bpv
We use the same proof as the one given by Bao and Bonahon
\cite{bao-bonahon}; so we only sketch the proof and refer the reader to
\cite{bao-bonahon} for more details. We will use a result on compact
polyhedra (as in \cite{Andreev,RH}) to obtain a compact polyhedron
which is the truncated version of a hyperideal polyhedron with the right
combinatorics. 

\begin{figure}[h] \label{fig:transfo}
\centerline{\psfig{figure=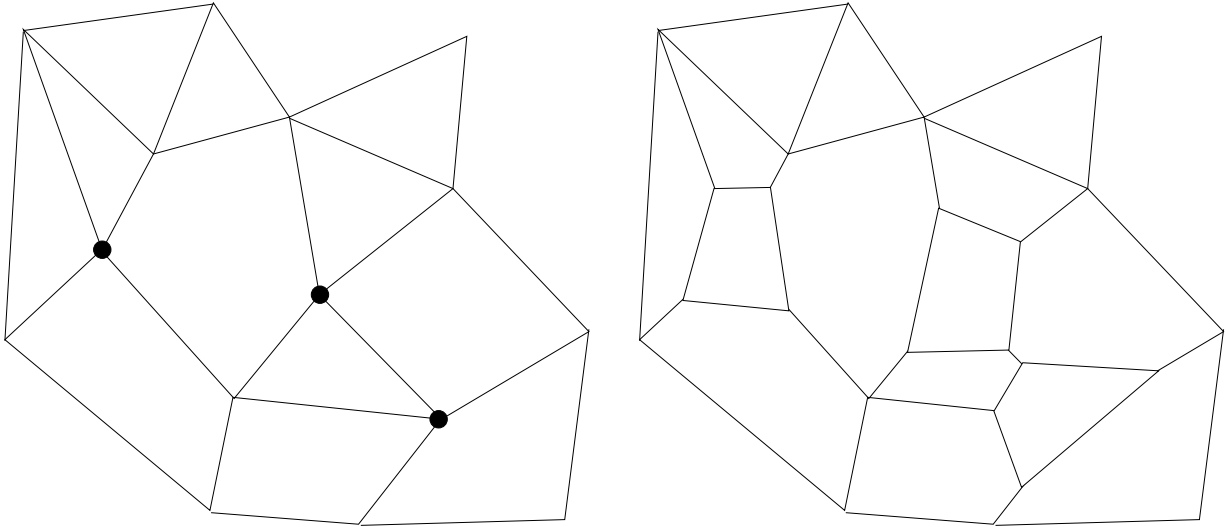,height=7cm}}
\caption{Truncating a cellulation with hyperideal vertices}
\end{figure}

We associate to $\sigma$ and $V_i$ another
cellulation, $\sigmab$, defined as follows. 
\begin{itemize}
\item $\sigmab$ has one face for each face of $\sigma$, and one for each
 vertex of $\sigma$ which is not in $V_i$ (those vertices correspond to
 the strictly hyperideal vertices of $\sigma$).
\item it has one vertex for each element of $V_i$, and one for each
  couple $(e,v)$, where $v$ is a vertex of $\sigma$ which is not in
  $V_i$ and $e$ is an edge containing $v$.
\item it has one edge for each edge of $\sigma$, and one for each couple
  $(f,v)$, where $v$ is a vertex of $\sigma$ not in $V_i$ and $f$ is a
  face of $\sigma$ containing $v$. 
\end{itemize}
Note that this transformation is the same as the transformation sending
the combinatorics of a hyperideal polyhedron to the combinatorics of its
trucated polyhedron. In figure \ref{fig:transfo} the hyperideal
vertices are represented as "bigger" dots. 
Now we can put a weight $2\pi/3$ on each edge of $\sigmab$ which is also
an edge of $\sigma$, and a weight $\pi/2$ on each edge of $\sigmab$
which is not in $\sigma$. Then apply the Andreev
theorem \cite{Andreev,RH} to obtain that
there is a unique compact polyhedron $P_i$ with combinatorics given by
$\sigmab$ and with those dihedral angles. We leave it to the reader to
check that the hypothesis of the Andreev theorem apply, as in
\cite{bao-bonahon}. 
\epv

\bpn{of lemma \ref{lm:angles-poly}}
First note that it is sufficient to prove the result when all faces
of $\sigma$ are triangles. Indeed, if $\sigma$ has some non-triangular
faces, one can add some edges so as to obtain triangles. By restricting
the possible deformations to those which do not change the angles at
those additional edges, the lemma as it is stated will  follow from the
result when the faces are triangles. 

Consider a triangulation $\sigma$ of $S^2$, with a subset $V_i$ of the set
of vertices of $\sigma$. Let $q:=|V_i|$, and let $p$ be the number of
vertices of $\sigma$ which are not in $V_i$. We call:
\begin{itemize}
\item $\cP_{\sigma, V_i}$ the space of hyperideal polyhedra with
  combinatorics given by $\sigma$, with ideal vertices exactly the
  elements of $V_i$.
\item $\cA_{\sigma, V_i}$ the space of angle assignations on the edges
  of $\sigma$, satisfying the hypothesis of lemma \ref{lm:angles-poly},
  such that, for each vertex $v$ of $\sigma$, the sum of the angles on
  the edges adjacent to $v$ is $2\pi$ if and only if $v\in V_i$. 
\item $F:\cP_{\sigma, V_i}\rightarrow \cA_{\sigma, V_i}$ the map sending
  a hyperideal polyhedron to the set of its dihedral angles. 
\end{itemize}

Both $\cP_{\sigma, V_i}$ and $\cA_{\sigma, V_i}$ are locally smooth
manifolds of dimension $3p+2q-6$. For $\cP_{\sigma, V_i}$ it is clear,
because $\cP_{\sigma, V_i}$ is parametrized by the position of the $p$
hyperideal vertices and the $q$ ideal vertices, modulo the action of the
6-dimensional group of isometries of $H^3$. $\cA_{\sigma, V_i}$ is
parametrized by the dihedral angles on the edges of $\sigma$, subject to
the condition that the sum of the angles at each ideal vertex is
$2\pi$. But a simple computation using the Euler formula shows that the
number of edges of $\sigma$ is $3p+3q-6$. So we only have to prove that
the conditions at the ideal vertices are independent; or, in other
terms, that if $a_1, \cdots, a_p\in \R$ are coefficients associated to
the ideal vertices such that, for any edge $e$, the sum of the coefficients
of the ideal vertices in $e$ is $0$, then all $a_i$ are zero. Then, for each
face $f$ of $\sigma$: 
\begin{itemize}
\item either $f$ has at least one strictly hyperideal vertex, and then
  the coefficients of the two other vertices of $f$ are $0$.
\item or all the vertices are ideal, and we use the fact that all faces
  of $\sigma$ are triangles, so an elementary parity argument shows
  again that all coefficients are $0$.
\end{itemize}
This ends the proof that $\dim(\cA_{\sigma, V_i})=\dim(\cP_{\sigma,
  V_i})$. 

Corollary \ref{cr:rig-poly} shows that $F$ is a local
homeomorphism. Moreover, proposition \ref{pr:compact-poly} shows that
$F$ is proper, so that it is a covering. But $\cA_{\sigma, V_i}$ is the
interior of a convex polytope in some $\R^N$ for some $N$, so this
covering is a homeomorphism between each connected component of
$\cP_{\sigma,V_i}$ and $\cA_{\sigma, V_i}$. 
Finally proposition \ref{pr:connect-poly} shows
that one element of $\cA_{\sigma, V_i}$ has a unique inverse image, and
therefore $F$ is a homeomorphism.
\epn

Finally, lemma \ref{lm:poly-vol} is a direct consequence of the
construction which we have made, and of the following elementary remark,
which we have taken from \cite{ideal}.

\brk \label{rk:concavity}
Let $\Omega\in \R^N$ be a convex subset, and let $f:\Omega\rightarrow
\R$ be a smooth, strictly concave function. Let
$\rho:\R^N\rightarrow \R^p$ be a 
linear map, with $p<N$, and let $\Omegab:=\rho(\Omega)$. Define a
function:
$$
\begin{array}{cccc}
\fb: & \Omegab & \rightarrow & \R \\
& y & \mapsto & \max_{x\in \rho^{-1}(y)} f(x) 
\end{array}
$$
Then $\Omegab$ is convex, and $\fb$ is a smooth, strictly
concave function on $\Omegab$.
\erk

An important point is that this strict concavity extends to the case
where some exterior dihedral angles are zero, i.e. when some edges are
just segments drawn in a face.

\pg{From polyhedra to hyperideal manifolds}

We will show below that the proof just given for lemma
\ref{lm:angles-poly} applies, with some modifications, to the proof of
theorem \ref{tm:angles}. Just as we have used a decomposition of
hyperideal polyhedra into hyperideal simplices, the "building blocks"
for hyperideal manifolds will be hyperideal polyhedra. It would be more
natural {\it a priori} to use hyperideal simplices here too, but it
turns out to be a problem to prove that any hyperideal manifold admits
a non-degenerate triangulation by hyperideal simplices. Another solution
is used in \cite{ideal} for ideal manifolds --- the point there is that
any ideal manifold has a finite cover which admits an ideal
triangulation. 

\bdf \label{df:triang}
Let $M$ be a hyperideal hyperbolic manifold. A {\bf hyperideal
  cellulation} of $M$ is a decomposition of $M$ as the union of
hyperideal polyhedra with non-empty interior, isometrically glued along
their faces, such that 
the intersection of two polyhedra is always a face of each of them.
\edf

This definition implies a non-degeneracy condition: we do not allow
degenerate hyperideal polyhedra, or two polyhedra having in common a
triangle which is not a face of one of them. 
 
\blm \label{lm:triang}
Each hyperideal manifold admits a hyperideal cellulation.
\elm

\bpv
It is done along the ideas of Epstein and Penner \cite{epstein-penner};
the situation here is simpler since the action of $\pi_1M$ on
$S^2\setminus \Lambda$ is discrete. On the other hand, one has to take
into account the hyperideal vertices. Another, slightly more complicated
but maybe more geometric version of the same idea, was used in
\cite{ideal}. 

Let $M$ be a hyperideal hyperbolic manifold. Then $M$ is isometric to the
convex hull of a set $\{ x_1, \cdots, x_p\}$ of ideal points and a set
$\{ y_1, \cdots, y_q\}$ of hyperideal
points in $E(M)$, where $E(M)$ is the extension of $M$. 
By "hyperideal point" in $E(M)$, we mean here
the orbit of a point in $S^3_1$ under $\pi_1(E(M))$, seen as a group
acting by isometries on $H^3$ and thus also on its "extension" $S^3_1$. 

For each $i\in \{ 1,2,\cdots,p\}$, choose a ``small'' horoball
$b_i\subset E(M)$ centered at the ideal point $x_i$; we suppose that the
$b_i$ are small enough to be embedded and disjoint. Let $B_i$ be the set
of inverse images of $b_i$ under the quotient of $H^3$ by $\pi_1E(M)$. 

Let $b$ be one of the horoballs in the collection $B_i$. In the
Minkowski model of $H^3$, $b$ is the intersection with $H^3$ of an
affine light-like hyperplane $h$. There exists a unique vector $u(b)\in
\R^4_1$ such that:
$$ h:=\{ x\in \R^4_1 ~| ~ \langle x, u(b)\rangle =1\}~. $$
Note that, in the projective model of $H^3$, $u(b)$ projects to the
point at infinity of $b$. Let:
$$ \overline{B_i} := \{ u(b) ~ | ~ b\in B_i\}~. $$

For each $j\in \{ 1, 2, \cdots, q\}$, the hyperideal vertex $y_j$ of $M$
lifts to an equivariant set $C_j$ of hyperideal points in $H^3$, or, in
other terms, to an equivariant set of points in $S^3_1$. In the
Minkowski model of $H^3$ we also consider the de Sitter space as the
submanifold: 
$$ S^3_1 = \{x\in \R^4_1 ~ | ~ \langle x,x\rangle =1\}~, $$
with the induced metric. Again note that $C_j$ projects, in the natural
extension of the projective model described above to $S^3_1$, to a lift 
of $y_j$.

Let $C:=(\cup_{i=1}^p \Bb_i) \cup(\cup_{j=1}^q C_j)$. It is a discrete
set of $\R^4_1$ which is invariant under the action of $\pi_1M$, and the
projections on $\{ x_1=1\}$ (in the direction of $0$) of its points are
exactly the vertices of $\Mt$. 

Let $CH$ be the convex hull of $C$ in $\R^4_1$. The radial projection on $\{
x_1=1\}$ of $CH$ is the convex hull of the radial projections of the
elements of $C$, so its intersection with the radius $1$ ball is exactly
$\Mt$. Moreover, $CH$ is polyhedral, and its boundary has a
decomposition into 3-dimensional polyhedra. Since the radial projection
sends the 
geodesics of $\R^4_1$ to lines in $\{ x_1=1\}$, this decomposition
projects to $\Mt$ as a cellulation $\Sigma$ of $\Mt$, i.e. a
decomposition of $M$ into 3-dimensional polyhedra. By construction
$\Sigma$ is invariant under the action of $\pi_1M$, so $\Sigma$
determines a cellulation of $M$, which we also call $\Sigma$. 

Also by construction, for each $h>0$, there is a finite set of points of $C$
with first coordinate $x_1\leq h$. So there is also a finite set of
polyhedra in $\dr CH$ which contain a point with $x_1\leq h$.
If $K\subset \Mt$ is a compact
subset, its radial projection on $\dr CH$ has first coordinate bounded
by some $h>0$, so $\Sigma$ is locally finite ---
each compact subset of $\Mt$ intersects only a finite number of
polyhedra. 

This implies that each polyhedron has only a finite number of
vertices. Otherwise, one of the polyhedra, say $P$, having an infinite set of
vertices, would intersect an infinite set of disjoint fundamental
domains of $M$, and it would follow that each fundamental domain of $M$
intersects an infinite set of images of $P$ under elements of $\pi_1M$.
\epv

Note that this construction is far from clear if there is a closed curve in
$\dr M$ which is a geodesic of $M$, because in this case the endpoints of that
geodesic are in the limit set $\Lambda$.


\section{Rigidity}

\pg{An outline}

We can now proceed to lemma \ref{lm:rigidity}, which is
the main step in the proof of theorem \ref{tm:angles}. We recall the
statement for the reader's convenience.

\blm \label{lm:rigidity}
Let $M$ be a hyperideal manifold. Any first-order deformation of its
dihedral angles is obtained by a unique first-order deformation of $M$. 
\elm

To prove it we will consider a hyperideal manifold $M$. By lemma
\ref{lm:triang}, it admits a decomposition into hyperideal polyhedra. 
Moreover, if those polyhedra have non-triangular faces, we can further
subdivide those faces, so as to obtain only polyhedra with triangular faces
--- 
they will then have some (exterior) dihedral angles equal to $0$, but it
won't make any difference in the proof. 

By lemma \ref{lm:angles-poly}, the dihedral angles of the hyperideal
polyhedra provide a parametrization of their deformations. However,
since some of the faces are hyperideal triangles, it is not always
possible to glue them to obtain a hyperbolic structure --- even singular
--- on $M$; there are compatibility conditions related to the lengths of
the edges. Of course, if such a gluing is possible, the hyperbolic
structure obtained in this way will 
in general be singular at the interior edges of $\Sigma$, where the
total angle will in general be different from $2\pi$ and where a "shear"
might occur along an edge having two ideal endpoints. 

But, just as in section 4, the condition that the total angle is $2\pi$
is an affine condition on the space of dihedral angles of the polyhedra,
while the lengths are compatible and the shear vanishes at all interior edges
of $\Sigma$ if and 
only if the volume, restricted to the deformation which do not change
the total angles at the edges of $\Sigma$, is critical. Since the volume
of each polyhedron is strictly concave, the volume of $M$ (as a function
of the dihedral angles of the polyhedra) is also strictly concave, and
this will allow the same argument as the one given in the proof of
corollary \ref{cr:rig-poly}. The existence of deformations inducing a
given first-order variation of the dihedral angles will follow from this
rigidity statement and a dimension-counting argument.

\pg{The space of angle assignations}

From now on, we consider a hyperideal manifold $M$, along with a
decomposition $\Sigma$ of $M$ into a finite set of polyhedra, for
instance as provided 
by lemma \ref{lm:triang} if $M$ is supposed to be a hyperideal manifold. 
We also suppose given a subset $V_i$ of the set of vertices of
$\Sigma$. $V_i$ will appear later as the set of vertices of $\Sigma$
which correspond to ideal points for the hyperbolic metrics that will
appear on $M$.

\bdf \label{df:angles}
An {\bf angle assignation on $\Sigma$} is the choice, for each
polyhedron $p$ of
$\Sigma$, of a set of (exterior) dihedral angles on the edges, subject
to the condition that:
\begin{itemize}
\item it satisfies the hypothesis of lemma \ref{lm:angles-poly}.
\item for each vertex $v$ of 
  $p$ which is in $V_i$, the sum of the angles assigned to the edges of $p$
  containing $v$ is equal to $2\pi$. 
\item for each vertex $v$ of $p$ which is not in $V_i$, the sum of the
  angles assigned to the edges of $p$ containing $v$ is strictly larger
  than $2\pi$. 
\end{itemize}
The
set of angle assignations on $\Sigma$ will be denoted by $\cA_\Sigma$.

An {\bf angle assignation on $\dr \Sigma$} is the choice, for each exterior
edge of $\Sigma$, of an angle in $(0, \pi)$, in a way that satisfies the
conditions in theorem \ref{tm:angles}. The set of angle assignations on
$\dr \Sigma$ will be denoted by $\cA_{\dr \Sigma}$.
\edf

By lemma \ref{lm:angles-poly}, each angle assignation on $\Sigma$
defines, for each polyhedron $p$ of $\Sigma$, a homeomorphism from $p$
to a hyperideal polyhedron (up to isotopies fixing the vertices), 
which we can also consider as a "hyperideal
metric" on $p$. Each vertex of $p$ is then either ideal --- if the sum
of the angles on the adjacent edges is $2\pi$, i.e. if $v\in V_i$ --- or
strictly  hyperideal. 

For some angle assignations, the hyperideal polyhedra associated to the
polyhedra of $\Sigma$ can be glued along their faces. This happens
exactly when, for each edge $e$ of $\Sigma$ with two endpoints which are not
in $V_i$, the lengths of $e$ for the hyperideal metrics on the
polyhedra of $\Sigma$ containing $e$ are the same. 
The next definition describes assignations such that the total angle
around each interior edge of $\Sigma$ is as needed.

\bdf 
Let $\theta\in \cA_\Sigma$ be an angle assignation on $\Sigma$. $\theta$ is:
\begin{itemize}
\item {\bf coherent}, if the lengths assigned to the edges with no
  endpoint in $V_i$ by all the polyhedra 
  that contains it are the same. The set of coherent angle
  assignations will be denoted by $\cA_\Sigma^c$.
\item {\bf regular} if the total angle around each interior edge is
  $2\pi$. The set of regular angle assignations will be denoted by
  $\cA_\Sigma^r$. 
\end{itemize}
\edf

Note that the set of regular angle assignations is an affine subspace
of $\cA_\Sigma$. When $\theta$ is coherent, it is possible to glue
isometrically the faces of the polyhedra (according to the combinatorics given
by $\Sigma$).

\bdf 
An angle assignation $\theta$ on $\Sigma$ is {\bf exact} if it is
coherent and, in addition, the singular 
hyperbolic structure obtained by gluing the polyhedra has zero shear
at all the interior edges of $\Sigma$. The space of exact angles
assignations on $\Sigma$ will be denoted by $\cA_\Sigma^e$.
\edf

The condition that the shear is non-zero is non-trivial only for the
edges with two endpoints which are in $V_i$. 

The smooth hyperbolic structures on $M$ correspond to the angle
assignations which are both exact and regular; in that case, at each
interior edge of $\Sigma$, the shear vanishes because the structure is
exact, while the total angle is $2\pi$ since the structure is regular.

\pg{The volume}

As a consequence of proposition \ref{lm:concave-1}, the volume is
concave over the space of angle assignations on $M$.

\bprop \label{pr:concave-2}
$V$ is a strictly concave function on $\cA_\Sigma$, and also on
$\cA_\Sigma^r$. 
\eprop

\bpv
The concavity on $\cA_\Sigma$ is clear since $V$ is the sum of the
volumes of the simplices, which are concave functions. The concavity
over $\cA_\Sigma^r$ is a consequence since $\cA_\Sigma^r$ is an affine
subset of $\cA_\Sigma$. 
\epv

\bdf
Let $\alpha:\sigma_1\rightarrow (0, \pi)$ be an angle assignation on
$\dr \Sigma$. We call $\cA_\Sigma(\alpha)$ the set of angle
assignations on $M$ such that the total (interior) angle at each
boundary edge $e$ is $\pi-\alpha(e)$, and $\cA^r_\Sigma(\alpha)$
(resp. $\cA^c_\Sigma(\alpha)$) the space of those angle assignations
which are regular (resp. coherent).
\edf

\bprop \label{pr:critical}
Let $\theta\in \cA^r_\Sigma$, and let $\alpha:\sigma_1\rightarrow (0, \pi)$ be
the function sending a boundary edge to its exterior dihedral
angle. Then $\theta$ is exact if and only if it is a critical point of $V$
restricted to $\cA^r_\Sigma(\alpha)$. 
\eprop

\bpv
Suppose first that $\theta$ is exact. By construction, $\theta$
determines a hyperbolic metric on each polyhedron of $\Sigma$. $\theta$
is coherent, so that, by definition, the hyperideal polyhedra obtained
can be glued along their faces. Since $\theta$ is exact and regular, the
resulting singular hyperbolic structure on $M$ is actually
smooth. Choose a horosphere centered at each of the ideal vertices of
$\Sigma$, and let $\thetad\in T_\theta \cA^r_\Sigma(\alpha)$ be an
infinitesimal variation of $\theta$.
The Schl{\"a}fli formula shows that the total first-order variation of the
volume vanishes. This shows that $\theta$ is a critical point of $V$
restricted to $\cA^r_\Sigma(\alpha)$.

\medskip

Suppose now that $\theta$ is not exact. Then either $\theta$ is not
coherent, or it is coherent but not exact. 

Suppose first that $\theta$ is not coherent. Then there is an edge $e$
of $\Sigma$ which is contained in two polyhedra $p_1$ and $p_2$ of
$\Sigma$ which, for the angle assignations determined by
$\theta$, give different lengths to $e$. This implies that both
endpoints of $e$ are strictly hyperideal. Consider the first-order
deformation $\thetad$ of $\theta$ which 
\begin{itemize}
\item increases the angle of $p_1$ at $e$ at speed $1$.
\item decreases the angle of $p_2$ at $e$ at speed $1$.
\item does not change any other angle. 
\end{itemize}
A short check shows that $\thetad\in T_\theta
\cA^r_\Sigma(\alpha)$. Moreover, the Schl{\"a}fli formula shows that
$dV(\thetad)\neq 0$, so $\theta$ is not a critical point of $V$
restricted to $\cA^r_\Sigma(\alpha)$.

Suppose now that $\theta$ is coherent but not exact. There is then an
interior edge $e$ of $\Sigma$ at which the shear of the singular
hyperbolic structure defined by $\theta$ is not zero. We call $e_-$ and
$e_+$ the endpoints of $e$, which are both ideal.
Let $p_1, \cdots, p_r$ be the polyhedra of $\Sigma$ containing $e$, in
cyclic order; we set $p_0:=p_r, p_{r+1}=p_1$. 
For each $i\in \{1, \cdots, r\}$, let $e^+_i$ be the edge of $\Sigma$
which is common to $p_{i-1}$ and $p_i$ and contains $e_+$, and let
$e^-_i$ be the edge of $\Sigma$ 
which is common to $p_{i-1}$ and $p_i$ and contains $e_-$. For each
$i\in \{1, \cdots, r\}$, let $f_i$ be the 2-face of $\Sigma$ containing
$e, e^+_i$ and $e^-_i$.

For each couple $(p, s)$, where $p$ is a polyhedron of $\Sigma$ and $s$
is an ideal vertex of $\Sigma$ contained in $p$, choose a horosphere in
$p$ centered at $s$, in such a way that:
\begin{itemize}
\item for all $i\in \{ 2, \cdots, r\}$, the intersections with $f_i$ of
  the horospheres centered at $e_-$ and $e_+$ on both sides of $f_i$
  coincide. 
\item for all $i\in \{1, \cdots, r\}$, the intersections with $f_i$ of
  the horospheres centered at the endpoints of $e^+_i$ and $e^-_i$
  opposite to $e_+$ and $e_-$ coincide. 
\end{itemize}
Since the shear of $\theta$ at $e$ does not vanish, the horospheres
centered at $e_+$ (resp. $e_-$) on both sides of $f_1$ do not coincide:
they are at a constant distance equal to the shear of $\theta$ at $e$.
Consider the first-order variation $\thetad$ of $\theta$ which, for each
$i\in \{1, \cdots, r\}$: 
\begin{itemize}
\item increases at speed one the angles of $p_i$ at $e^+_i$ and at
  $e^-_{i+1}$. 
\item decreases at speed one the angles of $p_i$ at $e^+_{i+1}$ and at
  $e^-_i$. 
\item does not change any other angle. 
\end{itemize}
This deformation is in $T_\theta \cA_\Sigma$, i.e. it respects the
condition that, for each polyhedron $p$ of $\Sigma$, the sum of the
dihedral angles at each vertex of $p$ which is in $V_i$ remains
$2\pi$. It is also in $T_\theta \cA_\Sigma (\alpha)$, because the
dihedral angle on the boundary edges of $\Sigma$ do not change. It is
even in $T_\theta \cA_\Sigma^r (\alpha)$, because the total angle around
the interior edges of $\Sigma$ do not change --- the variations
corresponding to $p_i$ and $p_{i+1}$ always cancel.

Now apply the Schl{\"a}fli formula (lemma \ref{lm:schlafli-2}) to $\thetad$ with
this choice of horosphere. There are two contributions for each edge
$e^+_i$ and $e^-_i$, one for each side of $f_i$. But it is clear that
those contributions cancel except for $e^+_1$ and $e^-_1$, where the
fact that the horospheres centered at $e_-$ and $e_+$ on both side do
not coincide on $f_1$ means that there is a discrepancy which is
proportional to the shear at $e$. This shows that, when $\theta$ is not
exact, it is not a critical point of $V$ restricted to
$\cA^r_\Sigma(\alpha)$. 
\epv

\bpn{of lemma \ref{lm:rigidity}}
We start from a hyperideal manifold $M$, with a decomposition into
hyperideal polyhedra, as given by lemma \ref{lm:triang}. 

The dihedral angles of the polyhedra of $\Sigma$, for the hyperbolic
metric on $M$, define an element $\theta$ of $\cA^r_\Sigma\cap \cA^c_\Sigma$,
and, by proposition \ref{pr:critical}, it is a critical point of $V$
restricted to $\cA^r_\Sigma(\alpha)$, where $\alpha$ is the function
sending a boundary edge of $\Sigma$ to its exterior dihedral angle.

Let $\alphad$ be a first-order deformation of $\alpha$, and let
$(\alpha_t)_{t\in [0,1]}$ be a 1-parameter deformation of $\alpha$
with $\alphad=d\alpha_t/dt$. For $t$ small enough,
$\cA^r_\Sigma(\alpha_t)$ is an affine subspace of $\cA^r_\Sigma$ close
to $\cA^r_\Sigma(\alpha)$, and the strict concavity of $V$ shows that
there is a unique maximum $\theta_t$ of $V$ restricted to
$\cA^r_\Sigma(\alpha_t)$. By proposition \ref{pr:critical}, $\theta_t$
determines a hyperideal structure on $M$.

Again by the strict concavity of $V$, $(\theta_t)$ is a smooth 1-parameter
deformation of $\theta$; if $\thetad:=d\theta_t/dt$, $\thetad$ is a
first-order deformation of $\theta$ such that the induced variation of
the boundary dihedral angles is $\alphad$. 
\epn

\brk \label{rk:concave}
$V$ is a strictly concave function on $\cA^r_\Sigma\cap \cA^c_\Sigma$,
parametrized by the boundary dihedral angles. 
\erk

\bpv
As in section 4, the proof is a direct consequence of remark
\ref{rk:concavity}. 
\epv


\section{Compactness}

This section contains the proof of the basic compactness result which we
need for the proof of theorem \ref{tm:angles}; we recall it first. 

\blm \label{lm:compact}
Let $(g_n)_{n\in \N}$ be a sequence of hyperideal structures on $M$,
with the same  boundary combinatorics. For each $n$, let $\alpha_n$ be
the function which associates to each boundary edge of $(M, g_n)$ its
exterior dihedral angle, and suppose that $\alpha_n\rightarrow \alpha$,
where $\alpha$ still satisfies the hypothesis of theorem
\ref{tm:angles}. Then, after taking a subsequence, $g_n$ converges to
a hyperideal structure $g$ on $M$.
\elm

This result uses the following natural notion of convergence of
manifolds with hyperideal
boundary. Let $(g_n)_{n\in \N}$ be a sequence of hyperbolic metrics on
$M$, such that the $(M, g_n)$ are manifolds with hyperideal
boundary.  
Let $(M, g)$ be a hyperbolic manifold with hyperideal
boundary. 
We say that the $(M, g_n)$ {\bf converges} to $(M, g)$ if, for each
compact subset $K\subset E(M,g)$ in the extension of $E(M,g)$ and 
each $\epsilon>0$, for each $n$ large enough, there exists a compact
subset $K_n\subset E(M, g_n)$ such that $K_n\cap (M, g_n)$ is 
$\epsilon$-close to $K\cap (M, g)$ in the Gromov-Hausdorff distance. 

The proof uses another compactness lemma, concerning sequences of
hyperbolic manifolds with a convex, polyhedral boundary (which is not
hyperideal). It states that if the third fundamental forms of the
boundary converge to a reasonable limit, then the sequence converges. 
Lemma \ref{lm:compact} will follow by truncating the
hyperideal ends, so as to obtain a sequence of manifolds with
polyhedral, non hyperideal boundary. It will be necessary later for some
applications to circle packings. 

We will also give a slightly more general compactness result, in which
some geodesics for the limit of the third fundamental forms have length
$2\pi$. 

\pg{Manifolds with polyhedral boundary}

Lemma \ref{lm:compact} will follow from the next lemma, which is also of
independent 
interest and might be useful when dealing with manifolds with
a boundary that is locally like a compact or ideal polyhedron. 

\bdf 
Let $(M, g)$ be a hyperbolic 3-dimensional manifold with convex
boundary, and let $E(M)$ be its extension. 
We say that $(M, g)$ is a manifold with polyhedral
boundary if: 
\begin{itemize}
\item for each convex ball $\Omega\subset H^3$ and each isometric
  embedding $\phi:\Omega\rightarrow E(M)$, the
intersection of $M$ 
with $\phi(\Omega)$ is the image by $\phi$ of the intersection with
$\Omega$ of a semi-ideal polyhedron $P\subset H^3$.
\item $\dr M$ contains no closed curve which is a geodesic of $M$.
\end{itemize}
\edf

A semi-ideal polyhedron in $H^3$ is a polyhedron which has vertices
which can be either in hyperbolic space, or on its boundary (ideal
points). For instance, compact polyhedra and ideal polyhedra are
semi-ideal. In this definition, the second condition is necessary
because, otherwise, the boundary of the convex hull of the vertices
could intersect the boundary of the convex core of $M$; this is a case
we want to exclude because our rigidity proof then fails. 

We use the same notion of convergence as for manifolds with hyperideal
boundary, as described near the end of this section,
i.e. Gromov-Hausdorff convergence on compact subsets. 

\blm \label{lm:compact-polyhedral}
Let $(g_n)_{n\in \N}$ be a sequence of metrics on $M$, such that
$(M, g_n)$ are manifolds with polyhedral boundary. Suppose that:
\begin{itemize}
\item for all $n$, $(M, g_n)$ has the same boundary combinatorics, and
  the ideal vertices of $\dr M$ remain the same.
\item the third fundamental forms $\III_n$ of $\dr M$ for the $g_n$
  converge to a limit $\III_\infty$, which is a spherical metric with
  conical singularities.
\item the closed geodesics of $\III_\infty$ which are contractible in
  $M$ have length $L\geq 2\pi$, and $L>2\pi$ except when they bound a
  hemisphere. 
\end{itemize}
Then, after taking a
subsequence, $(M,g_n)_{n\in \N}$ converges to a manifold with polyhedral
boundary $(M, g)$. 
\elm

The proof will be done below, after we show how the proof of lemma
\ref{lm:compact} follows from this lemma. 

\pg{Truncating hyperideal ends}

Let $(M, g)$ be a manifold with hyperideal boundary. By definition it
has a finite number of hyperideal ends, and, if $E$ is one of them, 
there is a totally geodesic plane $P$ in the extension $E(M,g)$ which is
orthogonal to all the faces and edges adjacent to $E$. $P$ is the plane
dual to the hyperideal vertex at $E$. 

One can then {\bf truncate} $(M, g)$ by each of the planes dual to its
hyperideal ends. One obtains in this manner a manifold with polyhedral
boundary, which we call the {\bf truncated manifold} $(M, g_T)$ 
associated to $(M,
g)$. It has two kinds of faces: the {\bf "black"} faces, which are what
remains of the faces of $(M, g)$ after truncation, and the {\bf "red"}
faces, which are where the truncation happened. Each "red" face is
adjacent to "black" faces only.
When a "black" face shares an edge
with a "red" face, the dihedral angle between them is always $\pi/2$.

We will also call {\bf "red" edges} the edges between a black and a red
face, and {\bf "black" edges} the edges between two black faces. So the
"black" edges are what remains of the edges of $(M, g)$ after truncation. 

Consider the universal cover $\Mt$ of $(M, g_T)$ as a convex subset of
$H^3$, and then its boundary $\dr\Mt$. Let $\dr^*\Mt$ be the dual
polyhedron in the de Sitter space, which is invariant under the natural
action of $\pi_1M$ on $S^3_1$. Taking the quotient, we find a compact
polyhedron (in a quotient  of $S^3_1$) for each boundary component of
$M$. This polyhedron has:
\begin{itemize}
\item a "red" vertex for each "red" face of $(M, g_T)$, i.e. for each
  strictly hyperideal vertex of $(M, g)$.
\item a "black" vertex for each "black" face of $(M, g_T)$, i.e. for
  each face of $(M, g)$.
\item a "red" edge of length $\pi/2$ between any red vertex and any adjacent
  black vertex. 
\item a "black" edge between each two adjacent black vertices, of length
  equal to the exterior dihedral angle between the corresponding faces
  of $(M, g)$.
\end{itemize}
Of course, the faces have a metric of constant curvature $1$. Thus we
see that, as for hyperideal polyhedra (see \cite{bao-bonahon,rousset1}),
the dual metric of the boundary of $(M, g_T)$ has a very special
form. Its "black" vertices, and the "black" edges between them, form a graph,
which is combinatorially and 
metrically the dual graph of the boundary of $(M, g)$ (i.e. before
truncation). Each face of this graph 
has a boundary of length at least $2\pi$
(and strictly larger except for faces corresponding to an ideal vertex
of $(M, g)$). 

Moreover, each of those faces with boundary length strictly larger than $2\pi$
contains exactly one "red" vertex, which is connected to each of the 
"black" boundary vertices by a "red" edge of length exactly $\pi/2$. 

Thus we have a pretty simple picture of the third fundamental form of
the boundary of the truncated manifold $(M, g_T)$. It has:
\begin{itemize}
\item one hemisphere for each ideal vertex of $(M, g)$. 
\item one "singular hemisphere", obtained as a quotient of the universal
  cover of a hemisphere minus its "center", for each strictly hyperideal
  vertex of $(M, g)$. The boundary length is then strictly larger than
  $2\pi$. The center of the "singular hemispheres" are the "red"
  vertices. Those "singular hemispheres" have geodesic boundary. 
\item the "singular hemispheres" (and the hemisphere corresponding to
  ideal vertices of $(M, g)$) are glued along the dual graph of the
  boundary of $(M, g)$. 
\end{itemize}
The main points of the discussion above are in the first 3 columns of
the table below; the fourth column is for later reference. 
\medskip

\begin{tabular}{|c|c|c|c|}
\hline
Boundary of $(M,g)$ & Boundary of $(M, g_T)$ & $\III$ of $(M, g_T)$ &
Circle packing limit \\
\hline \hline
Strictly hyperideal vertices & faces & "Central"
vertices & "Red" circles \\ 
\hline
Ideal vertices & Ideal vertices & Hemispheres & Tangency points \\
\hline 
Faces & "Black" faces & "boundary" vertices & "Black" circles
\\ 
\hline
Edges & "Black" edges & Length = exterior angle & Intersection
of $\geq 3$ circles \\
\hline
(vertex, face) & "Red" edges & Length $L = \pi/2$ & Orthogonal
intersection \\
\hline
(vertex, edge) & Vertices & Faces & \\
\hline
\end{tabular}
\medskip

A fundamental remark, made in \cite{rousset1} for hyperideal polyhedra
and hyperideal fuchsian polyhedra, is that the length condition which
appears in lemma \ref{lm:compact-polyhedral}, when applied to $(M,
g_T)$, is equivalent to a statement on the dihedral angles of the
boundary of $(M, g)$. We only outline the proof, since the proof given
in \cite{rousset1} extends to the situation we consider with only
minimal modifications.

\blm \label{lm:translation}
The following statements are equivalent.
\begin{enumerate}
\item Each closed geodesic of $(\dr M, \III)$, which is contractible in
  $M$, has length $L>2\pi$. 
\item The dihedral angles of the boundary of $(M, g)$ satisfy the
  conditions of theorem \ref{tm:angles}.
\end{enumerate}
\elm

\begin{proof}[Outline of the proof]
A key point is that the intersection of a geodesic with the interior of
a singular hemisphere (or a hemisphere) has length exactly
$\pi$. Therefore, a closed geodesic of length $L\leq 2\pi$ can either
remain on the graph made of the "black" edges, or it can enter
only one singular hemisphere. If it enters one singular hemisphere, the
remaining path outside it has length at most $\pi$. 
It follows that (2) $\Rightarrow$ (1).

Conversely, any path made of "black" edges is a geodesic of $\III$,
since, at each vertex, each side is made of at least one singular
hemisphere, so that the total angle on each side is at least $\pi$. So any circuit
in the dual graph of the boundary of $(M, g)$ is a geodesic of
$\III$. Moreover, two vertices in the boundary of a singular hemisphere
$f$ are either at distance less than $\pi$ along the boundary of $f$,
or are joined by a geodesic of length exactly $\pi$
going through the center of $f$. Therefore
the simple paths in the dual graph of $(M,g)$ also correspond to closed
geodesics of $\III$. This shows that condition (1) implies (2).
\end{proof}

We can now prove lemma \ref{lm:compact}, admitting lemma 
\ref{lm:compact-polyhedral}, which we will prove below.

\noindent 
\begin{proof}[Proof of lemma] {\it \ref{lm:compact} from lemma
  \ref{lm:compact-polyhedral}.} ~ 
We consider a sequence $(M, g_n)$ of hyperideal
manifolds. Let $(M, g_{T,n})$ be the truncated manifold associated to
$(M, g_n)$, and let $\III_n$ be its dual metric. 

Consider the graph $G$ on $\dr M$ dual to the combinatorics of 
$\dr (M, g)$, with, for each
edge, a length equal to the limit exterior 
angle $\alpha$. Define a metric $\III_\infty$ on $\dr M$ by gluing in
each face of $G$, a singular hemisphere (there is a unique singular
hemisphere with the right boundary length). Then
$\III_\infty=\lim_{n\rightarrow \infty} \III_n$. 

The hypothesis of lemma \ref{lm:compact}, along with lemma
\ref{lm:translation}, show that $\III_\infty$ satisfies the hypothesis
of lemma \ref{lm:compact-polyhedral}: its closed geodesics which are
contractible in $M$ have length $L\geq 2\pi$, and $L>2\pi$ unless they bound a
hemisphere. Applying lemma
\ref{lm:compact-polyhedral} then shows that the sequence of manifolds
with polyhedral boundary $(M, g_{T,n})$ converges to a limit $(M,
g_T)$. 

Since $\III_\infty$ has a very special form --- it is made by gluing
singular hemispheres as described above --- it is the dual metric of a
manifold with polyhedral boundary which is obtained by truncating a
hyperideal manifold $(M, g)$. So $(M, g_T)$ is the truncated manifold
associated to $(M, g)$.
\end{proof}

\pg{Metrics on the boundary}
We now turn to the proof of lemma \ref{lm:compact-polyhedral}.
We consider a sequence $(g_n)$ of
hyperbolic metrics with polyhedral boundary on $M$, with the same
boundary combinatorics 
$\Sigma$.  

Since, for each $n$, $(M, g_n)$ is a hyperbolic manifold with polyhedral
boundary, the faces are semi-ideal triangles --- each is isometric to
a semi-ideal hyperbolic triangle, i.e. a hyperbolic triangle having
vertices which can be either hyperbolic or ideal points. There are 3
types of edges: those joining two "finite" vertices, which we call
"finite edges", those which connect a finite vertex to an ideal vertex,
which will be called "semi-ideal", and the "ideal" edges connecting two
ideal vertices. 

The metric
induced by $g_n$ on $\dr M$ is however not completely determined in general by
the metrics on those triangles; this is particularly clear when all the
vertices are ideal. To reconstruct the metric, one needs some additional
information, related to the following two definitions. 

\bdf 
Let $T$ be a semi-ideal triangle, let $e$ be an edge of $T$, let $v$ be
the vertex of $T$ opposite to $e$, and let
$\eb$ be the complete hyperbolic geodesic containing $e$. There is
a unique point $c\in \eb$ such that the normal to $\eb$ at $c$ contains
the vertex $v$. We will call $c$ the {\bf projection} of $v$ on $e$.
\edf

\bdf \label{df:shift}
Let $e$ be an edge of $\Sigma$, and let $\eb$ be the complete geodesic
containing $e$. Orient it, and let $T_+$ and $T_-$ be
the triangles on the "positive" and "negative" sides of $e$,
respectively. The {\bf shift} of the metric $h$ at $e$ is the oriented
distance between the projections on $\eb$ of the vertices opposite to $e$
in $T_-$ and $T_+$. 
\edf

Note that this definition does not depend on the orientation chosen on
$\eb$. The shifts along the ideal edges, along with the metrics on the
triangles, are the data necessary to understand the metric induced on
$\dr M$, since they describe how the ideal edges are glued. 
This implies in particular the next proposition, which implicitly uses the
hypothesis of lemma \ref{lm:compact-polyhedral}.

\bprop \label{pr:convergence}
Suppose that:
\begin{itemize}
\item the lengths of all compact edges of the boundary converge. 
\item the shifts of all ideal edges of the boundary converge. 
\item each dihedral angle converges.
\end{itemize}
Then, after taking a subsequence, the sequence of metrics $(g_n)$ converges. 
\eprop

\bpv
Let $h_n$ be the sequence of metrics induced on $\dr M$. When two
triangles share an edge which is either finite or semi-ideal, there is a
unique way of gluing them isometrically along their common boundary. When
they share an ideal edge, the gluing is uniquely determined by the shift
along this edge. Moreover, the lengths of the edges uniquely determine
the metric on the triangles. So the $h_n$ are determined uniquely by the
edge lengths and the shifts along the ideal edges.

So, under the hypothesis of the proposition, $(h_n)$ converges to a
metric $h$. We have also supposed that the dihedral angles of all edges
converge. So, for each connected component $\dr_i M$ of $\dr M$, after we
compose with a sequence of hyperbolic isometries, the
lift to $H^3$ of the universal cover of $\dr_i M$ converges on compact
subsets, as a sequence of convex, polyhedral surfaces in $H^3$. 

We now consider the conformal structure at infinity on each connected
component of $\dr E(M)$. Those conformal structures can be reconstructed from
the induced metric $h_n$ and the dihedral angles, by a procedure known as
"grafting": one should "open" each edge, and glue in a strip of width equal to
the exterior dihedral angle. One should also glue in the hole corresponding to
each vertex the interior of a spherical polygon with edge lengths given by the
dihedral angles of the adjacent edges. The conformal structure of the
resulting metric is the conformal structure at infinity of the corresponding
boundary component of $E(M)$. 

By our hypothesis, this shows that the sequence of conformal structures $c_n$
at infinity of $E(M)$ converges. By the Ahlfors-Bers theorem \cite{ahlfors},
$E(M)$ converges. Thus, after extracting a subsequence, each boundary
component of $M$ converges, as a polyhedral surface in a converging sequence
of hyperbolic ends.  

This means that the universal cover of $M$ for the metrics $g_n$, seen
as convex subsets of $H^3$, converge on compact subsets to a convex
domain, on which $\pi_1M$ acts by isometries. The quotient provides the
limit hyperbolic metric $g$ on $M$.  
\epv

\pg{Consequences of short closed curves}

As already mentioned above, the proof of lemma
\ref{lm:compact-polyhedral} is based on the idea that the existence of
a short simple closed curve implies the appearance of a sequence of geodesics
for $\III$ which are either contractible in $M$, with length converging
to $2\pi$, 
or non-contractible in $M$, with length converging to $0$.
We state here two propositions clarifying both aspects of this
phenomenon.  

\bprop \label{pr:length-angles}
Let $c$ be a simple closed curve in $\dr M$ which is contractible in $M$. Let
$(c_n)_{n\in \N}$ be sequence of curves, each homotopic to $c$, such
that the length of $c_n$ for $h_n$ converges to $0$, and:
\begin{itemize}
\item either $c$ is non-contractible in $\dr \Mt$.
\item or $c_n$ bounds a disk in $\dr \Mt$ which contains a point at distance
  at least $1$ from $c_n$ (for the induced metrics $h_n$). 
\end{itemize}
Then, after extracting a subsequence, there exists a
sequence of closed curves $(\cb_n)_{n\in \N}$,
each homotopic to $c$, converging to a geodesic of length $2\pi$ for
$\III_{\infty}$. 
\eprop

\bpv
Choose $n\in \N$. $c_n$ is contractible in $M$, so it lifts to a closed curve
in $\dr \Mt$, which we also call $c_n$, with $\lim_{n\rightarrow
\infty}L(c_n)=0$. 

The universal cover $\Mt$ of $M$ has convex boundary. If $c$ is
non-contractible in $\dr \Mt$, there exists a
complete geodesic $\gamma$ in $\Mt$, which does not intersect $\dr \Mt$,
but such that $c$ is not contractible in $\Mt\setminus
\gamma$. Otherwise, $c_n$ bounds a disk in $\dr \Mt$ which contains a
point $x_n$ at distance at least $1$ from $c_n$; in this case let
$\gamma$ be a geodesic ray in $\Mt$ starting from $x_n$, such that $c_n$
is not contractible in $\Mt\setminus \gamma$.

Let $p_n$ be a point in $\gamma$ such that the distance between $p_n$ and
$c_n$ is minimal. Let $P_n$ be the hyperbolic plane orthogonal to
$\gamma$ at $p_n$. Then $\cb_n:=\dr \Mt\cap P_n$ is a
curve in $\dr \Mt$, homotopic to $c_n$, which for $n\rightarrow \infty$
is arbitrarily close to $c_n$. Moreover, $\lim_{n\rightarrow \infty}
L(\cb_n)=0$.   

For each $n$, let $\cb_n^*$ be the curve in $S^3_1$ made of points dual to the
support planes of $\dr \Mt$ along $\cb_n$. As $n\rightarrow \infty$, $\cb_n^*$
converges to the dual $\gamma^*$ of $\gamma$ in $S^3_1$, which is geodesic,
and therefore of length $2\pi$ like all geodesics of $S^3_1$. Therefore,
$(\cb_n^*)$ converges to a geodesic of $(\dr \Mt)^*$, so $(\cb_n)$ converges
to a geodesic of length $2\pi$ of $\III_{\infty}$. 
\epv

\bprop \label{pr:degen-closed}
Let $c$ be a closed curve in $\dr M$, which is not
contractible in $M$. Suppose that there exists a family $(c_n)$ of curves
homotopic to $c$, each $c_n$ being geodesic for $h_n$, such that the
length of $c_n$ for $h_n$ converges to $0$. Then the length of $c_n$
for $\III_n$ goes to $0$.
\eprop

\bpv
Fix $n$, and consider the universal cover $\Mt$ of $M$ as a subset of
$H^3$. Each $c_n$ lifts to a geodesic in $(\dr \Mt, h_n)$, with
endpoints in $\dr_\infty H^3$. The result is obtained by applying to
those lifted curves the
following elementary statement of hyperbolic geometry (see
\cite{thurston-notes,epstein-marden}): there exists a constant $C>0$
such that, if 
$\Omega\subset H^3$ is a convex set and $g$ is a complete geodesic in
$\dr \Omega$ with endpoints on $\dr H^3$, then, for each segment $s$ of
$g$ of length $l$, the total bending of $s$ is at most $C(l+1)$.
\epv

\pg{Appearance of short closed curves}

In order to prove lemma \ref{lm:compact-polyhedral}, we now investigates
two situations where the metric degenerates --- when the length of an
edge diverges,
and when a shift goes to infinity --- to show that, in each case, a
short geodesic appears. 

\bprop \label{pr:lengths}
Suppose that there is a finite edge $e$ in $\dr M$ whose length goes
to infinity for the induced metrics $(h_n)$. Then, after taking a subsequence,
there exists a closed curve $c$ in $\dr M$ whose length converges to
$0$, and which either is homotopically non-trivial, or bounds a disk 
$D$ which, for $n$ large enough, contains a non-ideal vertex at distance
at least $1$ from $\dr D$ (for $h_n$).  
\eprop

\bpv
Since the length of $e$ goes to infinity, and since the number of
vertices of the triangulation is bounded, we can find arbitrarily long
segments of $e$ which are far enough from all the vertices. More
precisely, for any $L>0$, after taking a subsequence, we can find a
segment $s$ of $e$ of length $L$ such that $s$ is at distance at least
$1$ from all the vertices. 

The triangulation of $\dr M$ is by a finite (i.e. bounded) number of
triangles; since the area of a hyperbolic triangle is bounded, the area
of $\dr M$ for the metrics $h_n$ is uniformly bounded. Therefore, the
infimum of the injectivity radius of the metrics $h_n$ at the points of
$s$ goes to $0$, and the existence of the requested curve follows.
\epv

\bprop \label{pr:shifts}
Suppose that the sequence of shifts of $(h_n)$ does not converge. Then,
after taking a sub-sequence, there is a closed curve $c_0$ in
$\dr M$ whose length goes to $0$, and which either is homotopically
non-trivial, or which bounds a disk $D$ containing at least two vertices
at distance at least $1$ from $\dr D$.
\eprop

\bpv
Same as for proposition \ref{pr:lengths}, but with a segment $s$ of $e$
which lies between the projections on $e$ of the vertices opposite to
$e$ in the triangles adjacent to $e$.
\epv

\pg{Proof of lemma} {\bf \ref{lm:compact-polyhedral}.}
Suppose that the sequence of metrics $(h_n)$ does not converge. 
Proposition \ref{pr:convergence} shows that either there is a finite
edge whose length goes to infinity, or there exists an ideal edge whose
shift goes to infinity. In the first case, proposition \ref{pr:lengths}
provides a short curve which, thanks to propositions
\ref{pr:length-angles} and \ref{pr:degen-closed}, is clearly seen to be
forbidden by the statement of lemma \ref{lm:compact-polyhedral}. In the
second case the same happens but with proposition \ref{pr:shifts}.

\pg{A more general compactness statement}

To prove theorem \ref{tm:koebe}, we will need in section 9 a compactness
statement more general than lemma \ref{lm:compact-polyhedral}. This is
because theorem \ref{tm:koebe} is based on the limit of circle
configurations when many angles converge to $0$, which, in terms of
truncated hyperideal manifolds, means that some edges converge to ideal
vertices --- a phenomenon which is absent from lemma
\ref{lm:compact-polyhedral}. We give here a more general compactness
statement, allowing the appearance of some geodesics of length $2\pi$
for $\III$. This leads to the appearance of cusps in the induced metric on
the boundary, with parts of the boundary disappearing by going to
infinity. 

\blm \label{lm:+general}
Let $(g_n)_{n\in \N}$ be a sequence of metrics on $M$, such that
$(M, g_n)$ are manifolds with polyhedral boundary. Suppose that:
\begin{itemize}
\item for all $n$, $(M, g_n)$ has the same boundary combinatorics $\Sigma$.
\item the third fundamental forms $\III_n$ of $\dr M$ for the $g_n$
  converge to a limit $\III_\infty$, which is a spherical metric with
  conical singularities.
\item there exists a finite family $D_1, \cdots, D_p$ of disjoint disks in
  $\dr M$ such that, for each $i$, $\dr D_i$ is a geodesic of $\III_\infty$ of
  length $2\pi$.
\item $\dr M\setminus (\cup_{i=1}^pD_i)$ contains no closed geodesic of
  $\III_\infty$ of length
  $L\leq 2\pi$, contractible in $M$, except its boundary components.
\end{itemize}
Then, after taking a
subsequence, $(M,g_n)_{n\in \N}$ converges to manifold with polyhedral
boundary $(M, g)$. The combinatorics of the boundary of $(M, g)$ is
obtained from $\Sigma$ by replacing each disk $D_i$ 
by an ideal vertex. Its third fundamental
form is obtained from $\III_{\infty}$ by replacing each $D_i$ by a hemisphere. 
\elm

As above, the convergence is as defined near the beginning of this
section, i.e. Gromov-Hausdorff convergence of compact subsets. 
The proof will be given below; it follows the proof of lemma
\ref{lm:compact-polyhedral}, the only additional point is that, if
$\III_{\infty}$ has a closed geodesic of length $2\pi$ which does not
bound a hemisphere, then the vertices it contains ``go to infinity'', so
that they are replaced in the limit manifold by only one ideal point. 

\pg{Cusps in the limit boundary}

We need to prove that, whenever a sequence of surfaces has a sequence of
closed geodesics of the third fundamental form with length going to $2\pi$,
there is a corresponding sequence of closed curves with length, for the
induced metrics, going to $0$. This is a converse to proposition
\ref{pr:length-angles}.

\bprop \label{pr:dual}
Let $(\phi_n)_{n\in \N}$ be a sequence of complete, convex
embeddings in $H^3$ of a surface $S$. Let $I_n$ and $\III_n$ be the
induced metric and the third fundamental form of $\phi_n$, and suppose
that $\III_n$ converges to a smooth limit $\III_\infty$. Let $c_n$ be
a sequence of curves in $S$, converging to a limit $c$, which is a
geodesic of $\III_\infty$ of length $2\pi$. Then:
\begin{itemize}
\item the length of $c_n$ for $I_n$ converges to $0$.
\item for any $\epsilon>0$, there exists $N\in \N$ such that, for $n\geq
  N$, the injectivity radius of $(S, I_n)$ is at most $\epsilon$ at all
  points within distance at most $1/\epsilon$ from $c_n$. 
\end{itemize}
\eprop

The proof of this proposition follows ideas which can be found in
\cite{RH,CD,these,shu}, so we only outline here. 
Consider the dual embeddings $(\phi^*_n)_{n\in \N}$ of
$S$ in the de Sitter space $S^3_1$. The images of the $\phi^*_n$ are convex
surfaces, and the induced metrics are the $\III_n$. It is proved in the
references cited above that the closed geodesics of convex surfaces in $S^3_1$
have length $L>2\pi$, and that, if a sequence of simple 
geodesics $(\gamma_n)$ with
$\gamma_n\subset \phi^*_n(S)\subset S^3_1$ has lengths going to
$2\pi$, then (after applying a sequence of isometries of $S^3_1$) its images
converge to a geodesic $\gamma_0$ of $S^3_1$. 

Then, by convexity, the surfaces $\phi_n(S)$ are contained in the interiors of
the "cylinders"
$C_n$ dual to the curves $\gamma_n$. But the $C_n$ converge, as $n\rightarrow
\infty$, to the geodesic $\gamma_0^*$ dual to $\gamma_0$. The estimate on the
injectivity radius follows.

\begin{proof}[Proof of lemma] {\it \ref{lm:+general}.~ }
The proof is based on lemma \ref{lm:compact-polyhedral}, but
uses also proposition \ref{pr:dual}.

By proposition \ref{pr:dual}, for each $i\in \{ 1, \cdots, p\}$, the boundary
of $D_i$ is homologous to a curve whose length goes to $0$ for the induced
metrics. Consider the boundary of the universal cover of $M$, which is a
convex surface in $H^3$.  Proposition \ref{pr:dual}, applied to the boundary
of each of the disks $D_i$, shows that each $D_i$ "goes to infinity", with a
thin tube connecting it to the complement of the $D_i$ in $\dr M$. 

For $n$ large enough, we can (equivariantly) "cut" the thin tubes by a plane
which is almost orthogonal to each of the edges which it intersects, and glue
a small 
polygon (of diameter going to $0$ as $n\rightarrow \infty$). The consequence on
the third fundamental forms $\III_n$ is to replace each disk $D_i$ by a
spherical polygon which is almost a hemisphere. 

Now we can apply lemma \ref{lm:compact-polyhedral} to the manifolds with
polyhedral boundary obtained after this surgery, and we obtain the result.
\end{proof}

\pg{Application to circle packings}

Lemma \ref{lm:+general} leads easily to the compactness statement which
will be useful later on. We will consider a notion of configuration of
circles, on which more details will be given in section 9. 

Given a surface $S$ with a $\C P^1$-structure $c$, a {\bf configuration of
  circles} on $(S,c)$ is a finite set of oriented circles $C_1, \cdots, C_n$
  on $S$ for $c$, such that: 
\begin{itemize}
\item each point of $S$ is in a most 2 closed disks bounded by the $C_i$.
\item for each interstice (i.e. connected component of the complement of the
  disks bounded by the $C_i$) $I$, there exists a circle $C$ orthogonal to
  each of the circles adjacent to $I$. 
\end{itemize}
The second condition is reminiscent of the statement of theorem
\ref{tm:koebe}. It is empty when all interstices are bounded by only 3
circles. 

\blm \label{lm:compact-packings}
Let $(c_n, C_n)$ be a sequence of couples, where, for each $n$, $c_n$ is
a $\C P^1$-structure on $\dr M$ induced by a complete, convex co-compact
hyperbolic metric, and $C_n$ is a circle configuration for
$c_n$. Suppose that the incidence graph of the $C_n$ remains the same,
and that all intersection angles go to $0$. Then, after taking a
subsequence, $(c_n)$ converges to a $\C P^1$-structure $c_\infty$ on $\dr
M$, and $(C_n)$ converges to a circle packing $C_\infty$ for
$c_\infty$. 
\elm

\begin{proof}
By definition of a configuration of circles, for each $n$, we can consider
another family $C'_n$ of circles, with one circle corresponding to each
interstice of $C_n$, and with the circles of $C'_n$ orthogonal to the circles
of $C_n$ when they intersect. 

Let $g_n$ be a complete, convex co-compact hyperbolic metric on $M$
corresponding to $c_n$. For each circle $\sigma$ of $C_n$ or $C'_n$, consider
the oriented totally geodesic plane in $(M,g)$ with boundary at infinity
$\sigma$. Since the disks bounded by the circles of $C_n$ and $C'_n$ cover
$\dr M$, the complement of the half-spaces bounded by those planes is a
compact submanifold of $(M,g)$ with polyhedral boundary. Let $g'_n$ be the
metric on this manifold.

The third fundamental form $\III_n$ 
of $\dr M$ for the metric $g'_n$ can quite easily
be described; $\dr M$ has two kind of faces (corresponding respectively to the
circles of $C_n$ and of $C'_n$) and the faces of one kind intersect the faces
of the other kind orthogonally. The dihedral angles between the faces
corresponding to the circles of $C_n$ are the intersection angles between the
corresponding circles. On the other hand, the circles of $C'_n$ are always
disjoint, although they might become tangent in the limit when $n\rightarrow
\infty$. This completely determines $\III_n$ since the vertices of $\dr M$ are
trivalent.

Let $\sigma, \sigmab$ be two intersecting circles of $C_n$. There are two
circles $\sigma', \sigmab'$ of $C'_n$ which intersect both $\sigma$ and
$\sigmab$ orthogonally. Consider a curve in $(\dr M, \III_n)$ which follows
the edges dual to the edge between the faces bounded by $\sigma,
\sigma', \sigmab, \sigmab'$. Its total length is $4\times\pi/2=2\pi$. As
$n\rightarrow\infty$, the restriction of $\III_n$ to the domain bounded by this curve
converges to the metric on a hemisphere. Moreover, by construction there is no
other geodesic path in $\III_\infty$ of length $2\pi$. 

Now apply lemma \ref{lm:+general}, taking as $D_i$ all the disks bounded by a
curve like the one we have just described; there is one such curve for each
edge of the incidence graph of the $C_n$. Lemma \ref{lm:+general} indicates
that the hyperbolic metrics $g'_n$ converge to the hyperbolic metric
$g'_{\infty}$ on a manifold with polyhedral boundary, with one ideal vertex
for each disk $D_i$, i.e. for each edge of the incidence graph of the $C_n$. 

So all vertices of $(M, g'_{\infty})$ are ideal vertices, and each is adjacent
to 4 faces. Two of those faces correspond to limits of circles of the $C_n$,
and two to limits of circles of the $C'_n$. Taking the boundary of those
faces, we find 4 circles, each tangent to exactly one other, each circle
orthogonal 
to the circles of the other pair. 

So the boundary of the faces of $(M, g'_{\infty})$ corresponding to the
circles of the $C_n$ constitute a circle packing on $\dr M$ for the $\C
P^1$-structure on $\dr M$ defined by the complete hyperbolic metric which is
the extension of $g'_n$. 
\end{proof}


\section{Spaces of polyhedra}

We are concerned in this section with the spaces of angle assignations
and of hyperideal manifolds which appear in theorem \ref{tm:angles}. 
It is basically necessary to show that
those spaces are connected. In practice, however, what we will prove is
a little weaker.

For the spaces of hyperideal manifolds, we will actually prove that, given two
hyperideal manifolds (with the same underlying 
topology) it is possible to connect them by a path in a space of
hyperideal manifolds with a larger number of ideal or hyperideal
vertices.

For the space of angle assignations,
the approach used here, as in \cite{ideal}, is
to remark that, when $M$ has incompressible boundary, the conditions on
the various connected components of the
boundary are "independent", so that it is sufficient to prove this in
the simpler, "fuchsian" situation. In this special case the analog of
theorem \ref{tm:angles} was proved by \cite{rousset1} using other
methods, and the connectedness follows.

\paragraph{Spaces of hyperideal manifolds}

We will later prove that the following spaces of polyhedra are
"weakly connected". First we introduce a handy notation. 

\bdf
We call $\cP$ the set of couples $(p, q)$, where $p=(p_1, \cdots, p_n),
q=(q_1, \cdots, q_n)$ are such that, for all
$1\leq i\leq n$, $p_i\in \N, q_i\in \N$, and $p_i+q_i\geq 1$.
\edf

\bdf \label{df:M}
Let $(p, q)\in \cP$. 
We call $\cM_{p,q}$ the space of hyperideal manifolds diffeomorphic to
$M$, with, for each $1\leq i\leq n$, have $p_i$ ideal vertices and $q_i$
hyperideal ends on $\dr_iM$.
\edf

\paragraph{Hyperideal manifolds and configurations of points and
  circles}
First we remark that hyperideal manifolds are associated to
configurations of points and disks in $\dr M$ with a $\C P^1$-structure,
and we will clarify to what extend the reverse correspondence is also
valid. 

Let $n$ be the number of connected components
of $\dr M$, for each $i\in \{1, \cdots, n\}$, we call $\dr_iM$ the
$i$-th connected component of $\dr M$.

\bdf \label{df:C}
Let $(p, q)\in \cP$.
We call $\cC_{p,q}$ the space of triples $(c, P, Q)$, where $c$ is the
$\C P^1$-structure on $M$ corresponding to a complete, convex co-compact
hyperbolic metric on $M$, $P$ is a finite set of distinct points in $\dr
M$, with $p_i$ points in $\dr_iM$, and $Q$ is a finite set of disjoint
closed disks in $\dr M$ 
for $c$, with $q_i$ disks in $\dr_iM$, such that no point in $P$ lies
in a disk in $Q$. 
\edf

Let $g\in \cM_{p,q}$. The extension $E(M)$ of $(M,g)$ is a complete,
convex co-compact hyperbolic manifold homeomorphic to $M$. Under the
isometric embedding of $(M, g)$ in $E(M)$, each ideal vertex of $M$ goes to
an ideal point in $E(M)$. Moreover, for each hyperideal vertex of $M$,
one can consider the dual plane; its boundary at infinity is a circle in
$\dr_\infty E(M)$. It is clear that the ideal vertices are outside the
disks in $\dr_\infty E(M)$ corresponding to the hyperideal vertices. So
$g$ determines an element of $\cC_{p,q}$, which we call $\Phi_{\cM,
  \cC}(g)$. This defines a map $\PMC:\cM_{p,q}\rightarrow
\cC_{p,q}$. 

Conversely, given $\PMC(g)$, one can reconstruct $g$ by taking the
convex hull in $E(M)$ of the ideal points in $\PMC(g)$ and of the
hyperideal points dual to the disks in $\PMC(g)$, so $\PMC$ is
injective. On the other hand, $\PMC$ is in general not surjective,
because, when one takes the convex hull $N$ of a set of ideal or hyperideal
points in $E(M)$, $N$ might not be a hyperideal manifold. The
reason is that $\dr N$ might have a non-empty intersection with the
convex core of $E(M)$, and thus have parts looking like a typical convex
core of hyperbolic 3-manifold, e.g. with a pleating lamination. Some
arguments showing that this indeed happens can be found in
\cite{ideal}. 

Note that taking the convex hull of the ideal and the hyperideal points
in the definition of an element $\gamma=(c, P, Q)\in \cC_{p,q}$ has an
interpretation in terms of the $\C P^1$-structure $c$, the points in $P$
and the disks in $Q$. Indeed, if $g$ is the hyperbolic metric on $M$
corresponding to $c$, it is not difficult to check that the faces
of the convex hull correspond to the maximal disks in $\dr_\infty M$ such
that:
\begin{itemize}
\item their interiors do not contain any point in $P$.
\item when they intersect a disk in $Q$, the intersection has its two
  angles at most $\pi/2$. 
\end{itemize}

Let $\gamma = (c,P,Q) \in \cC_{p,q}$. 
Let $g$ be the complete, convex co-compact
hyperbolic metric on $M$ corresponding to $c$, and let $\Lambda$ be the
limit set, in $\dr_\infty H^3$, of the action of $\pi_1M$ on $H^3$ with
quotient $M$. $\Mt$, with the lifted
metric, is isometric to $H^3$. Thus $\dr \Mt$ is projectively equivalent
to $S^2$, and $P$ and $Q$ lift to sets $\Pt$ and $\Qt$
of points and disks, respectively, in $S^2\setminus \Lambda$.  

\bdf \label{df:free}
$\gamma$ is {\bf free} if there is no closed disk $D\subset S^2$ such
that:
\begin{itemize}
\item $D\cap \Lambda \neq \emptyset$.
\item $\inte(D)\cap \Lambda = \emptyset$.
\item $\inte(D)$ contains no point of $\Pt$.
\item if $\inte(D)$ has non-empty intersection with a disk $D_1\in \Qt$,
then $D\cap D_1$ has both its angles acute (i.e. at most $\pi/2$). 
\end{itemize}
\edf

\bprop \label{pr:reverse}
Let $\gamma = (c,P,Q) \in \cC_{p,q}$. Let $g$ be the convex co-compact
hyperbolic metric on $M$ corresponding to $c$, and let $N$ be the convex 
hull in $(M,g)$ of the points in $P$ and of the hyperideal points
corresponding to the elements of $Q$. The following statements 
are equivalent.
\begin{enumerate}
\item $\dr N\cap C(M)=\emptyset$.
\item $\dr \Nt\cap \tilde{C(M)}$ contains no complete geodesic of $H^3$. 
\item $\gamma$ is free.
\item $\gamma$ is in the image of $\PMC$.
\end{enumerate}
\eprop

\bpv
{\bf (1) $\Leftrightarrow$ (2)}: clearly $\tilde{C(M)}\subset \Nt$,
$\tilde{C(M)}$ is the convex hull of $\Lambda$, while $\Nt$ is
convex. So the intersection between $\dr\Nt$ and $\dr \tilde{C(M)}$, if
it is not empty, contains a line, i.e. a complete hyperbolic geodesic. 

{\bf (1) $\Rightarrow$ (3)}: suppose that $\gamma$ is not free. Let
$D\subset \dr_\infty H^3$ be a disk as in definition \ref{df:free}.
Let $P$ be the corresponding hyperbolic plane, i.e. the plane in $H^3$
with boundary at infinity the boundary of $D$. Then $P\cap \Nt=\emptyset$. But $\dr P\cap \Lambda
\neq \emptyset$, so that $d(P, \dr\tilde{C(M)})=0$. Therefore, $d(\dr
\Nt, \tilde{C(M)})=0$, so that $d(\dr N, C(M))=0$. Since both $\dr N$
and $C(M)$ are compact, $\dr N\cap C(M)\neq\emptyset$. 

{\bf (3) $\Rightarrow$ (2)}: if $\dr \Nt\cap \tilde{C(M)}$ contains
a complete hyperbolic geodesic, then $\dr \Nt$ contains a complete
geodesic $\gamma_0$ with endpoints in $\Lambda\subset \dr_\infty
H^3$. Then by convexity there is a totally geodesic plane $P\subset H^3$
containing $\gamma_0$ but not intersecting the interior of $\Nt$. The
corresponding disk in $\dr_\infty H^3$ contains in its boundary the
endpoints of $\gamma_0$, and thus $\gamma$ is not free.

{\bf (4) $\Rightarrow$ (2)}: 
suppose that $\dr \Nt\cap \tilde{C(M)}$ contains a complete hyperbolic
geodesic $\gamma_0$. Let $\gamma_1$ be the projection of $\gamma_0$ on
the quotient $C(M)$. If $\gamma_1$ is a closed geodesic, it contradicts
point (2) of definition \ref{df:hyperideal}. Otherwise it contradicts
point (1), since $\gamma_1$ has some accumulation points in $\dr N$
where $\dr N$ is not locally polyhedral. 

{\bf (3) $\Rightarrow$ (4)}: 
suppose that $\gamma$ is free. Let $x\in \dr_\infty H^3\setminus
\Lambda$ be such that its projection in $(\dr_\infty H^3\setminus
\Lambda)/\pi_1M$ is
outside the points of $P$ and the disks in $Q$. Consider the disks in
$\dr_\infty H^3/\pi_1M$ which are maximal among the disks which:
\begin{itemize}
\item have interiors which do not contain any point in $P$.
\item when they meet a disk in $Q$, have an intersection which has
  angles at most $\pi/2$. 
\end{itemize}
The boundary of each of those maximal disks either:
\begin{itemize}
\item contains at least 3 points of $P$.
\item contains at least one point of $P$, and is orthogonal to at least
  one circle in $Q$.
\item is orthogonal to at least 2 circles in $Q$. 
\end{itemize}
Since $\gamma$ is free, there is a finite set of orbits of such maximal
disks under the action of $\pi_1M$. As
a consequence, each point of $\dr N$ has a neighborhood $U$ such that
$U\cap\dr N$ is contained in a finite number of
planes, so that $\dr N$ is locally polyhedral. Moreover, we already know
that $\dr \Nt\cap \tilde{C(M)}$ contains no complete geodesic, so that
$\dr N$ contains no closed geodesic, and thus point (2) of definition
\ref{df:hyperideal} is satisfied. So $N$ is a hyperideal manifold, and
this shows that $\gamma\in \PMC(\cM_{p,q})$.
\epv

\paragraph{A topology on the spaces of hyperideal manifolds}
For each choice of $p,q$, there is a natural topology on $\cC_{p,q}$,
which comes from the topology on the space of hyperbolic, convex co-compact
metrics on $M$, and of points and disks configurations on $\dr M$. 
In addition, given $p,q$, there are certain values of $p', q'$ such that
$\cC_{p',q'}$ can be naturally embedded in the boundary of
$\cC_{p,q}$. This happens if there exists $i_0\in \{ 1, \cdots, n\}$
such that:
\begin{itemize}
\item $p'_{i_0}=p_{i_0}-1\geq 1$, $p'_j=p_j$ for all $j\neq i_0$, and
  $q'_i=q_i$ for all $i\in \{ 1, \cdots, n\}$. This corresponds to two
  ideal vertices "collapsing" to one, or to an ideal vertex "collapsing"
  into a disk (e.g. a hyperideal vertex).
\item $p'_{i_0}=p_{i_0}+1$, $q'_{i_0}=q_{i_0}-1$, and $p'_j=p_j,
  q'_j=q_j$ for all $j\neq i_0$. This corresponds to a disk (e.g. a
  hyperideal vertex) "collapsing" to an ideal point. 
\end{itemize}
We call $\cC$ the union of the $\cC_{p,q}$, for the various values of
$p,q$, with the topology described above. 

Since $\PMC$ is an injective map from each $\cM_{p,q}$ to
$\cC_{p,q}$, the topology on $\cC$ determines a topology on $\cM$, the
union of the $\cM_{p,q}$ for the possible choices of $p,q$.

\paragraph{A weak connectivity property}

The main result of this section concerning the spaces of hyperideal
manifolds is the next lemma. For $p=(p_1, \cdots, p_n)$ and $q=(q_1,
\cdots, q_n)$ given, with $p_i+q_i\geq 1$ for all $i\in \{1, \cdots,
n\}$, we call:
$$ \cM_{\leq p, \leq q} := \cup \{ \cM_{p', q'} ~| ~ p'\leq p
~\mbox{and} ~ q'\leq q\}~, $$
with the topology induced by the topology on $\cC$ described above. 
Here $p'\leq p$ means that, for all $i\in \{1, \cdots, n\}$, $p'_i\leq
p_i$.

\blm \label{lm:weak}
Let $(p, q)\in \cP$. Let $m_0, m_1\in \cM_{p,q}$. There exists $p'\geq
q, q'\geq q$
such that $m_0$ and $m_1$ can be connected by a continuous path in
$\cM_{\leq p', \leq q'}$. 
\elm

The proof uses the following proposition. 

\bprop \label{pr:addition}
Choose $p,q$ and $\gamma=(c, P, Q)\in \cC_{p,q}$. Then:
\begin{enumerate}
\item suppose that
  $\gamma \in \PMC(\cM_{p,q})$. Let $p'\geq p, q'\geq q$, and let
  $\gamma'=(c, P', Q')$, with $P\subset P', Q\subset Q'$. Then $\gamma'
  \in \PMC(\cM_{p',q'})$. 
\item there exists $p'\geq p$ and $\gamma'=(c, P', Q)\in \cC_{p',q}$ such
  that $P\subset P'$ and $\gamma'\in \PMC(\cM_{p',q})$. 
\end{enumerate}
\eprop

\bpv
The first point is a direct consequence of proposition \ref{pr:reverse},
and in particular of its point (3), since adding ideal or hyperideal
points to an element $\gamma\in \cC_{p,q}$ which is free obviously
results in an element which is free.

The second point is also a consequence of point (3) of proposition
\ref{pr:reverse}, since adding enough ideal points to any element
$\gamma\in \cC_{p,q}$ eventually leads to an element which is free, and
thus in the image of $\PMC$.
\epv

\bpv[Proof of lemma \ref{lm:weak}]
Let $\gamma_0:=\PMC(m_0), \gamma_1:=\PMC(m_1)$. Since $\cC_{p,q}$ is
connected, there exists a path $(\gamma_t)_{t\in [0,1]}$ in $\cC_{p,q}$
connecting $\gamma_0$ to $\gamma_1$. Let $\gamma_t=(c_t, P_t, Q_t)$. 

Choose $t_0\in (0,1)$. Point (2) of proposition \ref{pr:addition} shows
that there exists a finite subset $P''_{t_0}$ of $\dr M$ such that
$\gamma'_{t_0}:= (c_{t_0}, P_{t_0}\cup P''_{t_0}, Q_{t_0})$ is in
$\PMC(\cC_{p', q})$, where $p':=p+\#(P''_{t_0})$. Adding some more ideal
vertices if necessary, this remains true for $t$ close enough to $t_0$,
so there exists an open interval $I=(a,b)\ni t_0$ and a family
$(P''_t)_{t\in (a,b)}$ such that, for all $t\in (a,b)$, $\gamma'_t:=(c_t,
P_t\cup P''_t, Q_t)\in \PMC(\cC_{p',q})$.

Doing this for all values of $t_0$ and using the compactness of $[0,1]$,
we find a finite sequence of intervals $I_k=(a_k, b_k), 0\leq k\leq N$,
covering $[0,1]$,
with both $(a_k)$ and $(b_k)$ increasing, and a sequence of families
$(P''_{k,t})_{t\in I_k}$, $0\leq k\leq N$, such that $(c_t, P_t\cup
P''_{t}, Q_t)\in \PMC(\cC_{p'(t), q})$. Since $\gamma_0, \gamma_1\in
\PMC(\cC_{p,q})$, we can suppose moreover that
$P''_{0,t}=P''_{N,t}=\emptyset$. 

Let $p_T:=\sum_{k=1}^{N-1} p_k$. Define a family 
$$ (\Pb_t)_{t\in [0,1]} = (\pb_{1,t}, \cdots, \pb_{p_T, t})_{t\in [0,1]} $$ 
of $p_T$-uples of points of $\dr M$ such that, for all $k\in \{1,
\cdots, N-1\}$, all $t\in I_k$, and all $j\in \{1, \cdots, p_k\}$:
$$ \{ \pb_{(j+\sum_{l=1}^{k-1} p_l), t}~| ~1\leq j\leq p_k\} = P''_{k,
t}~. $$ 

For all $t\in [0,1]$, let: 
$$ \gammab_t:=(c_t, P_t\cup \{\pb_{j,t}~| ~1\leq j\leq p_T\}, Q_t)~. $$
By the first point of proposition \ref{pr:addition}, $\gammab_t\in
\PMC(\cC_{p+p_T, q})$ for all $t\in [0,1]$. Taking the inverse image in
$\cM_{p+p_T, q}$ wields the results. 
\epv

\paragraph{Spaces of fuchsian hyperideal manifolds}
The relationship between the hyperideal metrics on $M$ and the elements
of $\cC_{p,q}$ is simpler when one considers fuchsian manifolds,
i.e. manifolds topologically of the form $S\times \R$, where $S$ is a
surface of genus $g\geq 2$, with a hyperbolic metric which admits an
isometric involution fixing a totally geodesic compact
surface. Restricting one's attention to such manifolds means, in terms
of $\cC$, that $\dr M=S_+\cup S_-$, with both $S_-$ and $S_+$
homeomorphic to $S$, and with $c, P$ and $Q$ invariant under a map
sending $S_-$ to $S_+$ (and conversely). We call $\cM^F_{p,q}$ the space
of fuchsian hyperideal manifolds with $p$ ideal and $q$ hyperideal
vertices (clearly $\cM^F_{p,q}=\emptyset$ unless $p$ and $q$ are even).

\bdf \label{df:Cfuchs}
Let $p,q$ be even numbers, $p=2p', q=2q'$. Let $M$ be diffeomorphic to
$S\times \R$, where $S$ is a closed surface of genus $g\geq
2$. $\cC^F_{p,q}$ is the subset of elements $(c, P, Q)\in \cC_{p,q}$ such
that there exists an involution $i:\dr M\rightarrow \dr M$ which
exchanges the two connected components while leaving invariant $c, P$
and $Q$. 
\edf

The following remark is easy and left to the reader.

\brk
Let $g \in \cM_{p,q}$, then $\PMC(g)\in \cC^F_{p,q}$ 
if and only if $\gamma \in \cM^F_{p,q}$. 
\erk

A more interesting fact is that, whenever one considers an element of
$\cC_{p,q}$ which is fuchsian, then the corresponding convex hull is
always a hyperideal hyperbolic manifold (which of course is fuchsian by
the previous remark). 

\bprop \label{pr:core}
Let $\gamma =(c, P, Q) \in \cC^F_{p,q}$, and let $N$ be the complete,
hyperbolic, convex co-compact metric on $M$ determined by $c$. Let $M$
be the convex hull in $N$ of the ideal points in $P$ and the hyperideal
points in $Q$. Then $M$ is a hyperideal hyperbolic manifold, i.e. $\dr
M\cap C(N)=\emptyset$, where $C(N)$ is the convex hull of $N$. 
\eprop

\bpv
By construction, $N=H^3/\Gamma$, where $\Gamma\subset \isom(H^3)$ is the
image of $\pi_1S$ by a morphism, and $\Gamma$ leaves invariant a totally
geodesic plane $P_0\subset H^3$. The limit set of $\Gamma$ in
$\dr_\infty H^3$ is the boundary at infinity of $P_0$, i.e. a circle in
$S^2$. Let $h_-$ and $h_+$ be the hyperbolic metrics on the two
connected components $D_-$ and $D_+$ of $S^2\setminus \dr_\infty P_0$
for which $\Gamma$ acts isometrically.

Let $D_1$ be a disk in $S^2$ which is tangent to $\dr_\infty P_0$,
suppose for instance that $D_1\subset D_+$. $D_1$ corresponds to a
horodisk in the hyperbolic metric $h_+$, so that, if $P\neq \emptyset$,
$D_1$ contains a point of the lift $\Pt$ of $P$ to $D_+$. 
If $P=\emptyset$ then $Q\neq
\emptyset$, and for the same reason $D_1$ contains a disk in the lift
$\Qt$ of $Q$ to $D_+$. So the proposition follows from proposition
\ref{pr:reverse}. 
\epv

We can sum up the previous two statements in the following lemma. 

\blm \label{lm:CMfuchs}
Let $p,q$ be even integers. Then $\cC^F_{p,q}=\PMC(\cM^F_{p,q})$. 
\elm

\paragraph{Results on fuchsian hyperideal manifolds}

First we define the natural notion of "fuchsian polyhedron".

\bdf 
A {\bf fuchsian subgroup} of $\isom(H^3)$ is a subgroup $\Gamma$ of
$\isom(H^3)$ such that:
\begin{itemize}
\item there exists a totally geodesic 2-plane $P_0$ which is globally
  invariant under all elements of $\Gamma$.
\item $P_0/\Gamma$ is a compact surface of genus $g\geq 2$. 
\end{itemize}
\edf

\bdf
A {\bf fuchsian polyhedron} in $H^3$ is a convex, complete polyhedral
surface $P\subset H^3$ 
(with an infinite number of faces) such that there exists a fuchsian
subgroup $\Gamma$ of $\isom(H^3)$, leaving $P$ globally invariant, and
such that $P/\Gamma$ has a finite number of faces. $P$ is {\bf
  hyperideal} if, for each end $E$ of $P$ of infinite area, there exists
a totally geodesic plane which is orthogonal to all the faces of $P$ at
$E$. 
\edf

Note that hyperideal fuchsian polyhedra are obviously related to 
hyperideal hyperbolic manifolds. More precisely, if $P$ is a hyperideal
fuchsian polyhedron, and $\sigma$ is the symmetry in the plane $P_0$
which is left invariant by the fuchsian subgroup of $\isom(H^3)$
associated to $P$, then $P$ and $\sigma(P)$ bound a convex domain
$\Omega\subset H^3$ whose quotient by $\Gamma$ is a fuchsian hyperideal
hyperbolic manifold.  

Given a surface $S$, a {\bf polyhedral embedding} of $S$ in $H^3$ is a
topological embedding whose image is locally like a polyhedron. 

\bdf
Let $S$ be a compact surface of genus $g\geq 2$. A {\bf hyperideal
  fuchsian embedding} of $S$ is a polyhedral embedding of $\St\setminus
C$ into $H^3$, where $C$ is the union of a finite number of orbits of
the action of $\pi_1S$ on $\St$, whose image is a hyperideal fuchsian
polyhedron.   
\edf

Note that this definition allows the existence of some ideal vertices in
the image of the embedding. 

\btm[M. Rousset \cite{rousset1}] \label{tm:fuchsian}
Let $S$ be a compact surface of genus at least 2, let $\sigma$ be a
cellulation of $S$, and let $w:\sigma_1\rightarrow 
(0,\pi)$ be a map on the set of edges of $\sigma$. There exists a
hyperideal fuchsian realization of $S$, with boundary combinatorics
given by $\sigma$ and exterior dihedral angles given by $w$, if and only
if:
\begin{itemize}
\item the sum of the values of $w$ on each circuit in $\sigma_1$ is
at least $2\pi$.
\item The sum of the values of $w$ on each simple path in $\sigma_1$ is
strictly larger than $\pi$.
\end{itemize}
This hyperideal realization is then unique. A vertex is ideal in this
realization if and only if the sum of the values of $w$ on the adjacent
edges is equal to $2\pi$. 
\etm

The proof given in \cite{rousset1} uses a reduction to the case of
manifolds with a boundary which has only ``compact'' points, for
which a result was obtained in \cite{ideal}. There are however some
subtle technical questions, in particular the infinitesimal rigidity of
those polyhedra and the extension to the limit case where some points
are ideal instead of ``compact'', which are proved in \cite{rousset1}
using methods essentially coming from the work of Pogorelov \cite{Po}.

\paragraph{Spaces of angle assignations}

We now have all the tools necessary to define the spaces of angle
assignations which we will need, and to state the relevant connectedness
properties.

\bdf 
Let $\sigma$ be a cellulation of $\dr M$, and let $P=(P_i)_{1\leq i\leq
n}$, where, for each $i\in \{1, \cdots, n\}$, $P_i$ is a subset of the
set of vertices of $\sigma$ in $\dr_iM$. We call $\cA_{\sigma,P}$
the set of functions $w$ from $\sigma_1$ to $(0,\pi)$ which satisfy the
hypothesis of theorem \ref{tm:angles}, and such that a vertex is in
$P_i$ if and only if the sum of the values of $w$ on the adjacent edges
is $2\pi$. For $(p, q)\in \cP$, we also call $\cA_{p,q}$ the union of
the $\cA_{\sigma, P}$  over all cellulations $\sigma$ with $p_i+q_i$
vertices in $\dr_iM$, and with 
$\card(P_i)=p_i$, and we call $\cA:=\cup_\sigma\cA_{p,q}$. 
\edf

Here $P$ corresponds to ideal vertices in the hyperideal
manifolds which will later turn out to be associated to the elements of
$\cA_{\sigma, P}$. 
There is a natural topology on $\cA$, corresponding to some natural
gluings of the ``cells'' $\cA_{\sigma, P}$:
\begin{itemize}
\item $\cA_{\sigma', P'}$ has an natural embedding in $\dr\cA_{\sigma,
p}$ if $\sigma'$ is obtained from $\sigma$ by collapsing an edge $e$ of
$\sigma$ with both endpoints in some $P_i$, and $P'$ is the same as $P$
except that the two vertices just mentioned are replaced by one. The
angle on $e$ is considered to be $\pi$ in the limit where elements of
$\cA_{\sigma, P}$ converge to $\cA_{\sigma', P'}$. 
\item $\cA_{\sigma', P'}$ has an natural embedding in $\dr\cA_{\sigma,
P}$ if $\sigma'$ is obtained from $\sigma$ by removing an edge $e$. The
angle on $e$ is considered to be $0$ for elements of $\cA_{\sigma', P'}$
considered as elements of $\dr\cA_{\sigma, P}$. 
\end{itemize}
It should be noted that, in each case, the result of the transformation
is indeed in $\dr \cA_{\sigma,P}$. 
We can now state and prove a proposition which will be necessary in the
proof of the main result. 

\bprop \label{pr:connex}
For each $\sigma$ and $P$, $\cA_{\sigma, P}$ is affinely equivalent to
the interior of a polytope in $\R^N$, for some $N$. 
If $M$ has incompressible boundary, then, for all $p, q\in \cP$,
$\cA_{p,q}$ is connected. 
\eprop

\bpv
The affine structure on each $\cA_{\sigma,p}$ comes from the parametrization
by the values of the function $w$; the fact that the $\cA_{\sigma,P}$ are
affinely equivalent to Euclidean polytopes is a direct consequence of
the conditions of theorem \ref{tm:angles}.

Note that the conditions on $w$ on each connected components of $\dr M$
are independent. Since we have supposed that $M$ has incompressible
boundary, the conditions on each connected component of $\dr M$ are the
same as in the fuchsian case. Therefore, to prove that $\cA_{p, q}$ is
connected, it 
is sufficient to prove that it is so in the fuchsian case. 

Now note that, in the fuchsian case, the space of hyperideal polyhedra
with a fixed number of ideal and hyperideal vertices and a fixed genus
is connected. This follows from lemma \ref{lm:CMfuchs} and the
connectedness of $\cC^F_{p,q}$, which one can readily check.
Therefore theorem \ref{tm:fuchsian} implies that the set of dihedral
angles assignations is also connected. 
\epv

\paragraph{Proof of the main theorem }

First we consider a fixed cellulation $\sigma$ of $\dr M$, along with a
subset $P$ of its vertices. Let
$\cM_{\sigma, P}$ be the space of hyperideal hyperbolic manifolds with
boundary combinatorics given by $\sigma$, and let
$\Phi_{\sigma,P}:\cM_{\sigma,P}\rightarrow \cA_{\sigma,P}$ be the map
sending a hyperideal hyperbolic manifold to the set of dihedral angles
of the edges of $\sigma$. If $\cM_{\sigma, P}$ is non-empty, then
$\cM_{\sigma, P}$ and $\cA_{\sigma, P}$ are 
manifolds with boundary of the same dimension, and lemma
\ref{lm:rigidity} shows that $\Phi_{\sigma, P}$ is a local homeomorphism
between them. 

Moreover, lemma \ref{lm:compact} shows that $\Phi_{\sigma, P}$ is
proper, so that $\Phi_{\sigma, P}$ is a covering of $\cA_{\sigma, P}$ by
$\cM_{\sigma, P}$; all the elements of $\cA_{\sigma, P}$ have the same
number $N_{\sigma, P}$ of inverse images, which can be $0$ (if
$\cM_{\sigma, P}=\emptyset$), $1$, or larger. 

By proposition \ref{pr:connex}, each $\cA_{p,q}$ is
connected. Since $M$ has incompressible boundary, lemma
\ref{lm:rigidity} shows that the number $N_{\sigma, 
P}$ remains the same when one moves from one cell of $\cA_{p,q}$ ---
corresponding to a space $\cA_{\sigma, P}$ --- to a neighboring cell ---
corresponding to a $\cA_{\sigma', P'}$. Therefore, the $N_{\sigma, P}$
are equal to a fixed number $N_{p,q}$ depending only on the number of
ideal and 
hyperideal vertices in the triangulation of each connected component of
$\dr M$.

When one goes from $\cA_{p,q}$ to $\cA_{p',q}$, with
$p'_{i_0}=p_{i_0}+1, q'_{i_0}=q_{i_0}-1$, and $p'_i=p_i, q'_i=q_i$ for
$i\neq i_0$, the number $N_{p,q}$ can only increase --- this follows from
lemma \ref{lm:+general}, because a sequence of hyperideal manifolds with a
hyperideal vertex which "becomes ideal" has a limit which is a hyperideal
manifold. 

Thus to prove the main theorem it is sufficient to remark that, by a result of
\cite{ideal}, $N_{p,q}=1$ when $q=0$, i.e. when all $q_i=0$, since this
correspond to ideal hyperbolic manifolds.


\section{Induced metrics}

This section contains the proof of lemma \ref{lm:rig-inf}, and then of
theorems \ref{tm:poly-metrics} and \ref{tm:fuchsian-metrics}. They will
follow from the tools introduced in the previous sections, once we have
given some simple definitions and properties of complete hyperbolic
metrics on triangulated surfaces. Those elements are
generalizations of those in section 9 of \cite{ideal}, where only ideal
triangles where considered. 

In all this section we consider a compact surface $S$, with a
triangulation $\sigma$ by a finite number of triangles, and a subset
$V_i$ of the set of vertices of $\sigma$. $S$ is not necessarily
connected, below it will be $\dr M$, $\sigma$ will be the triangulation
induced on $\dr M$ by a cellulation of $M$, and $V_i$ will be the set of
ideal vertices of a hyperideal metric on $M$. 

Recall from section 2 that we can define a hyperideal triangle as the
interior of a triangle in $\R^2$, with its vertices outside the open
unit disk $D^2$, and all its edges intersecting $D^2$. One can then
consider the metric
coming from the projective model of $H^2$ (maybe extended by the de Sitter
metric outside the disk). The triangle might have one, two or three
ideal vertices, which then sit on the boundary $S^1$ of $D^2$. 

If $T$ is a hyperideal triangle, and if $e$ is an edge of $T$ such that
none of the endpoints of $e$ is an ideal vertex, we can define the {\bf
  length} of $e$ as the distance between the hyperbolic geodesics which
are dual to the vertices of $e$. It is a positive number, see section 2
and in particular proposition \ref{pr:rig-1}.

\bdf 
Let $\cT_0$ be the space of hyperideal triangles. We call $\cN$ be the
space of maps from the set of triangles of $\sigma$ to $\cT_0$ such that, for
each edge $e$ of $\sigma$ with no endpoint in $V_i$, the lengths of $e$
for the hyperideal triangles corresponding to both sides of $e$ are equal.
\edf

When $V_i=\emptyset$, each element of $\cN$ determines a complete
hyperbolic metric on the complement of 
the vertices of $\sigma$ in $S$, such that the area of each end is
infinite. This metric is obtained by gluing adjacent hyperideal triangles in a
way such that the lines dual to their vertices have the same endpoints ---
this gluing condition will be used in all this section. 
But when there are some ideal vertices, however,
some additional care is needed because there might be more than one way to
glue the triangles, and the hyperbolic metrics obtained
might not be complete. 

\bdf \label{df-sh-som}
Let $g\in \cN$ and let $v$ be a vertex of $\sigma$. The {\bf shift} of
$g$ at $v$ is the sum of the shifts of $g$ at the edges containing
$v$. We will say that $g$ is {\bf complete} if its shift is zero at all
the vertices in $V_i$. 
The set of complete structures will be denoted by $\cN_c$.
\edf

The definition of the shift of $g$ at an edge can be found in definition
\ref{df:shift}. 
Of course the notion of completeness defined here is the same as the
usual, topological notion. Indeed if the shift of $g$ at a vertex $v$ is
non-zero, it is possible to use this --- and the fact that ideal
triangles are exponentially thin near their vertices --- to attain $v$ in
a finite time, by ``circling'' around it to take opportunity of the
shift. The converse is not difficult to prove either.

\paragraph{The lengths of the edges of $\sigma$}

Let $g\in \cN_c$. Let $v$ be a vertex of $\sigma$ which is in
$V_i$, and let $t_1, \cdots, t_n$ be the triangles of $\sigma$ adjacent
to $v$, in cyclic order. Then, by definition of $\cN$, $v$ is an ideal
vertex of each of the $t_i$. Consider $t_1$ as a triangle in $H^2$, and
choose a horocycle $H_1$ centered at $v$. Do the same thing in $t_2$,
with a horocycle $H_2$ such that $H_1$ and $H_2$ meet at a point of the
edge which is common to $t_1$ and $t_2$. Repeating this operation, one
finds a sequence $H_1, H_2, \cdots, H_n$ of horocycles in the $t_i$. The
fact that $g$ is complete implies that $H_n$ and $H_1$ meet at the edge
which is common to $t_n$ and $t_1$. 

One can do the same for all the vertices of $\sigma$ which are in $V_i$,
with the additional condition that the horocycles corresponding to
different ideal vertices do not intersect. There is then a well-defined
notion of length for all the edges of $\sigma$ for $g$:
\begin{itemize}
\item the length of an edge with no endpoint in $V_i$ was defined above,
  as the distance between the hyperbolic lines dual to the two vertices
  in one of the triangles. 
\item the distance between a vertex $v_i$ in $V_i$ and a vertex $v_h$
  which is not in $V_i$ is
  the distance, in either of the triangles bounded by the edge, between
  the hyperbolic line dual to $v_h$ and the horocycle attached to
  $v_i$. 
\item the distance between two vertices of $V_i$ is the distance between
  the horocycles attached to each. 
\end{itemize}
Of course the set of lengths of the edges depends on the choices of the
horocycles at the vertices in $V_i$. It does so, however, in a very
simple way, since changing the horocycles at a vertex $v\in V_i$ only
adds the same constant to the lengths of all the edges adjacent to $v$. 

\bdf 
Let $\cL_0$ be the space of maps from the set of edges of $\sigma$ to
$\R$. Let $v_i:=\# V_i$, and let $\cL:=\cL_0/\R^{v_i}$, where an element
$(r_1, \cdots, r_{v_i})\in \R^{v_i}$
acts on $\cL_0$ by adding $r_i$ to the numbers attached to all the edges
adjacent to the $i^{\mbox{th}}$ element of $V_i$. 
\edf

The considerations above show that there is natural map $l$ from $\cN_c$ to
$\cL$, defined by sending a complete hyperbolic metric to the lengths of
its edges (defined up to one additive constant for each ideal vertex).

\bprop \label{pr-bij}
$l$ is a homeomorphism between $\cN_c$ and its image in $\cL$.
\eprop

The proof uses the following elementary property of hyperideal triangles in
$H^2$, which is taken with a small generalization from \cite{ideal}.

\bsl
Consider two hyperideal triangles 
$x_1, x_2, x_3$ and $x_1, x_3, x_4$ with disjoint interior, sharing an
edge $(x_1, x_3)$ which has two ideal vertices at its endpoints. 
For each $i\in \{1,2,3,4\}$, let $h_i$ be a
horocycle centered at $x_i$,if $x_i$ ideal, and let $h_i$ be the
hyperbolic line dual to $x_i$ otherwise. 
For $i\neq j$, let $l_{ij}$ be the
distance between $h_i$ and $h_j$ along the geodesic going from $v_i$ to
$v_j$ --- which is negative if $h_i$ and $h_j$ overlap. Let $\pi_2$ and
$\pi_4$ be the orthogonal projections on $(x_3,x_1)$ of $x_2$ and $x_4$
respectively, and let $\delta$ be the oriented distance between $\pi_2$
and $\pi_4$ on $(x_3,x_1)$. Then:
$$ 2\delta = l_{12}-l_{23}+l_{34}-l_{41}~. $$
\esl 

\begin{figure}[h]
\centerline{\psfig{figure=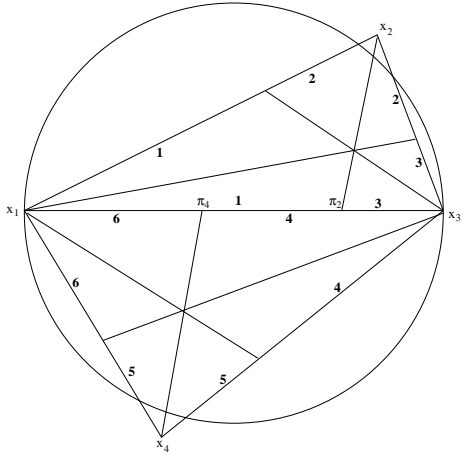,height=6cm}}
\caption{Hyperideal triangles (in the projective model of $H^2$)} 
\label{fig:htriangles}
\end{figure}

\bpv
It follows from figure \ref{fig:htriangles}, where numbers from 1 to 6 are
attached to lengths of segments. Elementary properties of the ideal
triangle show that:
\begin{eqnarray}
l_{12}-l_{23}+l_{34}-l_{41} & = & (1+2)-(2+3)+(4+5)-(5+6)\nonumber 
\\
& = & 1-3+4-6\nonumber \\
& = & (1-6) + (4-3) \nonumber \\
& = & 2\delta~. \nonumber
\end{eqnarray}
\epv

\bpv[Proof of proposition \ref{pr-bij}]
The sub-lemma shows that the lengths of the edges --- in the sense
defined above, with horocycles chosen around the ideal vertices ---
uniquely determine the shifts at the edges of $\sigma$ which have an
ideal vertex at each endpoint. The proposition clearly follows. 
\epv

Proposition \ref{pr-bij} shows that the metrics on $\dr M$ are determined
by the lengths of the edges, as they appear in the Schl{\"a}fli
formula. This means that the strict concavity of the volume translates
as an infinitesimal rigidity statement relative to the metric induced on $\dr
M$. To understand the situations where $\dr M$ is not triangulated, we
will need the next lemma (it is formulated in a more general context,
since it can be useful in different situations). 

Let $P$ be hyperbolic polygon, with vertices $x_1, x_2, \cdots, x_n$
which can be in $H^2$, ideal or hyperideal, and let $h$ be the induced
metric on $P$. Let $H_0$ be a totally
geodesic plane in $H^3$; consider $P$ as a polygon in $H_0$. Let
$\tau_1$ and $\tau_2$ be two triangulations of $P$, i.e. decompositions
into triangles with disjoint interiors and vertices the $x_i$. 

\blm \label{lm:polygone}
Let $\Pd$ be a first-order deformation of the $x_i$ as hyperbolic
(resp. ideal, hyperideal) points in $H^3$. Let $\hbu_1$ and $\hbu_2$ be
the first-order deformations of $h$ obtained through the deformations of
the triangulated surfaces in $H^3$ defined by $\tau_1$ and
$\tau_2$. Then $\hbu_1=\hbu_2$.
\elm

\bpv
It is quite easy to see that the first-order displacements of the $x_i$
orthogonally 
to $H_0$ induce no variation of the lengths of the edges of the
$\tau_i$, and of the shifts at the edges of the $\tau_i$ between two
ideal vertices. The arguments given above thus show that those
orthogonal displacements do not contribute to $\hbu_1$ and $\hbu_2$. 

But the displacements of the $x_i$ tangent to $H_0$ have the same
influence on the deformations $\hbu_1$ and $\hbu_2$, and the lemma follows.
\epv

\blm \label{lm:rig-inf}
Let $(M, g)$ be a hyperideal manifold. It has no first-order
deformation (among hyperideal manifolds with the same ideal vertices)
which does not change the induced metric on $\dr M$.
\elm

\bpv
Theorem \ref{tm:angles} shows that, if $\dr M$ is triangulated, the
deformations of $M$ are parametrized by the deformations of the dihedral
angles at the edges. If $\dr M$ is not triangulated, i.e. if some of its
faces have at least $4$ edges, then choose any triangulation of $\dr M$,
obtained by adding some edges to subdivide the non-triangular
faces. 

Then consider the first-order deformations of the hyperideal
metrics on $M$ which respect this triangulation, but for which $\dr M$
does not necessarily remain convex.
The first-order deformations of the hyperideal metrics 
are parametrized by the first-order deformations of the angles at the
edges of this triangulation, and lemma \ref{lm:polygone} shows that the
first-order variation of the metric on $\dr M$ does not depend on the
triangulation chosen.

By lemma \ref{lm:concave-1}, the volume is a strictly concave function
on the space of hyperideal manifolds, parametrized by the dihedral
angles at the edges of $\dr M$. By the Schl{\"a}fli formula (lemma
\ref{lm:schlafli-2}) this shows that, for each first-order deformation
of the hyperideal metric on $M$, the lengths of the edges have a
non-zero first-order variation. Proposition \ref{pr-bij} then implies
the lemma. 
\epv

An elementary dimension-counting argument then shows that each
first-order deformation of the hyperbolic induced on $\dr M$ can be
obtained from a first-order deformation of the hyperideal metric on
$M$. 

\bcr \label{cr:existence}
Let $(M, g)$ be a hyperideal manifold, and let $h$ be the induced metric
on $\dr M$. Let $\hbu$ be a first-order deformation of $h$ among the
complete hyperbolic metrics on $\dr M$, which have finite area at each
vertex of $M$ which is ideal. There is a
unique first-order deformation of $(M, g)$, among hyperideal manifolds
with the same number of ideal and strictly hyperideal vertices, such that the
variation of the induced metric on $\dr M$ is $\hbu$.
\ecr

\bpv
By theorem \ref{tm:angles}, the deformations of the hyperideal metrics
on $M$, with the same ideal vertices, are parametrized by the dihedral
angles at the edges, under the condition that the sum of the exterior
angles of the edges adjacent to each ideal vertex remains equal to
$2\pi$. The same argument as in the proof of lemma \ref{lm:compact}
shows that those conditions are independent. The dimension of the space
of hyperideal metrics on $M$, with 
the same ideal vertices, is thus $e-v_i$, where $v_i$ is the number of
ideal vertices. 

By proposition \ref{pr-bij}, the complete hyperbolic metrics on $\dr M$
(again with the same ideal vertices) are parametrized by the lengths of
the edges of a triangulation, which are defined up one additive constant
for each ideal vertex. The space of those metrics is therefore also
equal to $e-v_i$. 

Lemma \ref{lm:rig-inf} therefore shows that any first-order deformation
of the metric induced on $\dr M$ is obtained (uniquely) from a
first-order of the hyperideal metric on $M$.
\epv

As a consequence, we can find a result describing the induced metrics on
hyperideal polyhedra. It is a special case of results of \cite{shu}, but
the proof that we obtain here is different. 

\btm \label{tm:poly-metrics}
Let $h$ be a complete hyperbolic metric on $S^2$ minus a finite number
of points (at least $3$). There is a unique hyperideal polyhedron on
$H^3$ whose induced metric is $h$.
\etm

The same result can be obtained in the context of fuchsian hyperideal
manifolds. It was obtained in \cite{ideal} in the special case where all
vertices are ideal, i.e. the metrics which are considered have finite
area. There are also related results for complete, non compact smooth
surfaces, see \cite{rsc}.

\btm \label{tm:fuchsian-metrics}
Let $S$ be a compact surface with non-empty boundary of genus at least
$2$, and let $h$ be a complete hyperbolic metric on $S$ minus a finite
number of points. There is a
unique fuchsian hyperideal manifold $(M, g)$ such that the induced
metric on both connected components of $\dr M$ is $h$.  
\etm

The proofs use some compactness lemmas, which are distinct from lemma
\ref{lm:compact} since they deal with the induced metrics instead of the
dihedral angles. They are both easy to prove. We first deal with theorem
\ref{tm:poly-metrics}.

For each $p,q\in \N$ with $p+q\geq 3$, let:
\begin{itemize}
\item $\cP_{p,q}$ be the space of hyperideal polyhedra with $p$ strictly
  hyperideal and $q$ ideal vertices, up to the isometries of $H^3$. 
\item $\cN_{p,q}$ be the space of
complete hyperbolic metrics on $S^2$ minus $p+q$ points, with $q$ cusps
and $p$ ends of infinite area.
\end{itemize}
 
\blm \label{lm:comp-met-poly}
Let $(P_n)_{n\in \N}$ be a sequence of elements of $\cP_{p,q}$. 
Let $(h_n)_{n\in \N}$ be the
induced metrics, which are elements of $\cN_{p,q}$. 
Suppose that $(h_n)$ converges, as $n\rightarrow
\infty$, to a metric $h\in \cN_{p,q}$. Then, after taking a subsequence,
$(P_n)$ converges to a hyperbolic polyhedron with $p$ strictly
hyperideal and $q$ ideal vertices.
\elm

\bpv
Let $S$ be $S^2$ minus $p+q$ points, let $x\in S$, and let $x_n$ be the
points in the $P_n$ corresponding to $x$. Choose a sequence
$(\phi_n)_{n\in \N}$ of hyperbolic isometries such that $\phi_n(x_n)$
remains equal to a fixed point $y_0\in H^3$. 

Consider the projective model of $H^3$, with $y_0$ at the
center. Suppose that the $(P_n)$ do not converge to a hyperideal
polyhedron $P$, then, after taking a subsequence, at least two of the
vertices, say $v_1$ and $v_2$, converge to the same point in $\R^3$. 

Let $c$ be a closed curve in $S$ going through $x$ and such that $v_1$
and $v_2$ are in different connected components of the complement of
$c$. Then as $n\rightarrow \infty$ the curves in $P_n$ homeomorphic to
$c$ and containing $x$ have to go arbitrarily close to $\dr_\infty H^3$,
so that their lengths have to go to infinity. Therefore $(h_n)$ can not
converge. 
\epv

From here on, we choose a fixed number $g\geq 2$, and call $S_g$ the
closed (compact without boundary) surface of genus $g$. 
For each $p,q\in \N$ with $p+q\geq 1$, let:
\begin{itemize}
\item $\cP^F_{p,q}$ be space of fuchsian hyperideal polyhedra of genus
  $g$ with $p$ strictly
  hyperideal and $q$ ideal vertices, up to the isometries of $H^3$. 
\item $\cN^F_{p,q}$ be the space of 
complete hyperbolic metrics on $S_g$ minus $p+q$ points, with $q$ cusps
and $p$ ends of infinite area.
\end{itemize}

\blm \label{lm:comp-met-fuchs}
Let $(P_n)_{n\in \N}$ be a sequence of elements of $\cP^F_{p,q}$. 
Let $(h_n)_{n\in \N}$ be the
induced metrics, which are elements of $\cN^F_{p,q}$. 
Suppose that $(h_n)$ converges, as $n\rightarrow
\infty$, to a metric $h\in \cN^F_{p,q}$. Then, after taking a subsequence
$(P_n)$ converges to a fuchsian hyperbolic polyhedron with $p$ strictly
hyperideal and $q$ ideal vertices.
\elm

\bpv
The $P_n$ are fuchsian polyhedra with a representation that
fixes a given totally geodesic plane in $H^3$, say $H_0$. So by
definition the representation $\rho_n$ of each $P_n$ acts co-compactly
on $H_0$.  
 
Suppose that
$(\rho_n)$ has no converging subsequence. Then there would exist a
closed geodesic in $S_g$ whose lengths goes to infinity for the
hyperbolic metrics corresponding to the quotient of $H_0$ by
$\rho_n$. Since the projection from $P_n$ to $H_0$ is contracting, the
length of the same curve in $P_n$ would also go to infinity, which would
contradict the hypothesis. So, after taking a subsequence, $(\rho_n)$
converges. 

Suppose that the distance between $P_n$ and $H_0$ does not remain
bounded. Then the orthogonal projection from $P_n$ to $H_0$ would be
contracting by a factor $c_n\rightarrow \infty$. But, after extracting a
subsequence, the length of some closed curve in $H_0/\rho_n$ remains
bounded from below by some positive constant. Therefore the length of
the same curve in $P_n/\rho_n$ goes to infinity, a contradiction. So the
distance between $P_n$ and $H_0$ remains bounded. 

The lemma then follows from the same argument as in the proof of lemma
\ref{lm:comp-met-poly}.
\epv

\bpv[Proof of theorem \ref{tm:poly-metrics}]
Consider the map $F_{p,q}:\cP_{p,q}\rightarrow \cN_{p,q}$ sending a
hyperideal 
polyhedron to its induced metric. By corollary \ref{cr:existence},
$\cP_{p,q}$ and $\cN_{p,q}$ are manifolds of the same dimension, and
$F_{p,q}$ 
is a local homeomorphism. Moreover, by lemma \ref{lm:comp-met-poly},
$F_{p,q}$ is proper, thus it is a covering. 

Both $\cP_{p,q}$ and $\cN_{p,q}$ have retractions to the space of
configurations of $p+q$ points on $S^2$; for $\cN_{p,q}$ it follows from
considerations on the hyperbolic metrics in conformal classes, see
\cite{troyanov}. Both spaces have non-zero Euler characteristic since
$p+q\geq 3$ (this can checked directly by a recursion
argument). Therefore $F_{p,q}$ have degree $\pm 1$, and each element of
$\cP_{p,q}$ has a unique inverse image. 
\epv

\bpv[Proof of theorem \ref{tm:fuchsian-metrics}]
The proof is the same as for theorem \ref{tm:poly-metrics}, except that
the spaces $\cP^F_{p,q}$ and $\cN^F_{p,q}$ now have retractions on the
space of configurations of $p+q$ points in a surface of genus at least
$2$; the argument is otherwise the same. 
\epv

\section{Circle configurations}

\paragraph{Configurations of circles}

It is well known that the Andreev theorem can be formulated in terms of
configurations of circles on $S^2$. Namely, one considers the
decompositions of $S^2$ into a finite number of closed disks, such that each
point of $S^2$ is contained in the interior of at most $2$ of the disks. 

\begin{figure}[h]
\centerline{\psfig{figure=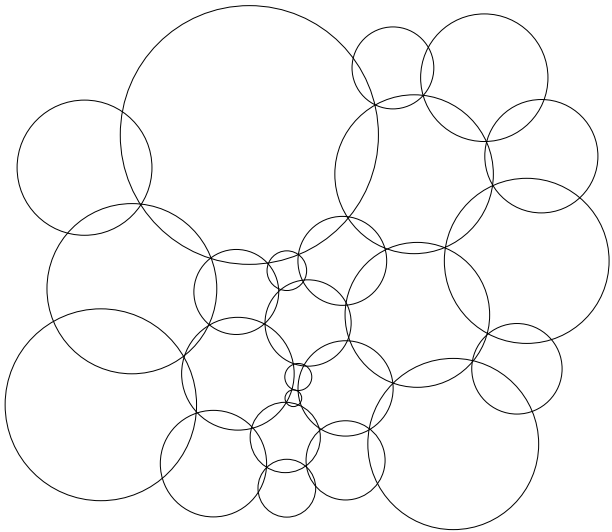,height=7cm}}
\caption{A "classical" configuration of circles}
\end{figure}

The Andreev theorem \cite{Andreev-ideal}, as extended by Rivin
\cite{rivin-annals}, provides an explicit description of the possible
crossing angles between the circles. The points is that one can consider
the ideal polyhedron in $H^3$ with faces the planes bounded by the
circles, and the dihedral angle between two faces is the equal to the
angle between the corresponding circles in $S^2$. 

This can also be done in hyperbolic surfaces, and corresponds to
fuchsian ideal polyhedra, see \cite{thurston-notes,CdeV}. 
It is also possible to give a more general
statement on circle configurations in the boundary of a 3-manifold which
admits a complete, convex co-compact hyperbolic metric \cite{ideal}.

Going from ideal to hyperideal polyhedra leads to another, more general
kind of circle configurations. One now considers two families of circles
on $S^2$:
\begin{itemize}
\item ``red'' circles, which never intersect one another.
\item ``black'' circles, which intersect the red circles orthogonally
and such that the each point of $S^2$ is contained in the interiors of
at most $2$ of the black circles. 
\end{itemize}
We also demand that the sphere is covered by the closed disks bounded by
those circles. 
We will call this setup a ``configuration of red and black
circles''. The red circles correspond to the hyperideal ends of the
polyhedron --- more precisely they are the boundary of the dual
hyperbolic planes --- while the black circles are the boundary of the
planes containing the faces of the polyhedron. 

To each configuration of red and black circles, one can associate a
graph with red vertices and black vertices, corresponding to the red and
to the black circles respectively, with an edge between two vertices if
and only if the corresponding circles intersect. No edge can have red
vertices at both endpoints, since we have specified that red circles do
not intersect. 

To each edge we can
also associate an angle, which we choose to be $\pi$ minus the interior
angle of the intersection between the corresponding two circles. For
edges which have a red vertex as one of their endpoints and a black
vertex at the other, the angle is
$\pi/2$ since we have specified that the black circles intersect the red
circles orthogonally. 

In the picture below, the ``red'' circle are drawn in thicker black
ink, and the ``red'' vertices are drawn with a bigger dot. 

\begin{figure}[h]
\centerline{\psfig{figure=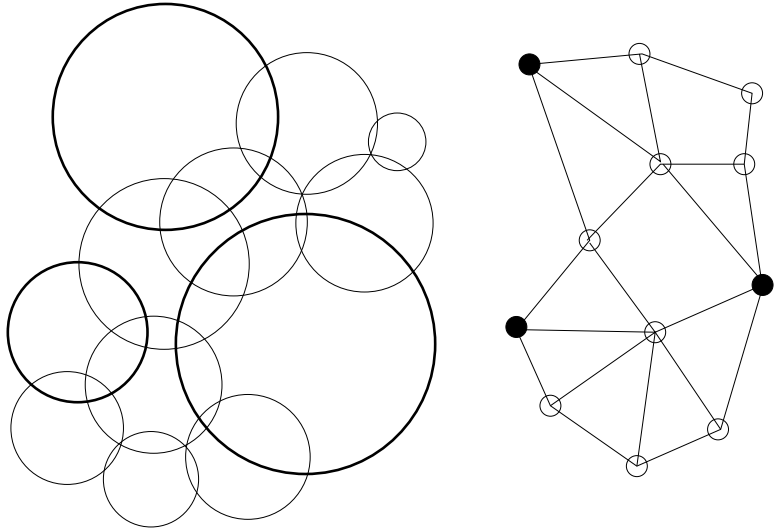,height=7cm}}
\caption{A configuration of red and black circles, and the corresponding
  graph}
\end{figure}

This is not limited to $S^2$, and can also be done in a hyperbolic
surface, or more generally in the boundary $\dr M$ of a 3-manifolds
$M$. To be able to speak about circles, we need to equip $\dr M$ with a
$\C P^1$-structure. We can then reformulate theorem \ref{tm:angles} as
in theorem \ref{tm:circles} below. 

To state it, one must define analogs of the notions of "circuits" and
"simple paths" defined for hyperideal manifolds in the introduction. The
translation is obvious once one remarks that the graph corresponding to
a red and black circle configuration is almost the graph dual to the
cellulation defined from an hyperideal manifold: the black vertices are indeed
associated to the faces of the polyhedron, but the red vertices
correspond to the hyperideal vertices. Therefore, one defines:
\begin{itemize}
\item a {\bf circuit} as a closed curve made of segments of the graph
  having black vertices at both endpoints, and which is contractible in
  $M$. 
\item an {\bf elementary circuit} as a circuit which is made of the
  segments bounding a face.
\item a {\bf simple path} is a closed curve made of segments of the graph, with
  exactly one red vertex, which is contractible in $M$. 
\end{itemize}

\btm \label{tm:circles}
Let $M$ be a 3-manifold with incompressible boundary, which admits a
complete, convex co-compact hyperbolic metric.  
Let $\sigma$ be a cellulation of $\dr M$, with a subset $R$ of the set
of its vertices, and let $w:\sigma_1\rightarrow
(0,\pi)$ be a map on the set of edges of $\sigma$ with no endpoint in
$R$. 
There exists a
complete, convex co-compact hyperbolic metric $g$ on $M$, inducing a $\C
P^1$-structure $c$ on $\dr M$, and a configuration of red and black
circles $C$ for $c$, with combinatorics 
given by $\sigma$, red vertices at the points of $R$, and intersection
angles by $w$, if and only if: 
\begin{itemize}
\item no edge of $\sigma$ has both vertices in $R$.
\item the sum of the values of $w$ on each circuit in $\sigma_1$ is
greater than $2\pi$, and strictly greater if the circuit is
non-elementary.
\item The sum of the values of $w$ on each simple path in $\sigma_1$ is
strictly larger than $\pi$.
\end{itemize}
There is then a unique possible choice of $g$ and of the circle packing
$C$. 
\etm

Given a configuration of red and black circles, one can associate to it
a volume, which is of course defined as the hyperbolic volume of the
corresponding hyperbolic (truncated) hyperideal manifold. From the
Schl{\"a}fli formula (lemma \ref{lm:schlafli-2}), the volume increases when
the exterior dihedral angles of the edges increase. 

Consider a cellulation $\sigma$ of $\dr M$. Add one red vertex for each
face of $\sigma$, with an edge going to each of the vertices of the
face. Let $\sigmab$ be the cellulation --- with red and black vertices ---
obtained, and let $\sigmab_1$ be its 1-skeleton. 
Attach to each edge an angle which is close enough to
$\pi$. Theorem \ref{tm:circles} shows that there is a unique
configuration of red and black circles in $\dr M$ (for a $\C
P^1$-structure coming from a complete convex co-compact hyperbolic
metric on $M$) with the incidence graph $\sigmab_1$ and the prescribed
angles. 

Moreover, if one increases the angles up to $\pi$, the conditions of
theorem \ref{tm:circles} remain satisfied, so that the circle
configurations still exist. In the limit case where all angles are equal
to $\pi$, the black circles are tangent (when they intersect) and the
volume is maximal. It is easy to check that this configuration is
exactly the one given by the extended Koebe theorem \ref{tm:koebe} ---
in the simple case where $M$ is a ball it is the Koebe theorem, as
classically extended to non-triangular interstices.

The nice point about this proof is that it gives directly the two
families of circles in the Koebe theorem; one, the black circles, are
the circles in the circle packing, while the other, the red circles, are
the circles orthogonal to the circles in the packing.

\section{Remarks}

\pg{The strictly hyperideal case}
An interesting remark that came up during a conversation with Francis
Bonahon is that, in the strictly hyperideal case, theorem
\ref{tm:angles} can be proved very simply using the results of
\cite{bonahon-otal} on the convex cores of hyperbolic manifolds, which
are complete when the bending lamination is along closed curves. Indeed,
if $M$ is a strictly hyperideal manifold, one can consider the associated
truncated manifold $M_0$, and glue two copies of $M_0$ along the cuts;
one obtains the convex core of a hyperbolic manifolds, and the convex
cores that can be obtained in this way are characterized by an
elementary symmetry properties. The results of \cite{bonahon-otal} thus
lead to theorem \ref{tm:angles}.

It is interesting to compare those two ways of looking at hyperideal
manifolds with strictly hyperideal vertices, in particular concerning
the infinitesimal rigidity. The fact that both a volume argument and an
argument based on \cite{HK} work might indicate the possibility of
proving infinitesimal rigidity statements other geometric objects, like
hyperbolic cone-manifolds, by methods based on the volume (or some
version of the Hilbert-Einstein functional, which is the same for
hyperbolic metrics). 

\pg{The ideal case}
As already mentioned, the proof given here is slightly different from
the approach used in \cite{ideal}. The main difference is in the
decomposition of the ideal manifolds, which were cut in simplices in
\cite{ideal} and in polyhedra here. It was necessary in \cite{ideal} to
consider finite covers of those manifolds to ensure that the
decomposition was possible, while it is not necessary here. 

\pg{Smooth surfaces}
There is an analog of lemma \ref{lm:angles-poly} for smooth surfaces,
see \cite{rsc}; it
is a result describing the third fundamental form of some complete
surfaces in $H^3$ whose boundary at infinity is a disjoint union of
circles --- this seems to be a natural analog of the condition that, for
each end, all edges going to infinity meet at a hyperideal vertex. There
is also a result on the induced metric. 

It appears that this result should also hold in the setting of
hyperbolic manifolds with boundary. Then one considers complete, constant
curvature $K$ metrics on the boundary minus a finite number of
points. Each such metric should be realized uniquely as the induced metric (if
$K>-1$) or the third fundamental form (for $K\leq 1$) on the boundary,
for a hyperbolic metric on $M$ such that each connected component of
$\dr_{\infty} M$ is a circle for the M{\"o}bius structure of $\dr_{\infty}
H^3$, see \cite{bkc}.

\appendix

\section{The concavity of the volume at a regular simplex}

We give here the {\it maple} code to check proposition
\ref{pr:regulier}, i.e. to make sure that, for at least one regular
hyperideal simplex, the matrix $(\dr \theta_i/\dr l_j)_{i,j}$ is
positive definite. 

\begin{verbatim}
with(linalg);
al:=arccosh((cosh(l0)^2+cosh(l))/sinh(l0)^2);
bl:=arccosh(cosh(l0)*(cosh(l)+1)/(sinh(l)*sinh(l0)));
a0:=arccosh(cosh(l0)*(cosh(l0)+1)/sinh(l0)^2);
t1:=arccos((-cosh(a0)+cosh(bl)^2)/sinh(bl)^2);
t2:=arccos(cosh(bl)*(cosh(a0)-1)/(sinh(bl)*sinh(a0)));
t3:=arccos((cosh(a0)^2-cosh(al))/sinh(a0)^2);
a:=simplify(subs(al=a0,bl=b0,l=l0,diff(t1, l)));
b:=simplify(subs(al=a0,bl=b0,l=l0,diff(t2, l)));
c:=simplify(subs(al=a0,bl=b0,l=l0,diff(t3, l)));
ae:=evalf(subs(l0=1.2, a));
be:=evalf(subs(l0=1.2, b));
ce:=evalf(subs(l0=1.2, c));
M:=matrix([[ae,be,be,ce,be,be],[be,ae,be,be,ce,be],[be,be,ae,be,be,ce],
[ce,be,be,ae,be,be],[be,ce,be,be,ae,be],[be,be,ce,be,be,ae]]);
V:=eigenvals(M);
\end{verbatim}
Feeding this into {\it maple} returns a set of clearly positive
eigenvalues.

\section*{Acknowledgments}

I would like to thank Francis Bonahon and Igor Rivin for some very
interesting conversations related to this work.

\bibliographystyle{alpha}

\end{document}